\newif\ifarxived
\arxivedtrue
\ifx\Arxived\Kaba\else\arxivedtrue\fi
\newif\ifextended
\ifx\Extended\Kaba\else\extendedtrue\fi
\newif\ifprivate
\ifx\Private\Kaba\else\privatetrue\extendedtrue\fi
\newif\ifJapanese
\newif\iftesting
\ifextended\testingtrue\fi
\ifarxived
\extendedfalse
\testingfalse
\privatefalse
\fi
\documentclass[12pt]{article}
\ifarxived
\usepackage[whole]{bxcjkjatype}
\usepackage[bookmarks=true,backref=page,bookmarksnumbered=true, 
  bookmarkstype=toc,
  pdftitle={On Recurrence Axioms},%
  pdfauthor={Sakae Fuchino (渕野 昌) and Toshimichi Usuba (薄葉 季路)},%
  pdfsubject={},%
  pdfkeywords={Generic large cardinal, Laver-generic large cardinal, Maximality Principles}]{hyperref}
\hypersetup{
  colorlinks=true,
  linkcolor=blue,
  filecolor=blue,      
  citecolor=cyan,
  urlcolor=cyan,
  pdftitle={On Recurrence Axioms},%
  pdfauthor={Sakaé Fuchino (渕野 昌) and Toshimichi Usuba (薄葉 季路)},%
  pdfsubject={},%
  pdfkeywords={Generic large cardinal, Laver-generic large cardinal, Maximality Principles}
}
\else
\ifextended
\usepackage[dvipdfmx,backref=page,bookmarks=true,bookmarksnumbered=true,
  bookmarkstype=toc,
]{hyperref}
\hypersetup{
  colorlinks=true,
  linkcolor=blue,
  filecolor=blue,      
  citecolor=cyan,
  urlcolor=cyan,
  pdftitle={On Recurrence Axioms},%
  pdfauthor={Sakaé Fuchino (渕野 昌) and Toshimichi Usuba (薄葉 季路)},%
  pdfsubject={},%
  pdfkeywords={Generic large cardinal, Laver-generic large cardinal, Recurrence Axioms, Maximality Principles}
}
\else
\usepackage[whole]{bxcjkjatype}
\usepackage[bookmarks=true,backref=page,bookmarksnumbered=true, 
  bookmarkstype=toc,
  pdftitle={On Recurrence Axioms},%
  pdfauthor={Sakae Fuchino (渕野 昌) and Toshimichi Usuba (薄葉 季路)},%
  pdfsubject={},%
  pdfkeywords={Generic large cardinal, Laver-generic large cardinal, Maximality Principles}]{hyperref}
\hypersetup{
  colorlinks=true,
  linkcolor=black,
  filecolor=black,      
  citecolor=black,
  urlcolor=black,
  pdftitle={On Recurrence Axioms},%
  pdfauthor={Sakaé Fuchino (渕野 昌) and Toshimichi Usuba (薄葉 季路)},%
  pdfsubject={},%
  pdfkeywords={Generic large cardinal, Laver-generic large cardinal, Recurrence Axioms, Maximality Principles}
}
\fi
\fi 
\ifarxived
\else
\ifextended
\usepackage[dvipdfmx]{pxjahyper}
\else
\fi
\fi

\ifarxived
\usepackage{graphicx, color}
\else
\usepackage[dvipdfmx]{graphicx, color}
\fi
\usepackage{amsmath, amssymb}
\usepackage{bbm}
\usepackage{dsfont}
\usepackage{bbold}
\usepackage{accents} 
\usepackage{marginnote}
\usepackage[many]{tcolorbox}
\usepackage{mathtools}
\usepackage[normalem]{ulem}
\usepackage{setspace}

\definecolor{darkelectricblue}{rgb}{0.33, 0.39, 0.52}
\definecolor{darkgreen}{rgb}{0.31, 0.47, 0.26}
\ifextended
\newcommand{\extendedcolor}{\color{darkelectricblue}}

\newcommand{\privatecolor}{\color{darkgreen}}
\newcommand{\darkred}{\color[rgb]{0.8,0.1,0.1}}

\newcommand{\It}{\it\darkred{}}
\else

\newcommand{\darkred}{}

\newcommand{\It}{\it{}}
\fi
\usepackage{theorem}
\newcommand{\bbd}[1]{{\mathbb{#1}}}
\theorembodyfont{\ifJapanese\rm\else\it\fi}
\newcount\minute	
\newcount\hour		
\newcount\hourMins  
\ifJapanese
\def\today%
{
  \the\year 年\,\zeroPadTwo{\the\month}月\,\zeroPadTwo{\the\day}日%
}
\fi
\def\now%
{
  \minute=\time    
  \hour=\time \divide \hour by 60 
  \hourMins=\hour \multiply\hourMins by 60
  \advance\minute by -\hourMins 
  \zeroPadTwo{\the\hour}:\zeroPadTwo{\the\minute}%
}
\def\zeroPadTwo#1%
{
  \ifnum #1<10 0\fi    
  #1
}
\title{\scalebox{1.1}[1.1]{\bf On Recurrence Axioms}}
\author{\ifarxived%
\protect\scalebox{1}[1.4]{\begin{CJK}{UTF8}{goth}\quad 渕野 昌\qquad\qquad\ \ \ 
    薄葉季路\end{CJK}}\medskip\\
\else
\ifextended\protect\scalebox{1}[1.4]{\ 渕野 昌\qquad\qquad\ \ \ 薄葉 季路}\medskip\\
\fi\fi
  \protect\scalebox{1}[1.4]{\quad\ \ Sakaé Fuchino$^{\ast}$ \quad and\quad Toshimichi Usuba
    $^{\dagger}$}\\{}
  }
\date{}

\ifarxived
\else
\ifextended
\else
\fi
\fi
\renewcommand{\baselinestretch}{1.2}
\renewcommand{\thefootnote}{(\arabic{footnote})\,}
\iftesting
\fi
\iftesting
\newcommand{\Label}[1]{\label{#1}\marginnote{{\color{cyan}%
      \renewcommand{\baselinestretch}{0.4}\tiny 
		  \rlap{#1}}}}

\newcommand{\Labelx}[2]{\label{#1}\marginnote{\color{cyan}{%
      \renewcommand{\baselinestretch}{0.4}\tiny 
		  {}#2\rlap{#1}}}}
\else
\newcommand{\Label}[1]{\label{#1}}

\newcommand{\Labelx}[2]{\label{#1}}
\fi

\def\memo#1{\ifprivate\marginnote{{\darkred\normalsize%
      \renewcommand{\baselinestretch}{0.4}\tiny\mbox{}\vspace{-2.52ex}
      \par\relax%
			#1\par\mbox{}}}\else\fi}%
\def\memox#1{}
\def\imemox#1{}

\newcounter{frml}[section]
\newcounter{frmla}[section]
\def\thefrml{{\arabic{section}.\arabic{frml}}}
\def\thefrmla{{$\aleph$\arabic{section}.\arabic{frmla}}}
\def\frmlabel#1{\refstepcounter{frml}{\def\baka{#1}\ifx\baka\empty\else\label{#1}\fi}%
{\rm({\thefrml})\hfill\hfill\hfill}}
\def\ifrmlabel#1{\refstepcounter{frml}{\def\baka{#1}\ifx\baka\empty\else\label{#1}\fi}%
{\iftesting\darkred\fi\rm({\thefrml})\,:\hspace{0.6em}}}
\def\frmlabela#1{\refstepcounter{frmla}{\def\baka{#1}\ifx\baka\empty\else\label{#1}\fi}%
{\rm({\thefrmla})\hfill\hfill\hfill}}
\def\ifrmlabela#1{\refstepcounter{frmla}{\def\baka{#1}\ifx\baka\empty\else\label{#1}\fi}%
{\iftesting\darkred\fi\rm({\thefrmla})\,:\hspace{0.6em}}}
\def\xitem[#1]{\item[\frmlabel{#1}]\mbox{}%
	\iftesting\marginnote{{\renewcommand{%
				\baselinestretch}{0.6}\color{cyan}\tiny#1}}\fi\ignorespaces}
\def\xitemq[#1]{\item[\frmlabel{#1}]\mbox{}%
	\ignorespaces}
\def\xitemd[#1]#2{\item[(\ref{#1})$#2$\hfill\hfill\hfill]}
\def\xitemA[#1]{\item[\frmlabela{#1}]\mbox{}%
	\iftesting\marginnote{{\renewcommand{%
				\baselinestretch}{0.6}\tiny#1}}\fi\ignorespaces}
\def\xitemx[#1]{\item[]}
\def\xitemsub[#1]#2{\item[\frmlabel{#1}$_{#2}$]\mbox{}%
	\iftesting\marginnote{{\renewcommand{%
				\baselinestretch}{0.6}\tiny#1}}\fi\ignorespaces}
\def\xxitem[#1][#2]{\item[(\ref{#1}{\makebox[1.4ex][c]{#2}})]\mbox{}%
	\iftesting\marginnote{{\renewcommand{%
				\baselinestretch}{0.6}\tiny\{#1\}\{#2\}}}\fi\ignorespaces}
\def\xitemof#1{{\rm({\ref{#1}})}}
\def\Xitem[#1]{\item[{\makebox[7ex][l]{\rm(\ref{#1})}}]\iftesting\marginnote{{\renewcommand{%
				\baselinestretch}{0.6}\tiny#1}}\fi\ignorespaces}
\def\xitemAof#1{{\rm({\ref{#1}})}}

\def\xxitemof#1#2{(\ref{#1}#2)}

\newenvironment{xitemize}{\begin{list}{}{\parsep=0.5\smallskipamount%
			\itemindent=-0.4ex%
			\itemsep=0.5\smallskipamount\leftmargin=4em\labelwidth=3em\labelsep=0.7em}}%
							 {\end{list}}
\iftesting
\def\ixitem[#1]{\ifrmlabel{#1}\marginnote{{\color{cyan}\renewcommand{%
				\baselinestretch}{0.6}\qquad\qquad\tiny\rlap{#1}\mbox{}}}\ignorespaces}
\else
\def\ixitem[#1]{\ifrmlabel{#1}\ignorespaces}
\fi
\iftesting
\def\ixitemx[#1]#2{\ifrmlabel{#1}\marginnote{{\color{cyan}\renewcommand{%
				\baselinestretch}{0.6}{}#2{}\qquad\qquad\tiny\rlap{#1}\mbox{}}}\ignorespaces}
\else
\def\ixitemx[#1]#2{\ifrmlabel{#1}\ignorespaces}
\fi
\iftesting
\def\ixitema[#1]{\ifrmlabela{#1}\marginnote{{\color{cyan}\renewcommand{%
				\baselinestretch}{0.6}\tiny\mbox{}\hfill #1}}\ignorespaces}
\else
\def\ixitema[#1]{\ifrmlabela{#1}\ignorespaces}
\fi

\def\assert#1{\noindent\makebox[4.8ex][r]{\rm(\makebox[2.2ex][c]{#1})}\ \ \ignorespaces}
\def\wassert#1{\assert{#1}}
\def\wassertof#1{\makebox[4.8ex][r]{\rm(\makebox[2.2ex][c]{#1})\ }}%
\def\assertof#1{{\rm(#1)}}%

\def\daimaru#1{\makebox[1em][c]{\mbox{\leavevmode\lower.144ex\hbox{
        \rlap{\hbox to 
          0.76em{\hfil\mbox{}\hfill{}\raisebox{0.054ex}{\scalebox{1.2}{○}}\hfil}}
        \raise0.342ex\hbox to 1em{\hfil{\hspace{0.16em}\footnotesize#1}\hfil}}}}\,}

\newtheorem{Thm}{\ifJapanese{\bf 定理}\else {\bf Theorem}\fi}[section]
\ifextended\tcolorboxenvironment{Thm}{
  colback=blue!4!white,
  boxrule=0pt,
  boxsep=0pt,
  left=8pt,right=8pt,top=2pt,bottom=4pt,
  oversize=4pt,
  sharp corners,
  before skip=\topsep,
  after skip=\topsep,
  breakable
}\fi

\newtheorem{ThmA}{\ifJapanese{\bf 定理\,A\!}\else{\bf Theorem\,A\!}\fi}[section]
\ifextended\tcolorboxenvironment{ThmA}{
  colback=blue!4!white,
  boxrule=0pt,
  boxsep=0pt,
  left=8pt,right=8pt,top=2pt,bottom=4pt,
  oversize=4pt,
  sharp corners,
  before skip=\topsep,
  after skip=\topsep,
  breakable
}\fi

\ifextended\tcolorboxenvironment{Ex}{
  colback=blue!3!white,
  boxrule=0pt,
  boxsep=0pt,
  left=8pt,right=8pt,top=2pt,bottom=4pt,
  oversize=4pt,
  sharp corners,
  before skip=\topsep,
  after skip=\topsep,
  breakable
}\fi
\newtheorem{Prop}[Thm]{\ifJapanese{\bf 命題}\else{\bf Proposition}\fi}
\ifextended\tcolorboxenvironment{Prop}{
  colback=blue!4!white,
  boxrule=0pt,
  boxsep=0pt,
  left=8pt,right=8pt,top=2pt,bottom=4pt,
  oversize=4pt,
  sharp corners,
  before skip=\topsep,
  after skip=\topsep,
  breakable
}\fi

\newtheorem{Lemma}[Thm]{\ifJapanese{\bf 補題}\else{\bf Lemma}\fi}
\ifextended\tcolorboxenvironment{Lemma}{
  colback=blue!4!white,
  boxrule=0pt,
  boxsep=0pt,
  left=8pt,right=8pt,top=2pt,bottom=4pt,
  oversize=4pt,
  sharp corners,
  before skip=\topsep,
  after skip=\topsep,
  breakable
}\fi
\newtheorem{LemmaA}[ThmA]{\ifJapanese{\bf 補題\,A\!}\else{\bf Lemma\,A\!}\fi}
\ifextended\tcolorboxenvironment{LemmaA}{
  colback=blue!4!white,
  boxrule=0pt,
  boxsep=0pt,
  left=8pt,right=8pt,top=2pt,bottom=4pt,
  oversize=4pt,
  sharp corners,
  before skip=\topsep,
  after skip=\topsep,
  breakable
}\fi
\newtheorem{CorA}[ThmA]{\ifJapanese{\bf 系\,A\!}\else{\bf Corollary\,A\!}\fi}
\ifextended\tcolorboxenvironment{CorA}{
  colback=blue!4!white,
  boxrule=0pt,
  boxsep=0pt,
  left=8pt,right=8pt,top=2pt,bottom=4pt,
  oversize=4pt,
  sharp corners,
  before skip=\topsep,
  after skip=\topsep,
  breakable
}\fi
\newtheorem{Cor}[Thm]{\ifJapanese{\bf 系}\else{\bf Corollary}\fi}
\ifextended\tcolorboxenvironment{Cor}{
  colback=blue!4!white,
  boxrule=0pt,
  boxsep=0pt,
  left=8pt,right=8pt,top=2pt,bottom=4pt,
  oversize=4pt,
  sharp corners,
  before skip=\topsep,
  after skip=\topsep,
  breakable
}\fi

\ifextended\tcolorboxenvironment{Remark}{
  colback=blue!4!white,
  boxrule=0pt,
  boxsep=0pt,
  left=8pt,right=8pt,top=2pt,bottom=4pt,
  oversize=4pt,
  sharp corners,
  before skip=\topsep,
  after skip=\topsep,
  breakable
}\fi

\newtheorem{Claim}{{\bf Claim}}[Thm]
\ifextended\tcolorboxenvironment{Claim}{
  colback=pink!5!white,
  boxrule=0pt,
  boxsep=0pt,
  left=8pt,right=8pt,top=2pt,bottom=4pt,
  oversize=4pt,
  sharp corners,
  before skip=\topsep,
  after skip=\topsep,
  breakable
}\fi

\ifextended\tcolorboxenvironment{ClaimA}{
  colback=pink!5!white,
  boxrule=0pt,
  boxsep=0pt,
  left=8pt,right=8pt,top=2pt,bottom=4pt,
  oversize=4pt,
  sharp corners,
  before skip=\topsep,
  after skip=\topsep,
  breakable
}\fi

\newtheorem{Subclaim}{{\bf Subclaim}}[Claim]
\newtheorem{Subsubclaim}{{\bf Subsubclaim}}[Subclaim]

\newcommand{\prf}{\noindent\ifJapanese{\bf 証明．\ }\ignorespaces\else{\bf 
		Proof.\ \ }\ignorespaces\fi}

\newcommand{\prfofClaim}{\raisebox{-.4ex}{\Large $\vdash$\ \ }}

\newcommand{\prfof}[1]{\ifJapanese{\bf #1 の証明．\ \ }%
	\ignorespaces\else{\bf Proof of #1:}\ \ \ignorespaces\fi}
\newcommand{\Thmof}[1]{\ifJapanese{定理\,\ref{#1}}\else{Theorem~\ref{#1}}\fi}

\newcommand{\bfThmof}[1]{\ifJapanese{\bf 定理\,\ref{#1}}\else{\bf Theorem~\ref{#1}}\fi}
\newcommand{\Lemmaof}[1]{\ifJapanese{補題\,\ref{#1}}\else{Lemma~\ref{#1}}\fi}

\newcommand{\LemmaAof}[1]{\ifJapanese{補題\,A\,\ref{#1}}\else{Lemma\,A\,\ref{#1}}\fi}
\newcommand{\CorAof}[1]{\ifJapanese{系\,A\,\ref{#1}}\else{Corollary\,A\,\ref{#1}}\fi}

\newcommand{\Propof}[1]{\ifJapanese{命題\,\ref{#1}}\else{Proposition~\ref{#1}}\fi}

\newcommand{\Corof}[1]{\ifJapanese{系\,\ref{#1}}\else{Corollary~\ref{#1}}\fi}

\newcommand{\Claimof}[1]{{Claim \ref{#1}}}

\newcommand{\Subclaimof}[1]{{Subclaim \ref{#1}}}
\newcommand{\Subsubclaimof}[1]{{Subsubclaim \ref{#1}}}

\newcommand{\sectionof}[1]{\ifJapanese{第\ref{#1}節}\else{Section~\ref{#1}}\fi}

\newcommand{\footnoteof}[1]{Footnote \ref{#1}}
\newcommand{\Thmabove}{{\ifJapanese 定理\else Theorem\fi\ \number\theThm}}

\newcommand{\Corabove}{{\ifJapanese 系\else Corollary\fi\ \number\theThm}}

\newcommand{\Claimabove}{{Claim \number\theClaim}}

\newcommand{\ubecause}[3]{\underbrace{{}#1{}%
  \ifx\bakakaba#2\bakakaba\rule[-0.72ex]{0pt}{1pt}\else\rule[#2]{0pt}{1pt}\fi}_{\mbox{\footnotesize\clap{#3}}}}
\newcommand{\obecause}[3]{\overbrace{{}#1{}%
  \ifx\bakakaba#2\bakakaba\rule[1.62ex]{0pt}{1pt}\else\rule[#2]{0pt}{1pt}\fi}^{\mbox{\footnotesize\clap{#3}}}}
\newsavebox{\qedbox}\sbox{\qedbox}{
{\unitlength=0.05mm \begin{picture}(40,60)
\put(0,0){\framebox(30,44)[cc]{}}
\put(30,-7){\rule{7\unitlength}{44\unitlength}}
\put(10,-7){\rule{27\unitlength}{7\unitlength}}
\end{picture}}}
\newcommand{\qed}{\mbox{}\hfill\usebox{\qedbox}}
\newcommand{\smallqed}%
{\mbox{}\smallskip\hfill\raisebox{-.4ex}{\Large $\dashv$}}
\newcommand{\qedof}[1]%
{\mbox{} \hspace*{\fill}{\usebox{\qedbox}{\tiny~(#1)}}}
\newcommand{\Qedof}[1]%
{\mbox{} \hspace*{\fill}{\usebox{\qedbox}%
{\tiny~(#1~\number\theThm)}}}
\newcommand{\QedAof}[1]%
{\mbox{} \hspace*{\fill}{\usebox{\qedbox}%
{\tiny~(#1~\number\theThmA)}}}
\newcommand{\qedofThm}{\Qedof{\ifJapanese 定理\else Theorem\fi}}

\newcommand{\qedofCor}{\Qedof{\ifJapanese 系\else Corollary\fi}}
\newcommand{\qedofProp}{\Qedof{\ifJapanese 命題\else Proposition\fi}}
\newcommand{\qedofLemma}{\Qedof{\ifJapanese 補題\else Lemma\fi}}
\newcommand{\qedofLemmaA}{\QedAof{\ifJapanese 補題\,A\!\else Lemma\,A\!\fi}}
\newcommand{\qedofCorA}{\QedAof{\ifJapanese 系\,A\!\else Corollary\,A\!\fi}}

\newcommand{\qedskip}{\medskip}

\newcommand{\qedofClaim}%
{\mbox{}\hfill\raisebox{-.4ex}{\Large $\dashv$ }\nolinebreak%
\mbox{\tiny~(Claim~\number\theClaim)}}
\newcommand{\qedofClaimA}%
{\mbox{}\hfill\raisebox{-.4ex}{\Large $\dashv$ }\nolinebreak%
\mbox{\tiny~(Claim~A\,\number\theClaimA)}}
\newcommand{\qedofClaimAof}[1]%
{\mbox{}\hfill\raisebox{-.4ex}{\Large $\dashv$ }\nolinebreak%
\mbox{\tiny~(Claim~A\,\ref{#1})}}
\newcommand{\qedofSubclaim}%
{\mbox{}\hfill\raisebox{-.4ex}{\Large $\dashv$ }\nolinebreak%
\mbox{\tiny~(Subclaim~\number\theSubclaim)}}
\newcommand{\qedofSubsubclaim}%
{\mbox{}\hfill\raisebox{-.4ex}{\Large $\dashv$ }\nolinebreak%
\mbox{\tiny~(Subsubclaim~\number\theSubsubclaim)}}
\newcommand{\cardof}[1]{\mathopen{|\,}#1\mathclose{\,|}}

\newcommand{\Card}{{\it Card\/}}

\newcommand{\setof}[2]{\{#1\,:\,#2\}}
\newcommand{\ssetof}[1]{\{#1\}}

\newcommand{\subseteqand}[1]{\mathrel{\mathop{\subseteq}%
		\limits_{\scriptscriptstyle\hbox to 14pt{$\scriptscriptstyle #1$\hss}}}}

\newcommand{\mapping}[3]{#1:#2\rightarrow #3}

\newcommand{\Elembed}[4]{#1:#2\stackrel{\prec\hspace{0.8ex}}{\rightarrow}_{#4}#3}

\newcommand{\fnsp}[2]{\mbox{}^{{#1}\hspace{-0.02em}}#2}
\newcommand{\imageof}{{}^{\,{\prime}{\prime}}}

\newcommand{\seqof}[2]{\langle#1\,:\,#2\rangle}
\newcommand{\pairof}[1]{\langle#1\rangle}
\newcommand{\psof}[1]{{\mathfrak P}\/(#1)}
\newcommand{\forces}[2]{\,\|\hspace{-.35ex}\mbox{\sf--}_{\,#1\,}%
\mbox{\rm``}\,#2\,\mbox{\rm''}}

\newcommand{\checked}{^{\checkmark}}
\newcommand{\modelof}[1]{\models\!\mbox{\rm``\,}#1\mbox{\rm''}}

\newcommand{\crit}{\mbox{\it crit\/}}

\newcommand{\bbone}{{\mathord{\mathbb{1}}}}

\newcommand{\circleq}{\mathrel{{\leqslant}%
		\hspace{-0.86ex}{\lower-0.53ex\hbox{$\scriptscriptstyle\circ$}}}}

\newcommand{\ol}[1]{\overline{#1}}
\newcommand{\restr}{\restriction}

\newcommand{\cf}{\mathop{cf\/}}

\newcommand{\Col}{{\rm Col}}

\newcommand{\otp}{\mathop{\mbox{\it otp\/}}}
\newcommand{\trcl}{\mathop{\mbox{\it trcl\/}}}

\newcommand{\natnums}{{\bbd{N}}}

\newcommand{\poP}{\bbd{P}}
\newcommand{\poQ}{\bbd{Q}}

\newcommand{\poS}{\bbd{S}}
\newcommand{\On}{{\rm On}}

\newcommand{\genF}{\mathbb{F}}
\newcommand{\genG}{\mathbb{G}}

\newcommand{\utpoP}{\utilde{\mathbb{P}}}

\newcommand{\utpoQ}{\utilde{\mathbb{Q}}}
\newcommand{\utpoR}{\utilde{\mathbb{R}}}

\newcommand{\genH}{\mathbb{H}}

\newcommand{\condp}{\mathbbm{p}}
\newcommand{\condq}{\mathbbm{q}}

\newcommand{\conds}{\mathbbm{s}}

\newcommand{\LT}{{<}\,}
\newcommand{\LE}{{\leq}\,}

\newcommand{\GE}{{\geq}\,}

\newcommand{\ctentenc}{,{}\linebreak[0]\hspace{0.04ex}{{.}{.}{.}\hspace{0.1ex},\,}\linebreak[0]}

\newcommand{\xmbox}[1]{ $\relax{\rm #1}\relax$ }

\newcommand{\continuum}{2^{\aleph_0}}

\newcommand{\calC}{{\mathcal C}}
\newcommand{\calD}{{\mathcal D}}

\newcommand{\calH}{{\mathcal H}}

\newcommand{\calL}{{\mathcal L}}

\newcommand{\calP}{{\mathcal P}}
\newcommand{\calQ}{{\mathcal Q}}

\newcommand{\uta}{\utilde{a}}
\newcommand{\utb}{\utilde{b}}

\newcommand{\utx}{\utilde{x}}

\newcommand{\Lin}{{\calL}_{\in}}

\newcommand{\ZF}{{\sf ZF}}
\newcommand{\ZFC}{{\sf ZFC}}
\newcommand{\CH}{{\sf CH}}

\newcommand{\SCH}{{\sf SCH}}

\newcommand{\MA}{{\sf MA}}
\newcommand{\MM}{{\sf MM}}
\newcommand{\MP}{{\sf MP}}

\newcommand{\RcA}{{\sf RcA}}
\newcommand{\IMH}{{\sf IMH}}
\newcommand{\IGH}{{\sf IGH}}
\newcommand{\LGM}{{\sf LGM}}
\newcommand{\RcAp}{{\sf RcA}{$^+$}}

\newcommand{\PFA}{{\sf PFA}}

\newcommand{\FRP}{{\sf FRP}}

\newcommand{\refl}{{\mathfrak{r}\mathfrak{e}\mathfrak{f}\mathfrak{l}\,}}

\newcommand{\st}{such that}
\newcommand{\wrt}{with respect to}
\newcommand{\Wolog}{Without loss of generality}
\newcommand{\wolog}{without loss of generality}

\newcommand{\tfae}{the following are equivalent}

\newcommand{\po}{poset}
\newcommand{\pos}{posets}
\newcommand{\uniV}{{\sf V}}

\newcommand{\uniW}{{\sf W}}

\newcommand{\Pkl}[2]{\ifx\bakakaba#1\bakakaba\ifx\bakakaba#2\bakakaba{\mathcal 
    P}_\kappa(\lambda)\else{\mathcal P}_\kappa(#2)\fi\else{\mathcal P}_{#1}(#2)\fi}

\newcommand{\utildeT}[1]{%
  \underaccent{{\sim}}{#1}}
\newcommand{\utildeS}[1]{%
	\hbox to 0pt{\smash{$\mathop{\scriptstyle #1}\limits_{%
				\raisebox{0.6ex}[0pt]{$\scriptscriptstyle\sim$}}$}\hss}%
	\relax\phantom{\mathord{{#1}_{\rule[-0.6ex]{0pt}{1pt}}}}}
\newcommand{\utildeSS}[1]{%
	\hbox to 0pt{$\mathop{\scriptscriptstyle #1}%
		\limits_{\scriptscriptstyle\sim}$\hss}%
		\relax\phantom{\underline{#1}}}
\newcommand{\utilde}[1]{%
	\mathchoice{\utildeT{#1}}{\utildeT{#1}}{\utildeS{#1}}{\utildeSS{#1}}}

{\end{minipage}\end{trivlist}}

\newcommand{\baB}{\bbd{B}}

\begin{document}
\maketitle
\renewcommand{\thefootnote}{$\ast$\ }
  \footnotetext{Graduate School of System Informatics, Kobe University \\Rokko-dai 1-1, Nada, Kobe 657-8501 Japan
   \\
    \quad\scalebox{0.95}[1]{\tt fuchino@diamond.kobe-u.ac.jp}}
\renewcommand{\thefootnote}{$\dagger$\ }
  \footnotetext{Faculty of Fundamental Science and Engineering, Waseda University\\
    Okubo 3-4-1, Shinjyuku, Tokyo, 169-8555 Japan 
   \\
    \quad\scalebox{0.95}[1]{\tt usuba@waseda.jp}}
\ifextended
\phantomsection
\addcontentsline{toc}{section}{* On Recurrence Axioms} 
\addcontentsline{toc}{section}{**** by S.Fuchino and T.Usuba}
\fi

\ifextended
\addcontentsline{toc}{section}{Abstract}
\fi
\begin{abstract}
   \memox{The Recurrence Axiom asserts that everything that can happen in the near future happened 
     already at some moment in the past \wrt\ the time-flow along with set 
     forcing extension. The nearness of the future is measured in terms of extent of the class 
     of \pos\ we may consider. The extent of the set of parameters which we may use in the 
     descriptions of "events" in the future also differentiates the variations of the 
     Recurrence Axiom. }
   The Recurrence Axiom for a class  $\calP$ of \pos\ and a set $A$ of parameters is 
   an axiom scheme in the language of \ZFC\ asserting that if a statement 
   with parameters 
   from 
   $A$ is forced by a 
   \po\ in $\calP$, then there is a ground 
   containing the parameters and satisfying the statement
   . 



   The tightly 
   super-$C^{(\infty)}$-$\calP$-Laver generic hyperhuge continuum implies the 
   Recurrence Axiom for $\calP$ and $\calH(\continuum)$. The 
   consistency strength of this assumption can be decided 
   thanks to our main theorems asserting that 
   the minimal ground (bedrock) exists under a tightly $\calP$-generic 
   hyperhuge cardinal $\kappa$, and that $\kappa$ in the bedrock is genuinely hyperhuge, or 
   even super $C^{(\infty)}$ hyperhuge if $\kappa$ is a tightly 
   super-$C^{(\infty)}$-$\calP$-Laver generic hyperhuge definable cardinal. 

   The Laver Generic Maximum (\LGM), one of the strongest combinations of axioms 
   in our context, integrates 
   practically all known set-theoretic principles and axioms in itself, either as its 
   consequences or as theorems holding in (many) grounds of the universe. For instance, 
   double plus version of Martin's Maximum is a consequence of \LGM\ while Cichoń's Maximum is a 
   phenomenon in many grounds of the universe under \LGM. 
\end{abstract}

\iftrue
{
\newpage
\phantomsection
\addcontentsline{toc}{section}{Contents}
\newcommand{\myscalebox}[1]{\scalebox{0.88}[1.06]{#1}}
\begin{quotation}
	\footnotesize
	\noindent
	\centerline{
      \normalsize\tt\quad\ Contents\hspace{6em}\mbox{}}\mbox{}\\
  {\mbox{}\hspace{-1.6em}\tt\makebox[3.4ex][l]{\ref{theintro}.}%
   \hyperref[theintro]{\tt\myscalebox{Introduction}}}\ \ \dotfill\ \ {\pageref{theintro}}\\ 
  {\mbox{}\hspace{-1.6em}\tt\makebox[3.4ex][l]{\ref{intro}.}%
   \hyperref[intro]{\tt\myscalebox{Recurrence Axioms}}}\ \ \dotfill\ \ {\pageref{intro}}\\ 
  {\mbox{}\hspace{-1.6em}\tt\makebox[3.4ex][l]{\ref{Lg-RcA}.}%
   \hyperref[Lg-RcA]{\tt\myscalebox{Recurrence Axioms in restricted forms and the 
       Continuum Problem}}}\ \ \dotfill\ \ {\pageref{Lg-RcA}}\\  
  {\mbox{}\hspace{-1.6em}\tt\makebox[3.4ex][l]{\ref{c-infty}.}%
   \hyperref[c-infty]{\tt\myscalebox{Tightly super-$C^{(\infty)}$-Laver generic 
       large cardinal}}}\ \ \dotfill\ \ {\pageref{c-infty}}\\   
  {\mbox{}\hspace{-1.6em}\tt\makebox[3.4ex][l]{\ref{bedrock}.}%
   \hyperref[bedrock]{\tt\myscalebox{Bedrock of a tightly generic hyperhuge cardinal}}}\ \ 
  \dotfill\ \ {\pageref{bedrock}}\\ 
    {\mbox{}\hspace{-1.6em}\tt\makebox[3.4ex][l]{\ref{bedrock-Lg}.}%
   \hyperref[bedrock-Lg]{\tt\myscalebox{Bedrock and Laver genericity}}}
  \ \ \dotfill\ \ {\pageref{bedrock-Lg}}\\  
  {\mbox{}\hspace{-1.6em}\tt\makebox[3.4ex][l]{\ref{LGM}.}%
   \hyperref[LGM]{\tt\myscalebox{The Laver-Generic Maximum}}}\ \ 
  \dotfill\ \ {\pageref{LGM}}\\  
  {\mbox{}\hspace{-1.6em}\hyperref[ref]{\tt References}}\ \ \dotfill\ \ 
  {\pageref{ref}}\medskip\\ 
\end{quotation}}
\fi
\renewcommand{\thefootnote}{}
\footnotetext{{\it Date:} August 10, 2023
  \qquad {\it Last update:} 
  \today\ (\now\ \ifarxived UTC\else JST\fi)\vspace{-1\smallskipamount}
}
\footnotetext{{\it MSC2020 Mathematical Subject Classification:} 03E45, 03E50, 03E55, 03E57, 03E65
  	\vspace{-1\smallskipamount}}

\footnotetext{
  \hspace{-0.8em}
  $^*$ The research is supported by Kakenhi Grant-in-Aid for 
  Scientific Research (C) 20K03717. 
  \smallskip\\ \indent The authors would like to thank Joan Bagaria, Gunter Fuchs, 
  Hiroshi Sakai, as well as the anonymous referee for helpful remarks and comments. 
} 

\ifextended
\ifprivate
\footnotetext{\hspace{-1em}This is a private extended version of the paper in preparation.
  \par All additional 
  details not contained in the submitted version of the paper are either typeset in 
  {\extendedcolor ``dark electric blue''} in case the text is also included in the 
  (non-private) extended 
  version of the paper, or in {\darkred ``dark red''} if the comment are thought only for the 
  private version. 
  \iffalse
  \sout{The numbering of the assertions is kept identical with the submitted version.}
  Since the changes from the submitted version are now quite extensive, I am not trying to keep the 
  numbering of the theorems and assertions identical with the numbering in the submitted 
  version. \else 
  The numbering of the assertions is kept identical with the submitted version.
  \fi

  The most up-to-date pdf-file of this private extended version is downloadable as:\\
  \ifprivate\href{https://fuchino.ddo.jp/papers/recurrence-axioms-xx.pdf}{\tt 
  https://fuchino.ddo.jp/papers/recurrence-axioms-xx.pdf}\else
  \href{https://fuchino.ddo.jp/papers/recurrence-axioms-x.pdf}{\tt 
  https://fuchino.ddo.jp/papers/recurrence-axioms-x.pdf}\fi 
}
\else
\footnotetext{\extendedcolor This is an extended version of the paper with the 
  title ``On Recurrence Axioms'' in preparation. 
  \par\extendedcolor  All additional 
  details not contained in the submitted version of the paper are either typeset in 
  dark electric blue (the color in which this paragraph is typeset) or put in a separate appendices. 

  The most up-to-date pdf-file of this extended version is downloadable as:\medskip\\
  \qquad\qquad\href{https://fuchino.ddo.jp/papers/recurrence-axioms-x.pdf}{\tt 
  https://fuchino.ddo.jp/papers/recurrence-axioms-x.pdf} \medskip\\
}
\fi
\else
\footnotetext{A pdf file of an updated, and hyperlinked version of this paper 
  possibly with more details and 
  proofs is downloadable as:\quad \href{https://fuchino.ddo.jp/papers/recurrence-axioms-x.pdf}{\tt https://fuchino.ddo.jp/papers/recurrence-axioms-x.pdf}}
\fi

\renewcommand{\thefootnote}{\arabic{footnote})\,}

\section{Introduction}
\Label{theintro}
   The Recurrence Axiom for a class  $\calP$ of \pos\ and a set $A$ of parameters 
   is 
   an axiom scheme in the language of \ZFC\ asserting that if a statement
   $\varphi(\ol{a})$ 
   with parameters $\ol{a}$ in $A$ is forced by a \po\ $\poP\in\calP$, then there is a ground 
   (i.e.\ an inner model from which the universe $\uniV$ is attainable via set forcing)
   containing the parameters and satisfying the statement $\varphi(\ol{a})$
   . 

   Recurrence Axioms can be interpreted as statements about the (eternal?) recurrence in 
   the set generic multiverse in terms of the time-flow along with set forcing extension:  
   everything that can happen in the near future (in form of forcing extension) actually happened 
   already in the past (in a ground). Here, the nearness of the future is 
   measured in terms of extent 
   of the class $\calP$
   of \pos\ we may consider. The extent of the set of parameters which we may use in the 
   descriptions of "events" in the future also differentiates the strength of the axiom. 

   Recurrence Axioms are actually variations of known axioms and principles: they are 
   weakenings of Maximality Principles with corresponding parameters (\Propof{p-intro-0})
   while they can be 
   characterized as the set-generic versions of Sy Friedman's Inner Ground Hypothesis 
   (\Propof{p-intro-2}).  See the end of \sectionof{intro} for discussions about why we want 
   to keep these axioms in spite of this almost identity with other known principles. 

   In \sectionof{Lg-RcA} we show that the tight Laver-generic ultrahugeness implies
   $\Sigma_2$-fragment of Recurrence Axioms (\Thmof{p-Lg-RcA-0}), 
   and 
   $\Sigma_1$-fragments of Recurrence Axioms with strong enough combination of $\calP$ and
   $A$ decide the size of the  
   continuum: in case of $\calP$ being the class of all \pos\ with $A=\calH(\continuum)$ the 
   Continuum Hypothesis (\CH) holds while the Recurrence Axiom for stationary preserving
   $\calP$ with $A=\calH(\kappa_\refl)$ implies $\continuum=\aleph_2$ (\Thmof{p-Lg-RcA-1}). 

   In \sectionof{c-infty} we introduce the notion of tightly 
   super $C^{(\infty)}$-$\calP$-Laver generic hyperhuge cardinal $\kappa$ and show that 
   Recurrence Axiom for $\calP$ and $\calH(\kappa)$ follows from the existence of this 
   generic large cardinal (\Thmof{p-Lg-RcA-5}). The consistency of the 
   existence of this cardinal is strictly between that of hyperhuge cardinal and a 2-huge 
   cardinal (actually the lower bound can be still raised, see \Corof{p-bedrock-8-1}). This 
   follows from the main theorems 
   (\Thmof{p-bedrock-1}, \Thmof{p-bedrock-2}) in \sectionof{bedrock}, asserting that the
   minimal ground (bedrock) under a tightly $\calP$-generic 
   hyperhuge cardinal $\kappa$ exists 
   and that $\kappa$ in the bedrock is genuinely hyperhuge, or even super $C^{(\infty)}$ 
   hyperhuge if $\kappa$ is a tightly  
   super-$C^{(\infty)}$-$\calP$-Laver generic hyperhuge definable cardinal. 

   This result strengthen the theorem on the existence of bedrock by Usuba under a 
   hyperhuge cardinal in \cite{usuba}.

   After examining some of the consequences of the main theorems in \sectionof{bedrock-Lg}, 
   we discuss in \sectionof{LGM} the Laver Generic Maximum (\LGM), one of the strongest axiom 
   available in our context, which integrates 
   practically all known set-theoretic principles and axioms in itself, either as its 
   consequences or  
   as theorems in (many of) the grounds of the universe. So for example, double plus 
   version of Martin's Maximum (\MM$^{++}$)is a consequence of \LGM\ while Cichoń's Maximum is a 
   phenomenon in many grounds of the universe under \LGM. 

   We tried hard to make the present paper as accessible as possible for a wide audience. The 
   terminology and notations used here are either standard or explained fully in the text. 
   For some basic notions nevertheless left unexplained the reader may consult 
   \cite{millennium-book}, 
   \cite{kunen-2011}  and/or \cite{higher-inf}.  \ifarxived The extended version of this paper 
   mentioned in the footnote of the  
   front page includes some more additional details of the proofs and explanations. \fi

\section{Recurrence Axioms}
\Label{intro}
In the following, $\Lin$ denotes the language of \ZFC\ consisting of single binary relation 
symbol $\in$. 

In the language of \ZFC, we always identify a class $\calP$ with the $\Lin$-formula which 
defines the class. Thus, if a class $\calP$ is defined by an $\Lin$-formula $\psi(x)$, 
with $\forces{\poP}{\utx\in\calP\,\dots}$ for a \po\ $\poP$, we simply mean 
$\forces{\poP}{\psi(\utx)\,\cdots}$. We adopt model theoretic convention that (in 
connection with lower case letters) a letter with bar denotes a tuple of objects. Thus, 
$\ol{a}$ means $a_0\ctentenc a_{k-1}$ for some natural number $k$ and write $\ol{a}\in X$ 
for $a_0$\ctentenc $a_{k-1}\in X$. 

We call a class $\calP$ of \pos\ {\It normal\/} if it satisfies
\quad\ixitem[x-intro-0-0] 
$\ssetof{\bbone}\in\calP$, and\\
\ixitem[x-intro-1] 
$\calP$ is closed \wrt\ forcing 
equivalence (i.e.\ if $\poP\in\calP$ and $\poP\sim\poP'$ then $\poP'\in\calP$). 

In the following we assume that all classes $\calP$ of \pos\ we consider 
are normal. In particular when we say that $\calP$ is a class of \pos\ we assume that
$\calP\not=\emptyset$ and it contains the trivial \po.

Some natural classes of \pos\ are not closed under forcing equivalence --- notably the 
class of $\sigma$-closed \pos. For such classes we simply take the closure of the class 
\wrt\ forcing equivalence and replace the class with the closure without mentioning it. 

A (normal) class of \pos\ is {\It iterable} if it also satisfies 
\quad\ixitem[x-intro-2] 
  closed \wrt\ restriction 
  (i.e.\ if $\poP\in\calP$ then $\poP\restr\condp\in\calP$ for any $\condp\in\poP$), and
\quad\ixitem[x-intro-3] 
  for 
  any $\poP\in\calP$ and $\poP$-name $\utpoQ$, $\forces{\poP}{\utpoQ\in\calP}$ implies
  $\poP\ast\utpoQ\in\calP$.

For a 
class $\calP$ of \pos\ and a set $A$ (of parameters), {\It$\calP$-Recurrence Axiom  
  with parameters from $A$} ({\It$(\calP, A)$-\RcA}, for  
short) is  
the following assertion formally expressed as an axiom scheme in $\Lin$:
\begin{xitemize}
\xitem[x-intro-4] For any $\poP\in\calP$,  $\Lin$-formula $\varphi(\overline{x})$ and $\overline{a}\in A$, if
  $\forces{\poP}{\varphi(\overline{a}\checked)}$, then there is a ground 
  $M$ of $\uniV$ \st\ $\overline{a}\in M$ and $M\models\varphi(\overline{a})$. 
\end{xitemize}

Here, an inner model $\uniW_0$ of a universe $\uniW$ is said to be a {\It ground\/} of $\uniW$ if 
there are a  
\po\ $\poP\in \uniW_0$  
and $(\uniW_0,\poP)$-generic $\genG\in\uniW$ \st\  $\uniW_0[\genG]= \uniW$. $\uniW_0$ is a {\It$\calP$-ground\/} of
$\uniW$ if $\poP$  
as above can be taken \st\  $\uniW_0\modelof{\poP\in\calP}$. 

All such grounds are definable. More precisely, the following theorem holds: 
\begin{Thm}\Label{p-intro-1}{\rm(Reitz \cite{reitz}, Fuchs-Hamkins-Reitz \cite{fhr})}
  There is an $\Lin$-formula $\Phi(x,r)$ \st\ the following is provable in \ZFC:\smallskip
  \begin{xitemize}
  \xitem[x-intro-10] for all $r$, $\Phi(\cdot,r):=\setof{x}{\Phi(s,r)}$ is a ground 
    in $\uniV$, and  
  \xitem[x-intro-11] for any ground $M$ (in $\uniV$), there is $r$ 
    \st\ $M=\Phi(\cdot,r)$.\qed 
  \end{xitemize}
\end{Thm}
\memo{\xitemof{x-intro-11} is not formalizable in $\Lin$ in this form! should be further elaborated.}

In the following we use this fact often without explicitly mentioning it. 
As a corollary to \Thmof{p-intro-1}, we immediately see that being a $\calP$-ground for a 
class $\calP$ of \pos\ is a definable property. 

The {\It Strong $\calP$-Recurrence Axiom 
  with parameters from $A$} ({\It $(\calP, A)$-\RcAp}, for short)
holds, if \xitemof{x-intro-4} holds with $M$ which is a 
$\calP$-ground of $\uniV$.

Actually the Strong Recurrence Axiom is equivalent to an already known axiom: 
In the following \Propof{p-intro-0}, we show that, for any (normal) $\calP$,
$(\calP, A)$-\RcAp\  
is equivalent to the Maximality Principle 
$\MP(\calP,A)$ in the notation of \cite{future} (see below for definition --- the 
characterization of $\MP(\calP,A)$ corresponding to this proposition as well as the 
statements corresponding to \Propof{p-intro-2} and 
\Propof{p-intro-3} were also observed by Barton, Caicedo, Fuchs, Hamkins, Reitz, and 
Schindler \cite{5a}). 

For a class $\calP$ of \pos, an $\Lin$-formula $\varphi(\overline{a})$ with parameters
$\overline{a}$ ($\in\uniV$) is said to be a {\It$\calP$-button} if there is $\poP\in\calP$ \st\ 
for any $\poP$-name $\utpoQ$ of \po\ with $\forces{\poP}{\utpoQ\in\calP}$, we have
$\forces{\poP\ast\utpoQ}{\varphi(\overline{a}\checked)}$.

If $\varphi(\overline{a})$ is a $\calP$-button then we call $\poP$ as above a {\It push of 
  the  $\calP$-button $\varphi(\overline{a})$}.

For a class $\calP$ of \pos\ and a set $A$ (of parameters), the {\It Maximality Principle 
  for $\calP$ and $A$} ({\It $\MP(\calP,A)$}, for short)
is the following assertion which is  
formulated as an axiom scheme in $\Lin$: 
\begin{xitemize}
\item[\darkred$\MP(\calP,A)$: ] For any $\Lin$-formula $\varphi(\overline{x})$ and
  $\overline{a}\in A$, if $\varphi(\overline{a})$ is a $\calP$-button then
  $\varphi(\overline{a})$ holds.
\end{xitemize}


\begin{Prop}
  \Label{p-intro-0} Suppose that $\calP$  is a class of \pos\ and $A$ a set 
  (of parameters).\smallskip\\ \wassert{1} $(\calP, A)$-\RcAp\ is equivalent to
  $\MP(\calP,A)$.\smallskip 

  \wassert{2} $(\calP, A)$-\RcA\ is 
  equivalent to the following assertion:  
\begin{xitemize}
\xitem[x-intro-5-0] For any $\Lin$-formula $\varphi(\overline{x})$ and
  $\overline{a}\in A$, if $\varphi(\overline{a})$ is a $\calP$-button then
  $\varphi(\overline{a})$ holds in a ground of\/ $\uniV$. 
\end{xitemize}
\end{Prop}
\prf \assertof{1}: Suppose first that $(\calP, A)$-\RcAp\ holds. We show that 
$\MP(\calP,A)$ holds. Assume that $\poP\in\calP$ is a push of the $\calP$-button
$\varphi(\overline{a})$. Let $\varphi'(\overline{x})$ be the formula expressing
\begin{xitemize}
\xitem[x-intro-6] for any $\poQ\in\calP$,
  $\forces{\poQ}{\varphi(\overline{x}\checked)}$ holds. 
\end{xitemize}
Then we have $\forces{\poP}{\varphi'(\overline{a}\checked)}$. By 
$(\calP, A)$-\RcAp, there is a $\calP$-ground $M$ of $\uniV$ \st\ $\overline{a}\in M$ and
$M\models\varphi'(\overline{a})$ holds. By the definition \xitemof{x-intro-6} 
of $\varphi'$, it follows that $\uniV\models\varphi(\overline{a})$ holds. 

Now suppose that $\MP(\calP,A)$ holds and $\poP\in\calP$ is \st\
$\forces{\poP}{\varphi(\overline{a}\checked)}$ 
for $\overline{a}\in A$.

Let $\varphi''$ be a formula claiming that
\begin{xitemize}
\xitem[x-intro-7] there is a $\calP$-ground $N$ \st\ $\overline{x}\in N$ and 
  $N\models\varphi(\overline{x})$. 
\end{xitemize}
Then $\varphi''(\overline{a})$ is a $\calP$-button and $\poP$ is its push.

By $\MP(\calP,A)$, $\varphi''(\overline{a})$ holds in $\uniV$ and hence there is 
a $\calP$-ground $M$ 
of $\uniV$ \st\ $\overline{a}\in M$ and $M\models\varphi(\overline{a})$. This shows that 
$(\calP, A)$-\RcAp\ holds. \smallskip

\assertof{2}: can be proved similarly to \assertof{1}. {\ifextended\extendedcolor 
  Suppose first that $(\calP, A)$-\RcA\ holds. We show that 
  \xitemof{x-intro-5-0} holds. Assume that $\poP\in\calP$ is a push of the $\calP$-button
  $\varphi(\overline{a})$. Let $\varphi'(\overline{x})$ be the formula expressing
  \begin{xitemize}
  \xitem[x-intro-8] for any $\poQ\in\calP$,
    $\forces{\poQ}{\varphi(\overline{x}^{\scalebox{0.8}[0.5]{✓}})}$ holds. 
  \end{xitemize}
  Then we have $\forces{\poP}{\varphi'(\overline{a}^{\scalebox{0.8}[0.5]{✓}})}$. By 
  $(\calP, A)$-\RcA, there is a ground $M$ of $\uniV$ \st\ $\overline{a}\in M$ and
  $M\models\varphi'(\overline{a})$ holds. 
  Since $\calP\ni\ssetof{\bbone}$ (see the remark in the paragraph after 
  \xitemof{x-intro-0-0}), it follows that $M\models\varphi(\overline{a})$.  

  Now suppose that \xitemof{x-intro-5-0} holds and $\poP\in\calP$ is \st\
  $\forces{\poP}{\varphi(\overline{a}^{\scalebox{0.8}[0.5]{✓}})}$ 
  for $\overline{a}\in A$.
  Let $\varphi''$ be a formula asserting that
  \begin{xitemize}
  \xitem[x-intro-9] there is a $\calP$-ground $N$ \st\ $\overline{x}\in N$ and 
    $N\models\varphi(\overline{x})$. 
  \end{xitemize}
  Then $\varphi''(\overline{a})$ is a $\calP$-button and $\poP$ is its push.
  Thus, By \xitemof{x-intro-5-0}, $\varphi''(\overline{a})$ holds in a ground $M$ of $\uniV$ 
  with $\overline{a}\in M$. By the definition 
  \xitemof{x-intro-9} of $\varphi''$, there is a $\calP$-ground $N$ of $M$ \st\
  $\overline{a}\in N$ and $N\models\varphi(\overline{a})$. Since $N$ is also a ground of
  $\uniV$, this shows that 
  $(\calP, A)$-\RcA\ holds. 
  \fi}
\qedofProp\qedskip

Recurrence Axioms are also related to the Inner Model Hypothesis introduced by Sy 
Friedman in \cite{friedman-sy}. {\It The Inner Model Hypothesis} ({\It\IMH}) is the 
following assertion formulated in the language of second-order set theory (e.g.\ 
in the context of von Neumann-Bernays-Gödel set theory): 
\begin{xitemize}
\item[{\It\IMH} :] For any statement $\varphi$ without parameters, if $\varphi$ holds in an 
  inner model of an inner extension of $\uniV$ then $\varphi$ holds in an inner model of
  $\uniV$. 
\end{xitemize}
Here we say a (not necessarily first-order definable) transitive class $M$ an {\It inner 
  model}\/  
of $\uniV$ if $M$ is a model of \ZF\ and $\On^M=\On^\uniV$. In the perspective from 
such $M$, we call $\uniV$ an {\It inner extension} of $M$. 

We shall call a set-forcing version of this principle {\It Inner Ground Hypothesis} 
({\It\IGH}):  

For a (definable normal) class $\calP$ of \pos\ and a set $A$ (of parameters), 
\begin{xitemize}
\item[{\It$\IGH(\calP,A)$} :] For any $\Lin$-formula $\varphi=\varphi(\overline{x})$ and
  $\overline{a}\in A$, if $\poP\in\calP$ forces ``there is a ground $M$ with
  $\overline{a}\in M$ satisfying $\varphi(\overline{a})$'', then there is a ground $\uniW$ of
  $\uniV$ \st\ $\overline{a}\in\uniW$ and $\uniW\models\varphi(\overline{a})$. 
\end{xitemize}

\begin{Prop}
  \Label{p-intro-2} For a class $\calP$ of \pos\ and a set $A$ (of 
  parameters), $(\calP,A)$-\RcA\ holds if and only if\/ $\IGH(\calP,A)$ holds.
\end{Prop}
\prf Suppose that $(\calP,A)$-\RcA\ holds. Let $\varphi=\varphi(\overline{x})$ be 
an $\Lin$-formula, $\overline{a}\in A$, and $\poP\in\calP$  be \st\
$\forces{\poP}{\varphi(\overline{a}\checked)\xmbox{ holds in a ground}}$. 

Let $\varphi'(\overline{x})$ be the $\Lin$-formula asserting that $\varphi(\overline{x})$ 
holds in a ground. Then $\forces{\poP}{\varphi'(\overline{a}\checked)}$. 
By $(\calP,A)$-\RcA, it follows that there is a ground $\uniW$ of $\uniV$ \st\
$\uniW\models\varphi'(\overline{a}\checked)$. Since a ground of a ground is a ground, we 
conclude that there is a ground $\uniW_0$ of $\uniV$ \st\ $\overline{a}\in M_0$ and
$\uniW_0\models\varphi(\overline{a})$. This shows that $\IGH(\calP,A)$ holds.

Suppose now that $\IGH(\calP,A)$ holds. Assume that
$\forces{\poP}{\varphi(\overline{a}\checked)}$ for an $\Lin$-formula
$\varphi=\varphi(\overline{x})$, $\overline{a}\in A$, and $\poP\in\calP$. Then
$\forces{\poP}{\varphi(\overline{a}\checked)\xmbox{ holds in a }
\calP\xmbox{-ground (of the universe)}}$ since $\forces{\poP}{\ssetof{\bbone}\in\calP}$. 
Thus, by $\IGH(\calP,A)$, there is a ground $\uniW$ of $\uniV$ \st\
$\uniW\models\varphi(\overline{a})$. 
\qedofProp
\qedskip

$(\poP,A)$-\RcAp\ ($\Leftrightarrow$\ $\MP(\calP,A)$ by \Propof{p-intro-0},\,\assertof{1}) can be also 
characterized in terms of a strengthening of Inner Ground Hypothesis:
For a (definable) class $\calP$ of \pos\ and a set $A$ (of parameters), 
\begin{xitemize}
\item[{\It$\IGH^+(\calP,A)$} :] For any $\Lin$-formula $\varphi=\varphi(\overline{x})$ and
  $\overline{a}\in A$ if $\poP\in\calP$ forces ``there is a $\calP$-ground $M$ with
  $\overline{a}\in M$ satisfying $\varphi(\overline{a})$'', then there is a $\calP$-ground
  $\uniW$ of 
  $\uniV$ \st\ $\overline{a}\in\uniW$ and $\uniW\models\varphi(\overline{a})$. 
\end{xitemize}

The following proposition can be proved similarly to \Propof{p-intro-2}.

\begin{Prop}
  \Label{p-intro-3} For a class $\calP$ of \pos\ and a set $A$ (of parameters),
  $(\calP,A)$-\RcAp\ holds if and only if $\IGH^+(\calP,A)$ holds. 
  \ifextended\else\qed\fi
\end{Prop}
{\ifextended\extendedcolor
\prf Suppose that $(\calP,A)$-\RcAp\ holds and assume that $\varphi=\varphi(\overline{x})$ 
is an $\Lin$-formula, $\overline{a}\in A$, and $\poP\in\calP$ is \st\ 
\begin{xitemize}
\item[] $\forces{\poP}{\varphi(\overline{a}\checked)\mbox{ holds in a }\calP\mbox{-ground }M
  \mbox{ with }\overline{a}\in M}$
\end{xitemize}
Let $\varphi'(\overline{a})$ be the formula expressing ``$\varphi(\overline{x})$ holds in a
$\calP$-ground $M$ with $\overline{a}\in M$''. Then $\poP$ is a push of the $\poP$-button 
$\varphi'(\overline{a})$. Thus, by \Propof{p-intro-0},\,\assertof{1},
$\varphi'(\overline{a})$ holds in $\uniV$. By definition of $\varphi'$, there is 
a $\calP$-ground $\uniW$ of $\uniV$ \st\ $\overline{a}\in\uniW$ and
$\uniW\models\varphi(\overline{a})$. 
This shows that
$\IGH^+(\calP,A)$ holds. 

Suppose now that $\IGH^+(\calP,A)$ holds, and assume that $\varphi=\varphi(\overline{x})$ 
is an $\Lin$-formula, $\overline{a}\in A$ 
and $\forces{\poP}{\varphi(\overline{a}\checked)}$ then  
(since $\ssetof{\bbone}\in\calP$)
$\forces{\poP}{\varphi(\overline{a}\checked)
  \xmbox{ holds in a }\calP\xmbox{-ground }M\xmbox{ with }\overline{a}\in M}$.
By $\IGH^+(\calP,A)$, it follows that there is  $\calP$-ground $\uniW$ of 
$\uniV$ \st\  
$\overline{a}\in\uniW$ and $\uniW\models\varphi(\overline{a})$. 
This shows 
that $(\calP,A)$-\RcAp\ holds. 
\qedofProp
\qedskip\fi}

In spite of these almost identity with other known principles, we want to keep the 
Recurrence Axioms as axioms on their own. One of the reasons is that the formulation of the 
principle emphasizes the downward absoluteness feature of the principle. Another is 
that we have the following 
monotonicity for (non plus versions) these axioms which is not valid with the Maximality 
Principles (i.e. the plus versions of the Recurrence Axioms).  

\begin{Lemma}{\rm (Monotonicity of Recurrence Axioms)}\Label{p-5}
  For classes of \pos\ $\calP$, $\calP'$ and sets $A$, $A'$ of parameters, 
  if $\calP\subseteq \calP'$ and $A\subseteq A'$, then we have 
  \begin{xitemize}
  \item[] $(\calP', A')$-\RcA\ \ $\Rightarrow$\ \ $(\calP,A)$-\RcA. \qed
  \end{xitemize}
\end{Lemma}

If we decide that the Recurrence Axioms provide desirable extensions of the axioms of 
\ZFC, then we should try to take the maximal instance of these axioms. (i.e.\ the one with 
maximal strength among the 
instances consistent with \ZFC) 
By \Lemmaof{p-5}, this means we should try to take the instance of Recurrence Axioms  
with the maximal $\calP$ and $A$ (\wrt\ inclusion) among the consistent ones. 

\Thmof{p-Lg-RcA-1} \memox{\large !!!} in the next section suggests that the following 
\xitemof{x-intro-12} and \xitemof{x-intro-13} are candidates of such maximal instances. 

Let $\kappa_\refl$ denote the cardinal number
$\max\ssetof{\continuum,\aleph_2}$. $\kappa_\refl$ appears often as the reflection point of 
a strong structural reflection principle (see \cite{future}).  
\begin{itemize}
\xitem[x-intro-12] $\ZFC$ $+$ $(\calP,\calH(\kappa_\refl))$-RcA\ for the class 
  $\calP$ of all 
  stationary preserving \pos.
\end{itemize}\begin{itemize}
\xitem[x-intro-13] $\ZFC$ $+$ $(\calQ,\calH(\continuum))$-RcA\ for the class
  $\calQ$ of all \pos.  
\end{itemize}

The consistency of \xitemof{x-intro-13} follows from the consistency of \ZFC\ $+$ ``there are 
stationarily many inaccessible cardinals'' (\cite{hamkins}). The consistency of 
\xitemof{x-intro-12} follows from  \Lemmaof{p-Lg-RcA-2}, 
\Thmof{p-Lg-RcA-4},\,\assertof{$B'$}, and   
\Thmof{p-Lg-RcA-5}.\memox{\large!!!}

The maximality of \xitemof{x-intro-12} and \xitemof{x-intro-13} follows from 
\Thmof{p-Lg-RcA-1},\,\assertof{2'} and \assertof{5'} respectively.\memox{\large!!!}

By \Thmof{p-Lg-RcA-1},\,\assertof{4} and \assertof{5}, \xitemof{x-intro-12} implies
$\continuum=\aleph_2$, and \xitemof{x-intro-13} implies \CH. In  
particular, these two extensions of \ZFC\ are not compatible. However, as we are going to 
discuss in \sectionof{LGM}, we can combine the plus version of \xitemof{x-intro-12} with 
the following weakening of \xitemof{x-intro-13}: 

\begin{itemize}
\xitem[x-intro-14]  $\ZFC$ $+$ $(\calP,\calH(\kappa_\refl))$-\RcAp\ $+$ 
$(\calQ,{\calH(\omega_1)}^{\ol{\uniW}})$-\RcAp\ where $\calP$ is the class of all 
  semi-proper \pos, $\calQ$ the class of all \pos, and $\ol{\uniW}$ the 
  bedrock\footnotemark\ which is also assumed here to exist.
\end{itemize}
\footnotetext{For the definition of the bedrock see \sectionof{bedrock}.}

Note that $\continuum=\aleph_2$ follows from \xitemof{x-intro-14}. 
In \sectionof{LGM}, we give an axiom in terms of existence of strong variants of Laver 
generic large cardinals from which \xitemof{x-intro-14} follows. 

\section{Recurrence Axioms in restricted forms and the 
       Continuum Problem} 
\Label{Lg-RcA}

We consider the following restricted forms of Recurrence Axiom: For an iterable class $\calP$ 
of \pos, a set $A$ (of parameters), and a set $\Gamma$ of $\Lin$-
formulas, {\It$\calP$-Recurrence Axiom for formulas in $\Gamma$ with parameters from $A$} 
({\It$(\calP,A)_\Gamma$-\RcA}, for short) is the following assertion expressed as an axiom 
scheme in $\Lin$:
\begin{xitemize}
\xitem[x-Lg-RcA-a-0] For any $\varphi(\overline{x})\in\Gamma$ and $\overline{a}\in A$, if
  $\forces{\poP}{\varphi(\overline{a}\checked)}$, then there is a ground 
  $\uniW$ of $\uniV$ \st\ $\overline{a}\in \uniW$ and $\uniW\models\varphi(\overline{a})$. 
\end{xitemize}
{\It$(\calP,A)_\Gamma$-\RcAp} corresponding to $(\calP,A)$-\RcAp\ is defined similarly. 

Recall that a cardinal $\kappa$ is {\It ultrahuge}\Label{ultrahuge} if for any 
$\lambda>\kappa$, there are $j$, $M\subseteq\uniV$ \st\ $\Elembed{j}{\uniV}{M}{\kappa}$,
$j(\kappa)>\lambda$, $\fnsp{j(\kappa)}{M}\subseteq M$ and $V_{j(\lambda)}\subseteq M$. 

For an iterable class $\calP$ of \pos, a cardinal $\kappa$ is said to 
be ({\It tightly}) {\It $\calP$-Laver-generically ultrahuge} ((tightly) $\calP$-Laver-gen.\ 
ultrahuge, for short), if, for 
  any $\lambda>\kappa$ and $\poP\in\calP$ there is a $\poP$-name $\utpoQ$ with 
$\forces{\poP}{\utpoQ\in\calP}$, \st\  
  for $(\uniV,\poP\ast\utpoQ)$-generic $\genH$, there  
  are $j, M\subseteq\uniV[\genH]$ 
  \st\ $\Elembed{j}{\uniV}{M}{\kappa}$, $j(\kappa)>
  \lambda$, 
  $\poP,\genH$, {\ifprivate\privatecolor\else\darkred\fi$j\imageof{j(\kappa)}$},
  $(V_{j(\lambda)})^{\uniV[\genH]}\in M$ 
  (and
    $\poP\ast\utpoQ$ is forcing 
  equivalent to a \po\ of size $\LE j(\kappa)$).\footnote{\phantomsection\Label{fn-0} In the 
    following, we shall 
    denote this condition simply by 
  ``$\cardof{\poP\ast\utpoQ}\leq\lambda$''. 
  More generally, we shall always write ``$\cardof{\poP}\leq\lambda$'' for a \po\ $\poP$ to 
  mean that ``$\poP$ is forcing equivalent to a \po\ of size $\leq\lambda$".}

  By Theorem 6.5 in \cite{future}, any known Laver-generic large cardinal axiom 
  (formalizable in a single formula) does not imply the full version of Recurrence Axiom 
  (even for empty set as the set of parameters). However we have 
  the following:\ifextended\smallskip\fi 
  \memo{Adjust this to the latest version of \cite{future} (The theorem has the label: {\tt p-para-max-2})}

\begin{Thm}
  \Label{p-Lg-RcA-0} Suppose that $\kappa$ is tightly $\calP$-Laver-gen.\ ultrahuge for an 
  iterable class $\calP$ of \pos
  . Then
  $(\calP,\calH(\kappa))_{\Sigma_2}$-\RcAp\ holds. 
\end{Thm}

For the proof of \Thmof{p-Lg-RcA-0}, we use the following lemma which  will be also applied 
several times later in this article. 

Note that the ``sufficiently large finite fragment of ZFC'' in the following\ 
Lemma \ref{p-Lg-RcA-0-0} depends on how  
$\poP$-names for a \po\  $\poP$ are introduced. We will assume that they are introduced as 
in Kunen \cite{kunen-2011} and so that the arguments with "nice names" in the sense of \cite{kunen-2011} 
is available here.

\ifprivate\addcontentsline{toc}{section}{** Vα[G]=Vα\^{ }\{V[G]\}}\fi
\begin{Lemma}
  \Label{p-Lg-RcA-0-0} If $\alpha$ is a limit ordinal and $V_\alpha$ satisfies a 
  sufficiently large finite fragment of \ZFC, then for any $\poP\in V_\alpha$ 
  and $(\uniV,\poP)$-generic $\genG$,  
  we have $V_\alpha[\genG]={V_\alpha}^{\uniV[\genG]}$.\\
  \vspace{-1.2ex}\ifextended\else\qed\fi
\end{Lemma}
\memox{connect this Lemma with recurrence5? Corollary 25.
to Ikegami + Trang Theorem 1.8}
{
\ifextended\extendedcolor
\prf ``$\subseteq$'': This inclusion holds without the condition  
on the fragment of \ZFC. Also the condition ``$\poP\in V_\alpha$'' is irrelevant for this 
inclusion. 

We show by induction on $\alpha\in\On$ that
$V_\alpha[\genG]\subseteq {V_\alpha}^{\uniV[\genG]}$ holds for all $\alpha\in\On$.

The induction steps for $\alpha=0$ and limit ordinals $\alpha$ are trivial. So we assume that
$V_\alpha[\genG]\subseteq {V_\alpha}^{\uniV[\genG]}$ holds and show that the same inclusion 
holds for $\alpha+1$. Suppose $a\in V_{\alpha+1}[\genG]$. Then $a=\uta^\genG$ for 
a $\poP$-name $\uta\in V_{\alpha+1}$. Since $\uta\subseteq V_\alpha$, each
$\pairof{\utb,\condp}\in\uta$ is an element of $V_\alpha$. By induction hypothesis, it 
follows that $\utb^\genG\in {V_\alpha}^{\uniV[\genG]}$. It follows that 
$\uta^\genG\subseteq{V_\alpha}^{\uniV[\genG]}$. Thus 
$a=\uta^\genG\in{V_{\alpha+1}}^{\uniV[\genG]}$. \smallskip 

``$\supseteq$'':  
Suppose that $a\in {V_\alpha}^{\uniV[\genG]}$. Note that we can choose the ``sufficiently 
large finite fragment of \ZFC'' which should hold in $V_\alpha$, \st\ this implies that 
{\ifextended\extendedcolor $(*)$\fi}
${V_\alpha}^{\uniV[\genG]}$ still satisfies a large enough fragment of \ZFC, although the 
fragment may be different from the one $V_\alpha$ satisfies. In particular we find a 
cardinal $\lambda>\cardof{\poP}$ in ${V_\alpha}^{\uniV[\genG]}$ (and hence it is also a 
cardinal in
$\uniV[\genG]$) \st\
$a\in\calH(\lambda)^{V_\alpha^{\uniV[\genG]}}
\subseteq\calH(\lambda)^{\uniV[\genG]}\subseteq {V_\alpha}^{\uniV[\genG]}$. 
{\ifextended\extendedcolor {[}\,Note that
    $\calH(\lambda)^{{V_\alpha}^{\uniV[\genG]}}
    =\setof{a}{\cardof{\trcl(a)}<\lambda}^{{V_\alpha}^{\uniV[\genG]}}
    \subseteq\setof{a}{\cardof{\trcl(a)}<\lambda}^{\uniV[\genG]}
    =\calH(\lambda)^{\uniV[\genG]}.$ {]}
\fi}

Let
$a^*\in\calH(\lambda)^{\uniV[\genG]}$ be a transitive set \st\ $a\in a^*$. Then $a^*$ 
can be coded by a subset of $\lambda$. We can find the subset of $\lambda$ 
in $\uniV[\genG]$ and this subset has a nice $\poP$-name which is an element of
${V_\alpha}^\uniV$ since $\poP\in V_\alpha$. 
This shows that $a^*\in V_\alpha[\genG]$ and hence also $a\in V_\alpha[\genG]$. 
\qedofLemma\qedskip\fi}

\noindent
\prfof{\bfThmof{p-Lg-RcA-0}} Assume that $\kappa$ is tightly $\calP$-Laver gen.\ ultrahuge 
for an  iterable class $\calP$ of \pos. 

Suppose that $\varphi=\varphi(\overline{x})$ is $\Sigma_2$ formula (in $\Lin$),
$\overline{a}\in\calH(\kappa)$, and  
$\poP\in\calP$ is \st\ 
\begin{xitemize}
\xitem[x-Lg-RcA-a] 
  $\uniV\models\forces{\poP}{\varphi(\overline{a}\checked)}$. \
\end{xitemize}

Let $\lambda>\kappa$ be 
\st\ $\poP\in\uniV_\lambda$ and 
\begin{xitemize}
\xitem[x-Lg-RcA-0] $V_\lambda\prec^{}_{\Sigma_n}\uniV$ for a sufficiently large $n$. 
\end{xitemize}
In particular, we may assume that we have chosen the $n$ above so that a sufficiently large 
fragment of \ZFC\ holds in 
$V_\lambda$ (``sufficiently large fragment'' means here, in particular, in terms of \Lemmaof{p-Lg-RcA-0-0}). 

Let $\utpoQ$ be a $\poP$-name \st\ $\forces{\poP}{\utpoQ\in\calP}$, and for
$(\uniV,\poP\ast\utpoQ)$-generic $\genH$, there are $j$, $M\subseteq\uniV[\genH]$ with 
\begin{xitemize}
\xitem[x-Lg-RcA-1] 
  $\Elembed{j}{\uniV}{M}{\kappa}$, 
\xitem[x-Lg-RcA-1-0] 
  $j(\kappa)>\lambda$, 
\xitem[x-Lg-RcA-1-1] 
  $\poP\ast\utpoQ$, $\poP$, $\genH$, ${V_{j(\lambda)}}^{\uniV[\genH]}\in M$, and 
\xitem[x-Lg-RcA-1-2]
  $\cardof{\poP\ast\utpoQ}\leq j(\kappa)$. 
\end{xitemize}
By \xitemof{x-Lg-RcA-1-2}, we may assume that the underlying set 
of $\poP\ast\utpoQ$ is $j(\kappa)$ and $\poP\ast\utpoQ\in {V_{j(\lambda)}}^\uniV$. 

Let $\genG:=\genH\cap\poP$. 
Note that $\genG\in M$ by \xitemof{x-Lg-RcA-1-1} and 
we have 
\begin{xitemize}
\xitem[x-Lg-RcA-2] 
  ${V_{j(\lambda)}}^M\ubecause{=}{}{by \xitemof{x-Lg-RcA-1-1}}
  {V_{j(\lambda)}}^{\uniV[\genH]}
  \obecause{=}{1.44ex}{\qquad\qquad\qquad\qquad\qquad\qquad\qquad\qquad\qquad\qquad\qquad
  \vbox{\hbox{Since ${V_{j(\lambda)}}^M$ ($={V_{j(\lambda)}^{\uniV[\genH]}}$) satisfies a 
    sufficiently large fragment of 
    \ZFC\vspace{-0.18ex}}\hbox{by elementarity of $j$, and hence the equality follows by 
        \Lemmaof{p-Lg-RcA-0-0}}}} 
  {V_{j(\lambda)}}^\uniV[\genH]$. 
\end{xitemize}
Thus, by the definability of grounds and by \xitemof{x-Lg-RcA-1-1}, we have 
${V_{j(\lambda)}}^\uniV\in M$ and ${V_{j(\lambda)}}^\uniV[\genG]\in M$. 
{\addtocounter{Thm}{-1}
\begin{Claim}
  \Label{cl-Lg-RcA-0} ${V_{j(\lambda)}}^\uniV[\genG]\models\varphi(\overline{a})$. 
\end{Claim}
\noindent
\prfofClaim
By \Lemmaof{p-Lg-RcA-0-0}, ${V_\lambda}^{\uniV}[\genG]={V_\lambda}^{\uniV[\genG]}$, and
${V_{j(\lambda)}}^{\uniV}[\genG]={V_{j(\lambda)}}^{\uniV[\genG]}$ by \xitemof{x-Lg-RcA-0} and 
\xitemof{x-Lg-RcA-2}.  
By \xitemof{x-Lg-RcA-0}, both ${V_\lambda}^\uniV[\genG]$ and $V_{j(\lambda)}^\uniV[\genG]$ 
satisfy large enough fragment of \ZFC. In particular,
\begin{xitemize}
\xitem[x-Lg-RcA-2-0] 
  ${V_\lambda}^\uniV[\genG]\prec_{\Sigma_1}{V_{j(\lambda)}}^\uniV[\genG]$. 
\end{xitemize}
By \xitemof{x-Lg-RcA-a} and \xitemof{x-Lg-RcA-0}, we have
${V_\lambda}^\uniV[\genG]\models\varphi(\overline{a})$. By \xitemof{x-Lg-RcA-2-0} and since 
$\varphi$ is $\Sigma_2$, it follows that ${V_{j(\lambda)}}^\uniV[\genG]\models\varphi(\overline{a})$. 
\qedofClaim\qedskip
\addtocounter{Thm}{1}}

Thus we have 
\begin{xitemize}
\xitem[x-Lg-RcA-3] 
  $M\modelof{\mbox{there is a }
  \calP\mbox{-ground }N\mbox{ of }V_{j(\lambda)}\mbox{ with }N\models\varphi(\overline{a})}$.
\end{xitemize}

By the elementarity \xitemof{x-Lg-RcA-1}, it follows that 
\begin{xitemize}
\xitem[x-Lg-RcA-4] 
  $\uniV\modelof{\mbox{there is a }
  \calP\mbox{-ground }N\mbox{ of }V_{\lambda}\mbox{ with }N\models\varphi(\overline{a})}$.
\end{xitemize}

Now by \xitemof{x-Lg-RcA-0}, it follows that there is a $\calP$-ground $\uniW$ of $\uniV$ \st\
$\uniW\models\varphi(\overline{a})$. \qedof{\Thmof{p-Lg-RcA-0}}\qedskip

Some instances of weak forms of Recurrence Axioms decide the size of the continuum.
\memox{Erase this if
  $\kappa_\refl$  
  will be defined elsewhere in \rlap{the paper.}}
\begin{Thm}
  \Label{p-Lg-RcA-1} Assume that $\calP$ is an iterable class of \pos. \wassertof{1} If 
  $\calP$ contains a \po\ which adds a real (over 
  the universe), then 
  $(\calP,\calH(\kappa_\refl))_{\Sigma_1}$-\RcA\ implies $\neg\CH$.\smallskip

  \wassert{2} Suppose that $\calP$ contains a 
  \po\ which 
  forces ${\aleph_2}^\uniV$ to be equinumerous with ${\aleph_1}^\uniV$. Then 
  $(\calP,\calH(2^{\aleph_0}))_{\Sigma_1}$-\RcA\ implies $\continuum\leq\aleph_2$.
  \smallskip

  \wassert{2'} If $\calP$ contains a \pos\ which forces ${\aleph_2}^\uniV$ to be 
  equinumerous with 
  ${\aleph_1}^\uniV$, then $(\calP,\calH((\aleph_2)^+))_{\Sigma_1}$-\RcA\ does not hold. 
  \smallskip

  \wassert{3} If $(\calP,\calH(\kappa_\refl))_{\Sigma_1}$-\RcA\ holds then all 
  $\poP\in\calP$ preserve $\aleph_1$ and they are also stationary preserving.\smallskip

  \wassert{4} If $\calP$ contains a \po\ which 
  adds a real as well as a \po\ which collapses 
  ${\aleph_2}^\uniV$, then 
  $(\calP,\calH(\kappa_\refl))_{\Sigma_1}$-\RcA\ implies $\continuum=\aleph_2$.
  \smallskip

  \wassert{5} If $\calP$ contains a 
  \po\ which 
  collapses ${\aleph_1}^\uniV$, then $(\calP,\calH(2^{\aleph_0}))_{\Sigma_1}$-\RcA\ implies 
  \CH.\smallskip 

  \wassert{5'} If $\calP$ contains a \po\ which collapses ${\aleph_1}^\uniV$ then 
  $(\calP,\calH((\continuum)^+))_{\Sigma_1}$-\RcA\ does not hold. \smallskip

  \wassert{6} Suppose that all $\poP\in\calP$ preserve cardinals and $\calP$ contains \pos\ 
  adding at least $\kappa$ many reals for each $\kappa\in\Card$. Then
  $(\calP,\emptyset)_{\Sigma_2}$-\RcAp\ implies that $\continuum$ is very large. More 
  precisely, we have $\continuum>\aleph_\alpha$ for 
  any $\calP$-absolutely $\Sigma_2$-definable ordinal $\alpha$. 

  \wassert{6'} Suppose that $\calP$ is as in \assertof{6}. Then 
  $(\calP,\calH(\continuum))_{\Sigma_2}$-\RcAp\ implies that $\continuum$ is a limit 
  cardinal. Thus if $\continuum$ is regular in addition, then $\continuum$ is weakly inaccessible.

\end{Thm}\memox{Add \assertof{6} for $\calP=$ cardinal preserving adding arbitrary number of 
  reals. See RIMS2024-set-theory-fuchino.tex}
\prf \assertof{1}: Assume that $\calP$ is an iterable class of \pos\ containing a 
\po\ $\poP$ adding a real and $(\calP,\calH(\kappa_\refl))_{\Sigma_1}$-\RcA\ holds. If \CH\ 
holds, then $\psof{\omega}^\uniV\in\calH(\kappa_\refl)$. Hence 
``$\exists x\,(x\subseteq\omega\land x\not\in\psof{\omega}^\uniV)$'' is 
a $\Sigma_1$-formula with parameters from $\calH(\kappa_\refl)$ and $\poP$ forces (the 
formula in forcing language corresponding to) this 
formula. 

By $(\calP,\calH(\kappa_\refl))_{\Sigma_1}$-\RcA, the formula must hold in a ground. This 
is a contradiction. \smallskip

\assertof{2}: Assume that $(\calP,\calH(2^{\aleph_0}))_{\Sigma_1}$-\RcA\ holds and
$\poP\in\calP$ 
forces ${\aleph_2}^\uniV$ to be equinumerous with ${\aleph_1}^\uniV$. 
If $\continuum>\aleph_2$ then 
${\aleph_1}^\uniV$, ${\aleph_2}^\uniV\in\calH(2^{\aleph_0})$. Letting 
$\psi(x,y)$ a $\Sigma_1$-formula stating 
``$\exists f\,(f\mbox{ is a surjection from }x\mbox{ to }y)$'', 
we have 
$\forces{\poP}{\psi(({\aleph_1}^\uniV)\checked,
({\aleph_2}^\uniV)\checked)}$. 

By $(\calP,\calH(2^{\aleph_0}))_{\Sigma_1}$-\RcA, the formula 
$\psi({\aleph_1}^\uniV,{\aleph_2}^\uniV)$
must hold in a ground. This is a contradiction. \smallskip

\assertof{2'}: Assume that $\poP\in\calP$ is \st\
$\forces{\poP}{\cardof{{\aleph_2}^\uniV}=\cardof{{\aleph_1}^\uniV}}$, and
$(\calP,\calH({\aleph_2}^+))_{\Sigma_1}$-\RcA\ holds. Then, since $\aleph_1$,
$\aleph_2\in\calH({\aleph_2}^+)$ and ``$\cardof{x}=\cardof{y}$'' is $\Sigma_1$, there is a 
ground $\uniW$ of $\uniV$ \st\
$\uniW\models\cardof{{\aleph_2}^\uniV}=\cardof{{\aleph_1}^\uniV}$. This is a 
contradiction.\smallskip

\assertof{3}: Suppose that $\poP\in\calP$ is \st\
$\forces{\poP}{{\aleph_1}^\uniV\mbox{ is countable}}$. Note that
$\omega, \aleph_1\in\calH(\kappa_\refl)$. 
By
$(\calP,\calH(\kappa_\refl))_{\Sigma_1}$-\RcA, it follows that there is a ground $\uniW$ of 
$\uniV$ \st\ $\uniW\modelof{{\aleph_1}^\uniV\mbox{ is countable}}$. This is a contradiction.

Suppose now that $S\subseteq\omega_1$ is stationary and $\poP\in\calP$ destroys the 
stationarity of $S$. Note that 
$\omega_1$, $S\in\calH(\aleph_2)$. Let $\varphi=\varphi(y,z)$ be the $\Sigma_1$-formula
\begin{xitemize}
\item[] $\exists x\,(x\mbox{ is a club subset of the ordinal }y\mbox{ and }
  z\cap x=\emptyset)$. 
\end{xitemize}
Then we have $\forces{\poP}{\varphi(\omega_1, S)}$. By
$(\calP,\calH(\kappa_\refl))_{\Sigma_1}$-\RcA, it follows that there is a ground
$\uniW\subseteq\uniV$ \st\ $S\in\uniW$ and $\uniW\models\varphi(\omega_1,S)$. This is a 
contradiction to the stationarity of $S$. 

\smallskip

\assertof{4}: follows from \assertof{1}, \assertof{2} and \assertof{3}. 
\smallskip

\assertof{5}: 
If $\aleph_1<\continuum$, 
then ${\aleph_1}^\uniV\in\calH(2^{\aleph_0})$.

Let $\poP\in\calP$ be a \po\ collapsing ${\aleph_1}^\uniV$. That is, 
$\forces{\poP}{{\aleph_1}^\uniV\mbox{ is countable}}$. Since ``$\cdots$ is countable'' is
$\Sigma_1$, there is a ground $M$ \st\ $M\modelof{{\aleph_1}^\uniV\mbox{ is countable}}$ by
$(\calP, \calH(2^{\aleph_0}))_{\Sigma_1}$-\RcA.  
This is a contradiction. \smallskip

\assertof{5'}: Assume that $\poP\in\calP$ is \st\
$\forces{\poP}{{\aleph_1}^\uniV\mbox{ is countable}}$, 
and $(\calP,\calH((\continuum)^+))$-\RcA\ holds. Since $\aleph_1\in\calH((\continuum)^+)$, 
it follows that there is a ground $\uniW$ of $\uniV$ \st\
$\uniW\models{\aleph_1}^\uniV\mbox{ is countable}$. This is a contradiction.\smallskip

\assertof{6}: To prove e.g.\ that $\continuum>\aleph_\omega$, let
  $\poP\in\calP$ be \st\  
  $\forces{\poP}{\continuum>\aleph_\omega}$. Then by
  $(\calP,\emptyset)_{\Sigma_2}{\mbox{-}}$\RcAp, there is a $\calP$-ground $\uniW$ 
  of $\uniV$ \st\ $\uniW\models\continuum>\aleph_\omega$. Since $\uniV$ is $\calP$-gen.\ 
  extension of $\uniW$ and $\calP$ preserves cardinals, it follows that
  $\uniV\models\continuum>\aleph_\omega$.\smallskip

  \assertof{6'}: Suppose $\mu<\continuum$. Then $\mu\in\calH(\continuum)$. There is $\poP\in\calP$ 
  \st\ $\forces{\poP}{\continuum>\mu^+}$. By
  $(\calP,\calH(\continuum))_{\Sigma_2}\mbox{-}$\RcAp, it follows that there is a $\calP$-ground 
$\uniW$ of $\uniV$ which satisfies this statement. Since $\calP$ preserves cardinals it 
  follows that $\uniV\models\continuum>\mu^+$. 
\qedofThm\qedskip

In contrast to \Thmof{p-Lg-RcA-1}, $(ccc,\calH(\kappa_\refl))$-\RcA\  
does 
not decide the size of the continuum (see \Thmof{p-Lg-RcA-4} and \Thmof{p-Lg-RcA-5}).

Recurrence Axioms can be considered as natural requirements claiming that a 
reflection holds from the set-generic multiverse down to the set-generic archaeology. 
From the standpoint that we should adopt the maximal amount of the Recurrence in the 
(ultimate) extension of \ZFC, we arrive at either 
$(\mbox{all \pos}, \calH(2^{\aleph_0}))$-\RcA\  
or $(\mbox{semi-proper \pos}, \calH(\kappa_\refl))$-\RcA\ according to 
\Thmof{p-Lg-RcA-1}, and these  
axioms imply \CH\ or $2^{\aleph_0}=\aleph_2$, respectively. In particular, these two axioms 
(or axiom schemes to be more precise) are not compatible to each other. 

In \sectionof{LGM}, we shall 
examine an axiom(scheme) which implies the full 
$(\xmbox{semi-proper \pos}, \calH(\kappa_\refl))$-\RcA$^+$ and also a large fragment of 
$(\mbox{all \pos}, \calH(\aleph_1)$-\RcA$^+$ as well as $\MM^{++}$. 
\memo{What is the axiom saying all (definable?) cardinals in $\uniV$ are 
  super-$C^{(\infty)}$ something 
  in the bedrock?  
  What is the axiom claiming that stationarity of any subclass of $\On$ is reflected down 
  to a regular cardinal. The latter axiom follows from the former?}

\ifextended\else\def\texorpdfstring#1#2{#1}\fi
\section{Tightly super-\texorpdfstring{$C^{(\infty)}$}{C(∞)}-Laver generic ultrahuge cardinal}
\Label{c-infty}
In \cite{future}, it is shown that the existence of a (tightly) $\calP$-Laver-gen.\ large 
cardinal 
does not imply Maximality Principle even without parameters. The proof in \cite{future} 
can be modified to prove the non-implication of $(\calP,\emptyset)_{\Pi_3}$-\RcA\ from 
generic large cardinals of various sorts, and this also 
shows that ``$(\calP,\emptyset)_{\Sigma_2}$-\RcAp'' in \Thmof{p-Lg-RcA-0} is optimal. 

In this section we show that the existence of a strong variant of $\calP$-Laver generic 
large cardinal $\kappa$ we are  
going to introduce  below does imply $\MP(\calP, \calH(\kappa))$ (see \Thmof{p-Lg-RcA-5}). 
It is essential that the variant of Laver genericity (called ``the tightly super
$C^{(\infty)}$-Laver gen.\ large cardinal'' below) is formulated not in a single formula 
but rather in an axiom scheme. 

For a natural number $n$, we call a cardinal $\kappa$ {\It super-$C^{(n)}$-hyperhuge} if  
for any $\lambda_0>\kappa$ there are $\lambda\geq\lambda_0$ with 
$V_\lambda\prec_{\Sigma_n}\uniV$, and $j$, $M\subseteq\uniV$ \st\
$\Elembed{j}{\uniV}{M}{\kappa}$, $j(\kappa)>\lambda$, $\fnsp{j(\lambda)}{M}\subseteq M$ and
$V_{j(\lambda)}\prec_{\Sigma_n}\uniV$. 

$\kappa$ is {\It super-$C^{(n)}$-ultrahuge} if the condition above holds with
``$\,\fnsp{j(\lambda)}{M}\subseteq M$'' replaced by ``$\,\fnsp{j(\kappa)}{M}\subseteq M$ and 
$V_{j(\lambda)}\subseteq M$''. 

If $\kappa$ is super-$C^{(n)}$-hyperhuge then it is super-$C^{(n)}$-ultrahuge. This can be 
shown similarly to \Lemmaof{p-bedrock-0} in \sectionof{bedrock}. 


We shall also say that $\kappa$ is {\It super-$C^{(\infty)}$-hyperhuge} ({\It super
$C^{(\infty)}$-ultrahuge}, resp.) if it is super-$C^{(n)}$-hyperhuge (super
$C^{(n)}$-ultrahuge, resp.) for all natural number $n$. 

A similar kind of strengthening of the notions of large cardinals which 
we call here ``super-$C^{(n)}$'' appears also in
Boney \cite{boney}. It is called in \cite{tsaprounis3} and \cite{boney}
``$C^{(n)+}$'' and the notion is considered there in connection with extendibility. 

For a natural number $n$ and an iterable class $\calP$ of \pos, 
a cardinal $\kappa$ is {\It super-$C^{(n)}$-$\calP$-Laver-generically ultrahuge} (super
$C^{(n)}$-$\calP$-Laver-gen.\ ultrahuge, for short) if, for 
any $\lambda_0>\kappa$ 
and for any 
$\poP\in\calP$,
there are a $\lambda\geq\lambda_0$ 
with $V_\lambda\prec_{\Sigma_n}\uniV$, a $\calP$-name $\utpoQ$ with
$\forces{\poP}{\utpoQ\in\calP}$ \st\ for $(\uniV,\poP\ast\utpoQ)$-generic $\genH$, there 
are $j$, $M\subseteq \uniV[\genH]$ with $\Elembed{j}{\uniV}{M}{\kappa}$,
$j(\kappa)>\lambda$, $\poP$, $\genH$, {\ifextended\privatecolor\else\darkred\fi$j\imageof{j(\lambda)}$}, ${V_{j(\lambda)}}^{\uniV[\genH]}\in M$ and
${V_{j(\lambda)}}^{\uniV[\genH]}\prec_{\Sigma_n}\uniV[\genH]$. 

A super-$C^{(n)}$-$\calP$-Laver-generically ultrahuge cardinal $\kappa$ is {\It tightly 
    super-$C^{(n)}$-$\calP$-Laver-generically ultrahuge} (tightly 
    super-$C^{(n)}$-$\calP$-Laver-gen.\ ultrahuge, for short), if  
$\cardof{\poP\ast\utpoQ}\leq j(\kappa)$ (see \footnoteof{fn-0}).

{\It super-$C^{(\infty)}$- $\calP$-Laver-gen.\ ultrahugeness} and 
{\It tightly super-$C^{(\infty)}$- $\calP$-Laver gen.\ ultrahugeness}
are defined similarly to super-$C^{(\infty)}$-ultrahugeness. 

Note that, in general, super-$C^{(\infty)}$-hyperhugeness 
and super-$C^{(\infty)}$-ultrahugeness are notions unformalizable in the language 
of \ZFC\ without introducing a new constant symbol for $\kappa$ since we need infinitely 
many $\Lin$-formulas to formulate them. Exceptions are when we 
are talking about a cardinal in a set model being with one of these 
properties like in \Lemmaof{p-Lg-RcA-2} below (and in such a case ``natural number $n$'' 
actually refers to ``$n\in\omega$''), or when we are talking about a cardinal definable in
$\uniV$ having these properties in an inner model like in \Thmof{p-bedrock-2-0} or \Corof{p-bedrock-8}.\memox{???} In the latter case, 
the situation is formalizable with infinitely may $\Lin$-sentences.

In contrast, the super-$C^{(\infty)}$-$\calP$-Laver gen.\ ultrahugeness 
of $\kappa$ is expressible in infinitely many $\Lin$-sentences. 
This is because a $\calP$-Laver gen.\ large 
cardinal $\kappa$ for relevant classes $\calP$ of \pos\ 
is uniquely determined as $\kappa_{\refl}$ or $2^{\aleph_0}$ (see e.g.\ \cite{sfetal-II} or 
\cite{future}). 

\begin{Lemma}\Label{p-Lg-RcA-1-a} 
  Suppose that $\kappa$ is super-$C^{(n)}$-ultrahuge cardinal. Then we have
  $V_\kappa\prec_{\Sigma_{n+1}}\uniV$. In particular, in a context in which we can 
  express that $\kappa$ is super-$C^{(\infty)}$-ultrahuge cardinal in a (set or class) model $\uniW$, 
  we have ${V_\kappa}^\uniW\prec\uniW$.  
\end{Lemma}
\prf We check that $V_\kappa$ passes Vaught's test. Let $\varphi(x,\overline{y})$ be a
$\Pi_n$-formula. Suppose that $\overline{b}\in V_\kappa$ and $a\in\uniV$ are \st\
$\uniV\models\varphi(a,\overline{b})$. We want to show that there is $a'\in V_\kappa$ \st\
$\uniV\models\varphi(a',\overline{b})$.

Let $\lambda$ be \st\ 
\begin{xitemize}
\xitem[x-Lg-RcA-4-a-0] $a\in V_\lambda$,
\xitem[x-Lg-RcA-4-a-1] $V_\lambda\prec_{\Sigma_n}\uniV$,
\xitem[x-Lg-RcA-4-a-2] there are $j,M\subseteq\uniV$ \st\ 
  \begin{xitemize}
    \xxitem[x-Lg-RcA-4-a-2][a]\quad $\Elembed{j}{\uniV}{M}{\kappa}$,\qquad
    \xxitemof{x-Lg-RcA-4-a-2}{b}\quad $j(\kappa)>\lambda$,\qquad
    \xxitemof{x-Lg-RcA-4-a-2}{c}\quad $\fnsp{j(\kappa)}{M}\subseteq M$, \qquad
    \xxitem[x-Lg-RcA-4-a-2][d]\quad $V_{j(\lambda)}\subseteq M$, and\quad
    \xxitemof{x-Lg-RcA-4-a-2}{e}\quad $V_{j(\lambda)}\prec_{\Sigma_n}\uniV$.
  \end{xitemize}
\end{xitemize}
By \xxitemof{x-Lg-RcA-4-a-2}{e}, we have $V_{j(\lambda)}\models\varphi(a,\overline{b})$. By 
\xxitemof{x-Lg-RcA-4-a-2}{d}, it follows that
$M\modelof{V_{j(\lambda)}\models\varphi(a,\overline{b})}$. Noting that
$\overline{b}=j(\overline{b})$, it follows that 
\begin{xitemize}
\item[] $M\modelof{\mbox{\,there is }x\in V_{j(\kappa)}\mbox{ \st\ }
  V_{j(\lambda)}\models\varphi(x,j(\overline{b}))}$.
\end{xitemize}
By the elementarity \xxitemof{x-Lg-RcA-4-a-2}{a}, it follows
\begin{xitemize}
\item[] $\uniV\modelof{\mbox{\,there is }x\in V_{\kappa}\mbox{ \st\ }
  V_{\lambda}\models\varphi(x,\overline{b})}$.
\end{xitemize}
Let $a'\in V_\kappa$ be a witness of this. Then by \xitemof{x-Lg-RcA-4-a-1} we have
$\uniV\models\varphi(a',\overline{b})$. \qedofLemma\qedskip

\begin{Cor}
  \Label{p-Lg-RcA-1-a-0}
  Suppose that $\kappa$ is a super-$C^{(\infty)}$-ultrahuge (-hyperhuge, resp.) cardinal in 
  a (set or class)  
  model $\uniW$. Then for each $n\in\natnums$,
  $\uniW\modelof{\,\xmbox{there are stationarily many super-}C^{(n)}
\xmbox{-ultrahuge (-hyperhuge, resp.) cardinals}}$.
\end{Cor}
\prf Suppose that $\kappa$ is a super-$C^{(\infty)}$-ultrahuge cardinal in $\uniW$.  
By \Lemmaof{p-Lg-RcA-1-a}, we have ${V_\kappa}^\uniW\prec\uniW$. Suppose that
$\varphi=\varphi(x,y)$ is an $\Lin$-formula and $a\in {V_\kappa}^\uniW$ \st\
${V_\kappa}^\uniW\modelof{\varphi(\cdot,b)\mbox{ is a club }\subseteq\On}$. Then, by 
elementarity, $\uniW\modelof{\varphi(\cdot,b)\mbox{ is a club }\subseteq\On}$ and
$\uniW\models\varphi(\kappa,b)$. Thus
$\uniW\modelof{\,\xmbox{there is a super-}C^{(n)}
\xmbox{-ultrahuge cardinal }\mu\xmbox{ \st\ }\varphi(\mu,b)}$. By elementarity, it follows 
that ${V_\kappa}^\uniW\modelof{\,\xmbox{there is a super-}C^{(n)}
  \xmbox{-ultrahuge cardinal }\mu\xmbox{ \st\ }\varphi(\mu,b)}$.

Since $b$ was arbitrary, it follows that 
${V_\kappa}^\uniW\modelof{\forall y\,(\xmbox{if }\varphi(\cdot,y)
  \xmbox{ is a club in }\On,\xmbox{ then there is a super-}C^{(n)}
  \xmbox{-ultrahuge cardinal }\mu\xmbox{ \st\ }\varphi(\mu,y))}$.
By elementarity, the same statement also holds in $\uniW$. 
\qedofCor\qedskip

An ultrafilter $U\subseteq\psof{\psof{\lambda^*}}$ is said to be normal if 
\begin{xitemize}
\xitem[x-Lg-RcA-4-0] $\setof{x\in\psof{\lambda^*}}{\alpha\in x}\in U$ for all
  $\alpha\in\lambda^*$, and 
\xitem[x-Lg-RcA-4-1]  for any $\seqof{X_\alpha}{\alpha<\lambda^*}\in\fnsp{\lambda^*}{U}$,
  $\bigtriangleup_{\alpha<\lambda^*}X_\alpha\in U$.
\end{xitemize}

Under \xitemof{x-Lg-RcA-4-0}, the condition \xitemof{x-Lg-RcA-4-1} is equivalent to
\begin{xitemize}
\xitem[x-Lg-RcA-4-2] 
  For any $X\in U$ and any regressive\footnotemark\ $\mapping{f}{X}{\lambda^*}$, 
  there is $X'\in U$ \st\ $X'\subseteq X$ and $f$ is constant on $X'$.
\end{xitemize}
\footnotetext{$\mapping{f}{X}{\lambda^*}$ is {\It regressive} if $f(x)\in x$ for all
  $x\in X$. }

\begin{Lemma}
  \Label{p-Lg-RcA-1-0} For $\kappa<\lambda<\kappa^*<\lambda^*$, the following \assertof{A} 
  and \assertof{B} are equivalent. \smallskip

  \wassert{A} there are $j$, $M\subseteq\uniV$ \st\ $\Elembed{j}{\uniV}{M}{\kappa}$,
  $j(\kappa)=\kappa^*$, $j(\lambda)=\lambda^*$, $\fnsp{j(\lambda)}{M}\subseteq M$.\smallskip

  \wassert{B} there is a $\kappa$-complete normal ultrafilter
  $U\subseteq\psof{\psof{\lambda^*}}$ \st\vspace{-0.72ex}
  \begin{xitemize}
  \xitem[x-Lg-RcA-4-3] 
    $X^*:=\setof{x\in\psof{\lambda^*}}{x\cap\kappa\in\kappa,\,
      \otp(x\cap\kappa^*)=\kappa,\,\otp(x)=\lambda}\in U$. 
  \end{xitemize}

\end{Lemma}
\noindent
{\bf A Sketch of the Proof:}
\memo{See Lemma 14.2 {\tt \{p-hyper-2\}} in \cite{math-20}}
\assertof{A} $\Rightarrow$ \assertof{B}: For $j$ as in 
\assertof{A}, the ultrafilter $U$ defined by
\begin{xitemize}
\xitem[x-Lg-RcA-4-3-0] $U:=\setof{X\in\psof{\psof{\lambda^*}}}{j\imageof{\lambda^*}\in j(X)}$
\end{xitemize}
satisfies \assertof{B}.\smallskip

\assertof{B} $\Rightarrow$ \assertof{A}: For $U$ as in \assertof{B}, the elementary 
embedding $j_U$ defined by 

\begin{xitemize}
\xitem[x-Lg-RcA-4-4]
  $\Elembed{j_U}{\uniV}{N:=mcol(\fnsp{\psof{\lambda^*}}{\uniV}/U)}{\kappa}$;
  $x\ \mapsto\ [c_x]$ 
\end{xitemize}
satisfies the conditions in \assertof{A}. 

Here, $mcol$ denotes the Mostowski collapse, $[f]$ the element of $N$ which 
corresponds to the equivalence class of $f\in\fnsp{\psof{\lambda^*}}{\uniV}$ and $c_x$ the 
function on $\psof{\lambda^*}$ whose value is constantly $x$. 
\qedofLemma
\qedskip

{\ifextended\extendedcolor 
  Similarly to the case of  
measurable cardinals or supercompact cardinals we can define another function $k_j$ 
associated with $j$ by:
\begin{xitemize}
\xitem[x-Lg-RcA-4-5] 
  $\mapping{k_j}{N}{M}$; $[f]\ \mapsto\ j(f)(j\imageof\lambda^*)$.
\end{xitemize}

\memo{see math-notes-20 Lemma 14.4 {\tt p-hyper-3-0}}.
\begin{LemmaA}\extendedcolor
  \Label{p-Lg-RcA-1-0-0} Suppose that $\kappa<\lambda<\kappa^*<\lambda^*$ and
  $\Elembed{j}{\uniV}{M\subseteq\uniV}{\kappa}$ are as in 
  {\rm\Lemmaof{p-Lg-RcA-1-0}},\,\assertof{A}. Let $U$ be defined by \xitemof{x-Lg-RcA-4-3-0}. 
  Then   
  $U$ satisfies the condition \xitemof{x-Lg-RcA-4-3-0} in 
  {\rm\Lemmaof{p-Lg-RcA-1-0}},\,\assertof{B}.  

  Let  $j_U$ and $k_j$ be defined as above. Then 
  \smallskip

  \wassert{1} $k_j$ is well-defined.\smallskip

  \wassert{2} $k_j$ is an elementary embedding.
  \smallskip

  \wassert{3} $k_j\circ j_U=j$.\smallskip

  \wassert{4} 
  $k_j\restr\calH((\lambda^*)^+)=id_{\calH((\lambda^*)^+)}$. \smallskip

  \wassert{5} $j\restr\calH(\lambda^+)=j_U\restr\calH(\lambda^+)$. 
\end{LemmaA}\fi}

In the following we shall use a (local) notation
``$(\kappa,\lambda)^{(\kappa^*,\lambda^*)}$'' to denote the condition 
\assertof{B} in \Lemmaof{p-Lg-RcA-1-0}.  Clearly we have 

\begin{Lemma}
  \Label{p-Lg-RcA-1-1} For a cardinal $\kappa$ and a natural number $n$, $\kappa$ is a 
  super-$C^{(n)}$-hyperhuge cardinal\ \ $\Leftrightarrow$\ \ for any $\lambda_0>\kappa$ 
  there are $\lambda^*>\kappa^*>\lambda\geq\lambda_0$ \st\smallskip

  \wassert{C} $V_\lambda\prec_{\Sigma_n}\uniV$, 
  $V_{\lambda^*}\prec_{\Sigma_n}\uniV$, and $(\kappa,\lambda)^{(\kappa^*,\lambda^*)}$. 
  \qed
\end{Lemma}

We shall denote the condition \assertof{C} by $[\kappa,\lambda]^{[\kappa^*,\lambda^*,n]}$. 
This is also merely a local notation. 

\begin{Lemma}
  \Label{p-Lg-RcA-2} If $\kappa$ is $2$-huge with the $2$-huge elementary embedding $j$, 
  that is,\ there is $M\subseteq\uniV$ \st\
  $\Elembed{j}{\uniV}{M\subseteq \uniV}{\kappa}$,  and  
  \begin{xitemize}
  \xitem[x-Lg-RcA-5] 
    $\fnsp{j^2(\kappa)}{M}\subseteq M$, 
  \end{xitemize}
  then $V_{j(\kappa)}\modelof{\kappa
    \mbox{ is super-}C^{(\infty)}\mbox{-hyperhuge cardinal\/}}$, 
  and for each $n\in\omega$,\\
  $V_{j(\kappa)}\modelof{\mbox{there are stationarily
      many super-}C^{(n)}\mbox{-hyperhuge cardinals\/}}$.
\end{Lemma}
\prf 
Suppose that $j$ is as above and $n\in\omega$. 

By elementarity of $j$ and since $\kappa=\crit(j)$, we have 
$V_\kappa\prec{V_{j(\kappa)}}^M$. By \xitemof{x-Lg-RcA-5}, it follows that 
\begin{xitemize}
\xitem[x-Lg-RcA-6] $V_\kappa\prec V_{j(\kappa)}$. 
\end{xitemize}
By elementarity of $j$, this implies
  $M\modelof{V_{j(\kappa)}\prec V_{j(j(\kappa))}}$. 
By the closedness \xitemof{x-Lg-RcA-5}
  of $M$ (and since $j(j(\kappa))$ is inaccessible), it follows that
\begin{xitemize}
\xitem[x-Lg-RcA-6-a] 
  $V_{j(\kappa)}\prec V_{j(j(\kappa))}$.
\end{xitemize}

For a $\lambda_0$ with $\kappa<\lambda_0<j(\kappa)$, let 
$\lambda_0\leq\lambda<j(\kappa)$ be \st\ 
\begin{xitemize}
\xitem[x-Lg-RcA-6-0] 
  $V_\lambda\prec_{\Sigma_n}V_{j(\kappa)}$. 
\end{xitemize}
By \xitemof{x-Lg-RcA-6-a}, it follows that
\begin{xitemize}
\xitem[x-Lg-RcA-6-1] 
  $V_\lambda\prec_{\Sigma_n}V_{j(j(\kappa))}$.  
\end{xitemize}

We also have $M\modelof{V_{j(\lambda)}\prec_{\Sigma_n}V_{j(j(\kappa))}}$ by
\xitemof{x-Lg-RcA-6-0} and by elementarity of $j$. 
By the closedness \xitemof{x-Lg-RcA-5} (and since $j(j(\kappa))$ is inaccessible), it 
follows that 
\begin{xitemize}
\xitem[x-Lg-RcA-6-0-0] 
  $V_{j(\lambda)}\prec_{\Sigma_n}V_{j(j(\kappa))}$.
\end{xitemize}

Let
\begin{xitemize}
\xitem[x-Lg-RcA-7] $U:=\setof{X\subseteq\psof{j(\lambda)}}{j\imageof{j(\lambda)}\in j(X)}$.
\end{xitemize}
The following is easy to check:
\begin{Claim}
  In $V_{j(j(\kappa))}$, $U$ witnesses 
  $(\kappa,\lambda)^{(j(\kappa),j(\lambda))}$.\quad  Hence we have\\
  $V_{j(j(\kappa))}\modelof{[\kappa,\lambda]^{[j(\kappa),j(\lambda),n]}}$ for all
  $n\in\omega $. 
\end{Claim}
{\ifextended\extendedcolor
\noindent
\prfofClaim
Clearly $U$ is a 
ultrafilter.

$U$ is $\kappa$-complete: Suppose that $\seqof{X_\xi}{\xi<\mu}\in\fnsp{\mu}{U}$ for some
$\mu<\kappa$. Then $j(\seqof{X_\xi}{\xi<\mu})=\seqof{j(X_\xi)}{\xi<\mu}$ by the 
elementarity of $j$ and $\mu<\crit(j)$. Since $X_\xi\in U$, we have $j\imageof{\lambda}\in j(X_\xi)$ 
for all $\xi$ by the definition \xitemof{x-Lg-RcA-7} of $U$. It follows that
$j\imageof{ j(\lambda)}\in\bigcap_{\xi<\mu}j(X_\xi)=j(\bigcap_{\xi<\mu}X_\xi)$. Thus
$\bigcup_{\xi<\mu}X_\xi\in U$.

$U\models$\xitemof{x-Lg-RcA-4-0}: For $\alpha< j(\lambda)$,
$j(\setof{x\in\psof{ j(\lambda)}}{\alpha\in x})=\setof{x\in\psof{j( j(\lambda))}}{j(\alpha)\in x}\ni 
j\imageof{ j(\lambda)}$. Hence $\setof{x\in\psof{ j(\lambda)}}{\alpha\in x}\in U$. 

$U\models$\xitemof{x-Lg-RcA-4-1}: Suppose that
$\vec{X}:=\seqof{X_\alpha}{\alpha< j(\lambda)}\in\fnsp{ j(\lambda)}{U}$. Then, by elementarity of
$j$, 
$j(\bigtriangleup_{\alpha< j(\lambda)}X_\alpha)=\setof{x\in\psof{j( j(\lambda))}}{\forall\alpha\in 
x\,(x\in (j(\vec{X}))_\alpha)}$. For $\alpha\in j\imageof{ j(\lambda)}$,
$j\imageof j(\lambda)\in j(X_\alpha)=(j(\vec{X}))_\alpha$ by the definition 
\xitemof{x-Lg-RcA-7} of $U$. It  
follows that $j\imageof{ j(\lambda)}\in j(\bigtriangleup_{\alpha< j(\lambda)}X_\alpha)$. 
This means $\bigtriangleup_{\alpha< j(\lambda)}X_\alpha\in U$  
by the definition \xitemof{x-Lg-RcA-7} of $U$. 

$X^*\in U$: Note that
\begin{xitemize}
\item[] 
  $j(X^*)=\setof{x\in\psof{j( j(\lambda))}}{{}
  \begin{array}[t]{@{}l}
    x\cap j(\kappa)\in j(\kappa),\,
    \otp(x\cap j(j(\kappa)))=j(\kappa),\,\\[\jot]
    \otp(x)={j(\lambda)}}
  \end{array}$
\end{xitemize}
by elementarity of $j$. Hence $j\imageof j(\lambda)\in j(X^*)$ and thus $X^*\in U$. 

This shows $V_{j(j(\kappa))}\modelof{(\kappa,\lambda)^{(j(\kappa),j(\lambda))}}$. By 
\xitemof{x-Lg-RcA-6-1} and \xitemof{x-Lg-RcA-6-0-0}, it follows that
$V_{j(j(\kappa))}\modelof{[\kappa,\lambda]^{[j(\kappa),j(\lambda), n]}}$
\qedofClaim
\qedskip
\fi}

By \Claimabove, $V_{j(j(\kappa))}\models\exists x\exists y\,([\kappa,\lambda]^{[x,y,n]})$ 
for all $n$. By \xitemof{x-Lg-RcA-6-a}, it follows that
$V_{j(\kappa)}\models\exists x\exists y\,([\kappa,\lambda]^{[x,y,n]})$. 
Since $n\in\omega$ and $\kappa<\lambda_0<j(\kappa)$ were arbitrary, it follows that 
$V_{j(\kappa)}\modelof{\kappa\mbox{ is super-}C^{(\infty)}\mbox{-hyperhuge}}$.

For a fixed $n$ and club $C\subseteq j(\kappa)$, we have
$V_{j(j(\kappa))}\modelof{\kappa\in j(C)
  \xmbox{ and }\kappa\xmbox{ is super-}C^{(n)}\mbox{-hyperhuge}}$. Thus 
$V_{j(j(\kappa))}\modelof{\exists x\,(x\in j(C) 
  \xmbox{ and }x\xmbox{ is super-}C^{(n)}\mbox{-hyperhuge})}$. It follows
that $V_{j(\kappa)}\modelof{\exists x\,(x\in C 
  \xmbox{ and }x\xmbox{ is super-}C^{(n)}\mbox{-hyperhuge})}$. 
This shows that $V_{j(\kappa)}\modelof{\mbox{there are stationarily
      many super-}C^{(n)}\mbox{-hyperhuge cardinals}}$. \qedofLemma\qedskip

\begin{Thm} {\rm(Laver function for a super-$C^{(n)}$-hyperhuge cardinal)}
\Label{p-Lg-RcA-3}
\imemox{Laver function for ...hyperhuge. }  
Suppose that $\mu$ is an inaccessible cardinal and 
$\kappa$ is super-$C^{(\infty)}$-hyperhuge in $V_\mu$. Then there is
$\mapping{f}{\kappa}{V_\kappa}$ \st\ for any $n\in\omega$, $a\in V_\mu$, and
$\lambda_0\geq\kappa$, there are
$\lambda_0<\lambda<\kappa^*<\lambda^*<\mu$ \st\
$\kappa<\lambda$, $\cardof{\trcl(a)}\leq\lambda^*$, 
\begin{xitemize}
\xitem[x-Lg-RcA-7-a] 
  $V_\mu\models[\kappa,\lambda]^{[\kappa^*,\lambda^*,n]}$,  
\end{xitemize}
and there is a ultrafilter $U\subseteq\psof{\psof{\lambda^*}}$ witnessing 
\xitemof{x-Lg-RcA-7-a} \st\ $j_U(f)(\kappa)=a$ where $j_U$ is given by 
\xitemof{x-Lg-RcA-4-4}.\footnotemark 
\end{Thm}\memo{Scan\_2023-06-16--18.56-resurrection.pdf pp.66--67}
\prf The proof is just an adaptation of the proof of Theorem 20.21 in 
\cite{millennium-book}.
\footnotetext{In particular, we have $j_U(\kappa)=\kappa^*$, $j_U(\lambda)=\lambda^*$,
  $V_\lambda\prec_{\Sigma_n}\uniV$ and $V_{\lambda^*}$ ($=V_{j(\lambda)}$) $\prec_{\Sigma_n}\uniV$.} 

Suppose toward a contradiction, that the Theorem does not hold. Then for each
$\mapping{f}{\kappa}{V_\kappa}$, there are $a_f$, $n_f\in\omega$ and
$\kappa<\lambda_f<\kappa^*_f<\lambda^*_f<\mu$ \st\ 
\begin{xitemize}
\xitem[x-Lg-RcA-7-a-0] $\cardof{\trcl(a_f)}\leq\lambda^*_f$ ($\Leftrightarrow$\ \
  $a_f\in\calH((\lambda^*_f)^+)$);
\xitem[x-Lg-RcA-7-a-1] $V_\mu\models[\kappa,\lambda_f]^{[\kappa^*_f,\lambda^*_f,n_f]}$; but
\xitem[x-Lg-RcA-7-a-2] there is no witness $U$ of \xitemof{x-Lg-RcA-7-a-1} with
  $j_U(f)(\kappa)=a_f$. 
\end{xitemize}
We assume that, for each $\mapping{f}{\kappa}{V_\kappa}$, 
$\pairof{\lambda_f,\kappa^*_f,\lambda^*_f,n_f}$ is chosen to be the 
minimal possible (\wrt\ the lexicographical ordering) 
among those which satisfy 
\xitemof{x-Lg-RcA-7-a-0} $\sim$ \xitemof{x-Lg-RcA-7-a-2} together with $f$ and some $a_f$. 


For a $\alpha\in\On$, and $\mapping{g}{\alpha}{V_\alpha}$, let $(\ast)_{\alpha,g}$ be the 
assertion
\begin{xitemize}
\xitem[x-Lg-RcA-7-a-3] $\alpha$ is a cardinal and there are
  $\alpha<\delta<\alpha^*<\delta^*<\mu$, $n\in\omega$ and $a$ \st\ 
  \begin{xitemize}
  \xxitem[x-Lg-RcA-7-a-3][a] $V_\mu\models[\alpha,\delta]^{[\alpha^*,\delta^*,n]}$,\qquad
    $a\in\calH((\delta^*)^+)$
  \end{xitemize}
  but there is no witness $U$ of \xxitemof{x-Lg-RcA-7-a-3}{a} \st\ $j_U(g)(\alpha)=a$.
\end{xitemize}

Let $\mapping{f^\dagger}{\kappa}{V_\kappa}$ be defined recursively \st\ 
for $\alpha<\kappa$, if $(\ast)_{\alpha, f^\dagger\restr\alpha}$ holds, then 
$f^\dagger(\alpha)$  witnesses $(\ast)_{\alpha, f^\dagger\restr\alpha}$ as $a$ in 
\xitemof{x-Lg-RcA-7-a-3} together with $\pairof{\delta,\alpha^*,\delta^*,n}$ which is chosen 
to be minimal possible (\wrt\ the lexicographic ordering). Otherwise, we let
$f(\alpha)=\emptyset$.

Let 
\begin{xitemize}
\xitem[x-Lg-RcA-7-a-4]  $a^\dagger:= j_{U_{f^\dagger}}(f^\dagger)(\kappa)$
\end{xitemize}
where $U_{f^\dagger}$ is a witness of 
$[\kappa,\lambda_{f^\dagger}]^{[\kappa^*_{f^\dagger},\lambda^*_{f^\dagger}, n_{f^\dagger}]}$.
By definition of $f^\dagger$, elementarity and the assumption of this indirect proof, 
$a^\dagger$ together with $\kappa^*_{f^\dagger},\lambda^*_{f^\dagger}, n_{f^\dagger}$ is a 
witness of  
$(\ast)_{\alpha,
\ubecause{\scriptstyle j_{U_{\scriptscriptstyle f^\dagger}}(f^\dagger)\restr\kappa}%
         {}{$=f^\dagger$}}$. 
But the existence of the ultrafilter $U_{f^\dagger}$ in \xitemof{x-Lg-RcA-7-a-4} is a 
contradiction to this. \qedofThm 
\qedskip

In \cite{sfetal-II} and further in \cite{sfetal-III}, \cite{FuOt}, 
\cite{fuchino-sakai}, we studied in connection with Laver genericity, only 
classes $\calP$  
of \pos\ which are all stationary preserving. This is because the 
existence of a Laver-gen.\ large cardinal for one of such classes of \pos\ naturally 
extends the known reflection properties down to $\LT\kappa_\refl$, and strong versions of 
forcing axioms.  

However, we can also consider the  
class of all \pos\ in connection with Laver genericity. \assertof{5} in the following \Thmof{p-Lg-RcA-4} is 
such an instance of Laver genericity. We shall go more into this scenario in 
\sectionof{bedrock-Lg}.  

\memox{This must be rewritten!!}
\begin{Thm}
  \Label{p-Lg-RcA-4} \wassertof{1} Suppose that $\mu$ is inaccessible and 
  $\kappa<\mu$ is super-$C^{(\infty)}$-ultrahuge in $V_\mu$. Let 
  $\poP=\Col(\aleph_1,\kappa)$. Then,  
  in $V_\mu[\genG]$, for any ${V_\mu,\poP}$-generic $\genG$, $\aleph_2^{V_\mu[\genG]}$
  ($=\kappa$) is tightly super-$C^{(\infty)}$-$\sigma$-closed-Laver-gen.\ ultrahuge
    and $\CH$ holds.\smallskip

  \wassert{2}
  Suppose that $\mu$ is inaccessible and $\kappa<\mu$ is 
  super-$C^{(\infty)}$-ultrahuge with a 
  Laver function $\mapping{f}{\kappa}{V_\kappa}$ for super-$C^{(\infty)}$-ultrahugeness in
  $V_\mu$.   
  If $\poP$ is the CS-iteration of length $\kappa$ for forcing \PFA\ along with $f$, then,  
  in $V_\mu[\genG]$ for any $(V_\mu,\poP)$-generic $\genG$, $\aleph_2^{V_\mu[\genG]}$
  ($=\kappa$)  
  is tightly super-$C^{(\infty)}$-proper-Laver-gen.\ ultrahuge and 
  $\continuum=\aleph_2$ holds.\smallskip

  \wassert{2\rlap{$'$}\,}
  Suppose that $\mu$ is inaccessible and $\kappa<\mu$ is super-$C^{(\infty)}$-ultrahuge with a
  Laver function $\mapping{f}{\kappa}{V_\kappa}$ for super-$C^{(\infty)}$-ultrahugeness in
  $V_\mu$.  
  If $\poP$ is the RCS-iteration of length $\kappa$ 
  for forcing \MM\ along with $f$, then,   
  in $V_\mu[\genG]$ for any $(V_\mu,\poP)$-generic $\genG$, $\aleph_2^{V_\mu[\genG]}$
  ($=\kappa$)  
  is tightly super-$C^{(\infty)}$-semi-proper-Laver-gen.\ ultrahuge and 
  $\continuum=\aleph_2$ holds.\smallskip\memox{\assertof{3}, should be 
    also rewritten to assertions in the form of $V_\mu\models$ ... }

  \wassert{3}
  Suppose that $\mu$ is inaccessible and $\kappa$ is super-$C^{(\infty)}$-ultrahuge with a 
  Laver function $\mapping{f}{\kappa}{V_\kappa}$ for super-$C^{(\infty)}$-ultrahugeness in
  $V_\mu$. If 
  $\poP$ is a FS-iteration of length $\kappa$ for  
  forcing \MA\  
  along with $f$,   
  then, in $V_\mu[\genG]$ for any $(V_\mu,\poP)$-generic $\genG$, $\continuum$ ($=\kappa$) is 
  tightly super-$C^{(\infty)}$-c.c.c.-Laver-gen.\ ultrahuge, and $\continuum$ is 
  very large in $V_\mu[\genG]$.

  \wassert{4}
  Suppose that $\mu$ is inaccessible and $\kappa$ is super-$C^{(\infty)}$-ultrahuge with a 
  Laver function $\mapping{f}{\kappa}{V_\kappa}$ for super-$C^{(\infty)}$-ultrahugeness in
  $V_\mu$. If 
  $\poP$ is a FS-iteration of length $\kappa$ 
  along with $f$ enumerating ``all'' \pos,   
  then, in $V_\mu[\genG]$ for 
  any $(V_\mu,\poP)$-generic $\genG$, $\continuum$ ($=\aleph_1$) is  
  tightly super-$C^{(\infty)}$-all \pos-Laver-gen.\ ultrahuge, and \CH\ holds.\footnotemark
\end{Thm}
\footnotetext{Actually we can obtain a model with the desired property simply by Levy 
  collapsing $\kappa$ to $\omega_1$. We just chose this narrative to emphasize the parallelism 
  to the cases \assertof{2}, \assertof{2'}  and \assertof{3}.} 

\prf The proof can be done similarly to that of Theorem 5.2 in \cite{sfetal-II} using 
\Lemmaof{p-Lg-RcA-4-0} below. In the following we shall only check the case \assertof{4}.

Suppose that $\mapping{f}{\kappa}{V_\kappa}$ is a super-$C^{(\infty)}$-hyperhuge Laver 
function as in \Thmof{p-Lg-RcA-3}. 

Let $\seqof{\poP_\alpha,\utpoQ_\beta}{\alpha\leq\kappa,\beta<\kappa}$ be a FS-iteration 
defined by
\begin{xitemize}
\xitem[p-Lg-RcA-4-a-0] $\utpoQ_\beta={}
  \left\{\,\begin{array}{@{}ll}
    f(\alpha), &\mbox{ if }f(\alpha)\mbox{ is a }\poP_\beta\mbox{-name of a \po};\\[\jot]
    \ssetof{\bbone}, &\mbox{otherwise}
  \end{array}\right.
$
\end{xitemize}
for $\beta<\kappa$.

Let $\genG$ be a $(V_\mu,\poP_\kappa)$-generic filter. Clearly
$V_\mu[\genG]\modelof{2^{\aleph_0}=\kappa=\aleph_1}$. We show that $\kappa$ is tightly 
super-$C^{(\infty)}$-all \pos-Laver-gen.\ ultrahuge in $V_\mu[\genG]$.

Suppose that $\poP$ is a \po\ in $V_\mu[\genG]$, $\kappa<\lambda_0$ and
$n\in\omega$. Let $n'>n$ be sufficiently large and let $\utpoP$ be a $\poP_\kappa$-name of
$\poP$. 

Working in $V_\mu$, we can find 
\begin{xitemize}
\xitem[x-Lg-RcA-6-1-a] 
  $\cardof{\poP}<\lambda<\kappa^*<\lambda^*$ and $j$, $M\subseteq\uniV$
\end{xitemize}
\st\
\begin{xitemize}
\xitem[x-Lg-RcA-6-1-0] 
  $\Elembed{j}{\uniV}{M}{\kappa}$,
\xitem[x-Lg-RcA-6-1-1] $j(\kappa)=\kappa^*$, $j(\lambda)=\lambda^*$,
\xitem[x-Lg-RcA-6-1-2] $\fnsp{\lambda^*}{M}\subseteq M$,
\xitem[x-Lg-RcA-6-1-3] $V_\lambda\prec_{\Sigma_{n'}}\uniV$,
  $V_{\lambda^*}\prec_{\Sigma_{n'}}\uniV$, and 
\xitem[x-Lg-RcA-6-1-4] $j(f)(\kappa)=\utpoP$ 
\end{xitemize}
by definition of $f$. 


By elementarity (and by the definition \xitemof{p-Lg-RcA-4-a-0} of $\poP_\kappa$),
\begin{xitemize}
\xitem[x-Lg-RcA-6-2] 
  $j(\poP_\kappa)\ \sim_{\poP_\kappa}\ (\poP_\kappa\ast\utpoP)\ast\utpoR$ 
\end{xitemize}
for a $(\poP_\kappa\ast\utpoP)$-name $\utpoR$ of a \po.
Note that $(\poP_\kappa\ast\utpoP\ast\utpoR)/\genG$ corresponds to a \po\ of the form
$\poP\ast\utpoQ$.

Let $\genH^*$ be $(\uniV,(\poP_\kappa\ast\utpoP)\ast\utpoR)$-generic filter 
with $\genG\subseteq\genH^*$. $\genH^*$ corresponds to a $(\uniV,j(\poP_\kappa))$-generic 
filter $\genH\supseteq\genG$ via the equivalence \xitemof{x-Lg-RcA-6-2}.

Let $\tilde{j}$ be defined by
\begin{xitemize}
\xitem[x-Lg-RcA-6-3] 
  $\mapping{\tilde{j}}{\uniV[\genG]}{M[\genH]}$;\quad $\uta^\genG\mapsto j(\uta)^\genH$
\end{xitemize}
for all $\poP_\kappa$-names $\uta$. 

A standard proof shows that $f$ is well-defined, and
$\Elembed{j}{\uniV[\genG]}{M[\genH]}{\kappa}$.   
By \xitemof{x-Lg-RcA-6-1-1} and \xitemof{x-Lg-RcA-6-1-2}, we have
$\tilde{j}\imageof \tilde{j}(\lambda)=j\imageof j(\lambda)\in M[\genH]$.  
Since $\genH\in M[\genH]$, the $(\uniV[\genG],\poP\ast\utpoQ)$-generic filter corresponding 
to $\genH$ is also in $M[\genH]$. 

By \xitemof{x-Lg-RcA-6-1-a}, \xitemof{x-Lg-RcA-6-1-3}, by the choice of $n'$,
and by 
\Lemmaof{p-Lg-RcA-4-0},\,\assertof{1}, we have\\ 
${V_\lambda}^{\uniV[\genG]}\prec_{\Sigma_n}\uniV[\genG]$\quad and\quad
${V_{\tilde{j}(\lambda)}}^{\uniV[\genH]}
={V_{\lambda^*}}^{\uniV[\genH]}\prec_{\Sigma_n}\uniV[\genH]$.

Since $\poP$ and $n$ were arbitrary, this shows that $\kappa$ is tightly 
super-$C^{(\infty)}$-all \pos-Laver-gen.\ ultrahuge in $V_\mu[\genG]$.
\qedofThm\qedskip

\ifprivate\addcontentsline{toc}{section}{** Vλ≺n V 〜 Vλ[G]≺n' V[G]}\fi
\begin{Lemma}
  \Label{p-Lg-RcA-4-0} \wassertof{1} For a natural number $n$, if $n'>n$ is sufficiently 
  large and\smallskip\\ 
  \ixitem[x-Lg-RcA-7-0] $V_\lambda\prec_{\Sigma_{n'}}\uniV$,\smallskip\\
  then we have
  $V_\lambda[\genG]\prec_{\Sigma_n}\uniV[\genG]$ for any \po\ $\poP\in V_\lambda$ and
  $(\uniV,\poP)$-generic $\genG$. \smallskip

  \wassert{2} For a natural number $n$, if $n'>n$ is sufficiently large and\smallskip\\
  \ixitem[x-Lg-RcA-7-1] $V_\lambda[\genG]\prec_{\Sigma_{n'}}\uniV[\genG]$ for some \po\
  $\poP\in V_\lambda$ and $(\uniV,\poP)$-generic $\genG$,\smallskip\\ then we have
  $V_\lambda\prec_{\Sigma_n}\uniV$. \ifextended\else\qed\fi
\end{Lemma}
{\ifextended\extendedcolor
\prf \assertof{1}: 
Suppose that $\overline{a}\in V_\lambda[\genG]$ and $\varphi=\varphi(\overline{x})$ is 
a $\Sigma_n$-formula. 
There are $\poP$-names 
$\overline{\uta}\in V_\lambda$ \st\ $\overline{a}=\overline{\uta}[\genG]$.

If $V_\lambda[\genG]\models\varphi(\overline{a})$, there is $\condp\in\genG$ \st\
$V_\lambda\models\condp\forces{\poP}{\varphi(\overline{\uta})}$. By \xitemof{x-Lg-RcA-7-0} it 
follows that $\uniV\models\condp\forces{\poP}{\varphi(\overline{\uta})}$. 
Thus $\uniV[\genG]\models\varphi(\overline{a})$. 

The same argument also applies to $\neg\varphi$. \smallskip

\assertof{2}: We use the $\Lin$-formula $\Phi(x,y)$ of \Thmof{p-intro-1}. 
By \xitemof{x-Lg-RcA-7-1}, there is $r\in V_\lambda[\genG]$ \st\ 
\begin{xitemize}
\item[] $V_\lambda=\Phi(\cdot,r)^{V_\lambda[\genG]}
  =\Phi(\cdot,r)^{\uniV[\genG]}\cap V_\lambda[\genG]\subseteq\Phi(\cdot,r)^{\uniV[\genG]}$.
\end{xitemize}
For any $\Sigma_n$-formula $\varphi(\overline{x})$ and
$\overline{a}\in\Phi(\cdot,r)^{V_\lambda[\genG]}$. Since $\varphi^{\Phi(\cdot,r)}$ is a
$\Sigma_{n'}$-formula (by the choice of $n'$), we have
\begin{xitemize}
\item[] $V_\lambda\models\varphi(\overline{a})$\ \ 
  $\Leftrightarrow$\ \  
$V_\lambda[\genG]\models\varphi^{\Phi(\cdot,r)}(\overline{a})$\ \ 
  \smash{$\obecause{\Leftrightarrow}{}{by \xitemof{x-Lg-RcA-7-1}}$}\ \  
  $\uniV[\genG]\models\varphi^{\Phi(\cdot,r)}(\overline{a})$\ \ 
  $\Leftrightarrow$\ \  
  $V\models\varphi(\overline{a})$. 
\end{xitemize}
This shows that $V_\lambda\prec_{\Sigma_n}\uniV$.
\qedofLemma
\qedskip\fi}

In \Thmof{p-Lg-RcA-4},\,\assertof{4}, the model constructed in the proof satisfies \CH. 
The next lemma suggests that this does 
not depend on the specific construction of the model given there.

In the following, we consider a strengthening of tightness of Laver genericity: we 
say a cardinal $\kappa$ is tightly$^+$ $\calP$-Laver-gen.\ hyperhuge/ultrahuge if the 
definition of 
the tightly $\calP$-Laver gen.\ 
hyperhugeness/ultrahugeness of $\kappa$ holds with the clause
``$\cardof{\poP\ast\utpoQ}\leq j(\kappa)$'' replaced by
``there is a complete Boolean algebra $\baB$ of size $j(\kappa)$ \st\
$\poP\ast\utpoQ\sim \baB^+$''. Note that this condition is equivalent to 
$2^{\LT\kappa}=\kappa$, $\cardof{\poP\ast\utpoQ}\leq\kappa$ and $\poP\ast\utpoQ$ has the
$j(\kappa)$-cc. 
Note  also
that all models constructed in the proof of \Thmof{p-Lg-RcA-4} actually satisfy this 
stronger notion of tightness of the (super-$C^{(n)}$) $\calP$-Laver-gen.\ hyperhugeness or 
ultrahugeness of $\kappa$. Later with much more work we show in \Corof{p-bedrock-11} that this 
stronger version of tightness can be dropped from the following Lemma. 

\begin{Lemma}\Label{p-Lg-RcA-4-1}
If $\calP$ is the class of all \pos\ and $\kappa$ is tightly$^+$ $\calP$-Laver-gen.\ 
hyperhuge cardinal then $\kappa=2^{\aleph_0}=\aleph_1$. 
\end{Lemma}
\prf $\kappa\leq 2^{\aleph_0}$ holds even without tightness: see e.g.\ Lemma 3.7 in 
\cite{future}. \memo{The number might have been changed: it is p-Lg-RA-1-2-0 in \cite{future}}

To show $2^{\aleph_0}\leq\kappa$, let $\lambda$ and $\poQ$ be \st\ 
\begin{xitemize}
\xitem[x-Lg-RcA-8-a-0] 
  $\lambda>2^{\aleph_0}$,
  $\kappa$ and $\lambda$ is large enough \st\ \SCH\ holds above some $\mu<\lambda$ (this is 
  possible by \Corof{p-bedrock-3} ,\,\assertof{2} in \sectionof{bedrock}, and we need here 
  the Laver gen.\ hyperhugeness of $\kappa$),
\xitem[x-Lg-RcA-8-a-1] 
  $\poQ$ is positive elements of a complete Boolean algebra, 
  and, 
\xitem[x-Lg-RcA-8-a-2] 
  for $(\uniV,\poQ)$-generic $\genH$, there are $j$, $M\subseteq\uniV[\genH]$ \st\
  \wassertof{1} $\Elembed{j}{\uniV}{M}{\kappa}$,
  \wassertof{2} $j(\kappa)>\lambda$, \wassertof{3} $\cardof{\poQ}\leq j(\kappa)$,
  and \wassertof{4} ${V_{j(\lambda)}}^{\uniV[\genH]}\subseteq M$. 
\end{xitemize}

By \xitemof{x-Lg-RcA-8-a-1}, 
each $\poQ$-name $\utilde{r}$ of a 
real 
corresponds to a mapping $\mapping{f}{\omega}{\poQ}$. 
By \xitemof{x-Lg-RcA-8-a-0} and by \xitemof{x-Lg-RcA-8-a-2}, 
\assertof{3}, 
there are at most $j(\kappa)$ many such mappings. Thus 
we have $\uniV[\genH]\modelof{\continuum\leq j(\kappa)}$, By \xitemof{x-Lg-RcA-8-a-2}, 
\assertof{4}, it follows 
$M\modelof{\continuum\leq j(\kappa)}$. By elementarity, it follows that
$\uniV\modelof{\continuum\leq\kappa}$.\smallskip

$\kappa=\aleph_1$: Otherwise $\aleph_1<\kappa$. In $\uniV$, let $\poP$ be the standard \po\ 
collapsing $\aleph_1$ to be countable. Let $\utpoQ$ be a $\poP$-name \st\ for 
$(\uniV,\poP\ast\utpoQ)$-generic $\genH$ there are $j$, $M\subseteq\uniV[\genH]$ \st\
$\Elembed{j}{\uniV}{M}{\kappa}$, and $\poP$, $\genH\in M$.

By $\genH\cap\poP\in M$, we have
$M\modelof{{\aleph_1}^\uniV=j({\aleph_1}^\uniV)\mbox{ is countable}}$. This is a 
contradiction to the elementarity of $j$. 
\qedofLemma\qedskip

Recall that, for an iterable $\calP$, $(\calP,\calH(\kappa))$-\RcAp\ holds if and 
only if $\MP(\calP,\calH(\kappa))$ holds (\Propof{p-intro-0},\,\assertof{1}).

\begin{Thm}
  \Label{p-Lg-RcA-5} Suppose that $\calP$ is an iterable class of \pos\ and $\kappa$ is 
  tightly super-$C^{(\infty)}$-$\calP$-Laver-gen.\ ultrahuge. Then 
  $(\calP,\calH(\kappa))$-\RcAp\ holds. 
\end{Thm}
\prf A modification of the proof of \Thmof{p-Lg-RcA-0} works. 

Suppose that $\kappa$ is tightly super-$C^{(\infty)}$-$\calP$-Laver-gen.\ ultrahuge, 
$\poP\in\calP$, and $\forces{\poP}{\varphi(\overline{a}\checked)}$ for 
an $\Lin$-formula $\varphi$ and $\overline{a}\in\calH(\kappa)$. We want to show that
$\varphi(\overline{a})$ holds in some $\calP$-ground of $\uniV$.

Let $n$ be a sufficiently large natural number \st\ the following arguments go through. In 
particular, we assume that ${V_\alpha}^\uniV\prec_{\Sigma_n}\uniV$ implies that
``$\varphi(\overline{x})$'' and ``$\forces{\cdot}{\varphi(\overline{x}\checked)}$ are 
absolute between ${V_\alpha}^\uniV$ and $\uniV$, 
and ${V_\alpha}^\uniV\prec_{\Sigma_n}\uniV$ also implies that a sufficiently large fragment 
of \ZFC\ holds in $V_\alpha$. 

Let $\utpoQ$ be a $\poP$-name \st\ $\forces{\poP}{\utpoQ\in\calP}$ and, for
$(\uniV,\poP\ast\utpoQ)$-generic $\genH$, there are a $\lambda>\kappa$ 
with 
\begin{xitemize}
\xitem[x-Lg-RcA-9] 
  $V_\lambda\prec_{\Sigma_n}\uniV$, 
\end{xitemize}
and $j$, $M\subseteq \uniV[\genH]$ \st\ $\Elembed{j}{\uniV}{M}{\kappa}$,
$j(\kappa)>\lambda$, $\poP$, $\genH$, ${V_{j(\lambda)}}^{\uniV[\genH]}\in M$,
$\cardof{\poP\ast\utpoQ}\leq j(\kappa)$  ($<j(\kappa)$), and 
${V_{j(\lambda)}}^{\uniV[\genH]}\prec_{\Sigma_n}\uniV[\genH]$. 

Replacing $\poP\ast\utpoQ$ by an appropriate isomorphic \po\ (and replacing $\calH$ by 
corresponding filter), we may assume that $\poP\ast\utpoQ\in {V_{j(\lambda)}}^\uniV$. 

By the choice of $n$, we 
have $V_\lambda\models\forces{\poP}{\varphi(\overline{a}\checked)}$.
$j({V_\lambda}^\uniV)={V_{j(\lambda)}}^M\prec_{\Sigma_n}M$ by elementarity of $j$,   
and 
\begin{xitemize}
\xitem[x-Lg-RcA-9-0] 
  ${V_{j(\lambda)}}^M={V_{j(\lambda)}}^{\uniV[\genH]}$
\end{xitemize}
by the closedness of $M$. 
Since $V_\lambda\prec_{\Sigma_n}\uniV$, we have
$V_\lambda[\genH]\prec_{\Sigma_{n_0}}\uniV[\genH]$ for a still large 
enough $n_0\leq n$ by \Lemmaof{p-Lg-RcA-4-0},\,\assertof{1}. Since
${V_{j(\lambda)}}^{\uniV[\genH]}\prec_{\Sigma_n}\uniV[\genH]$, it  
follows 
that ${V_\lambda}^{\uniV[\genH]}\prec_{\Sigma_{n_0}}{V_{j(\lambda)}}^{\uniV[\genH]}$. 
Thus 
\begin{xitemize}
\xitem[x-Lg-RcA-10] 
  ${V_\lambda}^\uniV\prec_{\Sigma_{n_1}}{V_{j(\lambda)}}^\uniV$
\end{xitemize}
for a still large enough
$n_1\leq n_0$ by \Lemmaof{p-Lg-RcA-4-0},\,\assertof{2}. 

In particular, we have
${V_{j(\lambda)}}^\uniV\models\forces{\poP}{\varphi(\overline{a}\checked)}$,  and hence
$V_{j(\lambda)}[\genG]\models\varphi(\overline{a})$ where $\genG$ is the $\poP$-part of
$\genH$. Note that by \xitemof{x-Lg-RcA-9} and \xitemof{x-Lg-RcA-10}, $V_{j(\lambda)}$ 
satisfies a sufficiently large fragment of \ZFC. 

Thus we have
$V_{j(\lambda)}[\genH]\modelof{\,\mbox{there is a }
  \calP\mbox{-ground satisfying }\varphi(\overline{a})}$, and hence 
\begin{xitemize}
\item[] ${V_{j(\lambda)}}^{\uniV[\genH]}\modelof{\,\mbox{there is a }
  \calP\mbox{-ground satisfying }\varphi(\overline{a})}$
\end{xitemize}
by \Lemmaof{p-Lg-RcA-0-0}. By 
\xitemof{x-Lg-RcA-9-0} and elementarity, it follows that
\begin{xitemize}
\item[] ${V_{\lambda}}\modelof{\,\mbox{there is a }
  \calP\mbox{-ground satisfying }\varphi(\overline{a})}$.
\end{xitemize}
Finally, this implies ${\uniV}\modelof{\,\mbox{there is a }
  \calP\mbox{-ground satisfying }\varphi(\overline{a})}$ by \xitemof{x-Lg-RcA-9}.\\
\qedofThm

\section{Bedrock of a tightly generic hyperhuge cardinal}
\Label{bedrock}
Recall that a cardinal $\kappa$ is {\It hyperhuge}, if for every $\lambda>\kappa$, there is 
$\Elembed{j}{\uniV}{M\subseteq\uniV}{\kappa}$ \st\ $\lambda<j(\kappa)$ and
$\fnsp{j(\lambda)}{M}\subseteq M$. By \Lemmaof{p-Lg-RcA-1-0}, a hyperhuge cardinal $\kappa$ 
can be characterized in terms of existence of $\kappa$-complete normal ultrafilters with the 
property \xitemof{x-Lg-RcA-4-3}. 

For a class $\calP$ of \pos, 
a cardinal $\kappa$ is {\It tightly $\calP$-generic hyperhuge} (tightly $\calP$-gen.\ 
hyperhuge, for short) if for any 
$\lambda>\kappa$, there is $\poQ\in\calP$ \st\ for a $(\uniV,\poQ)$-generic $\genH$, there 
are $j$, $M\subseteq\uniV[\genH]$ \st\ $\Elembed{j}{\uniV}{M}{\kappa}$,
$\lambda<j(\kappa)$, 
$\cardof{\poQ}\leq j(\kappa)$, and
$j\imageof{j(\lambda)}, \genH\in M$. 

For a class $\calP$ of \pos, 
a cardinal $\kappa$ is {\It tightly $\calP$-Laver-generically hyperhuge} 
(tightly $\calP$-Laver-gen.\ hyperhuge, for short) if for any
$\lambda>\kappa$, and $\poP\in\calP$
there is a $\poP$-name $\utpoQ$ with $\forces{\poP}{\utpoQ\in\calP}$ \st\ for a
$(\uniV,\poP\ast\utpoQ)$-generic $\genH$, there  
are $j$, $M\subseteq\uniV[\genH]$ \st\ $\Elembed{j}{\uniV}{M}{\kappa}$,
$\lambda<j(\kappa)$, 
$\cardof{\poP\ast\utpoQ}\leq j(\kappa)$, and $j\imageof{j(\lambda)}, \genH\in M$. 

\memo{Scan\_2023-06-16--18.56-resurrection - annotated.pdf\\pp.83--84}

The following Lemma is easy to prove.
\begin{Lemma}
  \Label{p-bedrock-0} \wassertof{1} $\kappa$ is \hyperref[bedrock]{hyperhuge} if and only 
  if the following  
  holds: 
  \begin{xitemize}
  \xitem[x-bedrock-0] for every $\lambda>\kappa$, there is 
    $\Elembed{j}{\uniV}{M\subseteq\uniV}{\kappa}$ \st\ $\lambda<j(\kappa)$,\\
    ${V_{j(\lambda)}}\in M$, and $\fnsp{j(\lambda)}{M}\subseteq M$. 
  \end{xitemize}

  \wassert{2} If $\kappa$ is hyperhuge then it is \hyperref[ultrahuge]{ultrahuge}.\smallskip

  \wassert{3} Suppose that $\calP$ is a class of \pos\ \st\ the trivial \po\
  $\ssetof{\bbone}$ is in $\calP$.  
  If $\kappa$ is hyperhuge then $\kappa$ is tightly $\calP$-gen.\ hyperhuge.
  \smallskip 

  \wassert{4} For a class $\calP$ of \pos, $\kappa$ is tightly $\calP$-gen.\ hyperhuge if 
  and only if the following holds: 
  \begin{xitemize}
  \xitem[x-bedrock-1] for any
    $\lambda>\kappa$, there is $\poQ\in\calP$ \st\ for a $(\uniV,\poQ)$-generic $\genH$, 
    there  
    are $j$, $M\subseteq\uniV[\genH]$ \st\ $\Elembed{j}{\uniV}{M}{\kappa}$,
    $\lambda<j(\kappa)$, 
    $\cardof{\poQ}\leq j(\kappa)$, 
    and
     ${V_{j(\lambda)}}^{\uniV[\genH]}$, $j\imageof{j(\lambda)}, \genH\in M$. 
  \end{xitemize}

  \wassert{5} For a class $\calP$ of \pos, $\kappa$ is tightly $\calP$-Laver-gen.\ 
  hyperhuge if and only if the following holds: 
  \begin{xitemize}
  \xitem[x-bedrock-1-0] for any
    $\lambda>\kappa$ and $\poP\in\calP$, there is a $\poP$-name $\utpoQ$ with
    $\forces{\poP}{\utpoQ\in\calP}$ \st\ for a $(\uniV,\poP\ast\utpoQ)$-generic $\genH$,  
    there  
    are $j$, $M\subseteq\uniV[\genH]$ \st\ $\Elembed{j}{\uniV}{M}{\kappa}$,
    $\lambda<j(\kappa)$, 
    $\cardof{\poQ}\leq j(\kappa)$, 
    and
    ${V_{j(\lambda)}}^{\uniV[\genH]}$, $j\imageof{j(\lambda)}, \genH\in M$. 
  \end{xitemize}

  \wassert{6} For an iterable class $\calP$ of \pos, a 
  tightly $\calP$-Laver generically hyperhuge cardinal is 
  tightly $\calP$-gen.\ hyperhuge. 
  \smallskip

  \wassert{7} For a class $\calP$ of \pos, if $\kappa$ is (tightly 
  resp.) $\calP$-Laver-gen.\ hyperhuge, then it is (tightly resp.) $\calP$-Laver-gen.\ 
  ultrahuge. 
  \ifextended\else\qed\fi
\end{Lemma}
{\ifextended\extendedcolor\prf 
  \assertof{1}: It is clear that \xitemof{x-bedrock-0} implies that $\kappa$ is hyperhuge.

  Suppose that $\kappa$ is hyperhuge. We show that \xitemof{x-bedrock-0} holds. 
  For $\lambda>\kappa$, let $\lambda^*:=\cardof{V_\lambda}^+$ and
  $\Elembed{j}{\uniV}{M}{\kappa}$ be \st\ $j(\kappa)>\lambda^*$ and 
  \begin{xitemize}
  \xitem[x-bedrock-2] $\fnsp{j(\lambda^*)}{M}\subseteq M$.
  \end{xitemize}


  By elementarity, we have 
  \begin{xitemize}
  \xitem[x-bedrock-3] 
    $M\models j(\lambda^*)=\cardof{{V_{j(\lambda)}}^M}^+$.
  \end{xitemize}
  \begin{Claim}
    \Label{cl-bedrock-0} $M\ni{V_{j(\lambda)}}^M=V_{j(\lambda)}$. 
  \end{Claim}
  \noindent
  \prfofClaim We prove ${V_\alpha}^M=V_\alpha$ for all $\omega\leq\alpha\leq j(\lambda)$. 

  For $\alpha=\omega$, this clearly holds.

  Suppose that ${V_\alpha}^M=V_\alpha$ for $\omega\leq\alpha<j(\lambda)$. Then
  $j(\lambda^*)>\cardof{{V_\alpha}^M}=\cardof{V_\alpha}$ by \xitemof{x-bedrock-3}. 
  By \xitemof{x-bedrock-2}, we have $V_{\alpha+1}=\psof{V_\alpha}\in M$:  
  $\psof{V_\alpha}\subseteq M$ by \xitemof{x-bedrock-2}. Since $M\models\ZFC$ by 
    elementarity of $j$, it follows that $\psof{V_\alpha}=\psof{V_\alpha}^M\in M$.  

  If $\alpha\leq j(\lambda)$ is a limit and the Claim holds for $\beta<\alpha$, then
  $\seqof{V_\beta}{\beta<\alpha}\in\fnsp{\alpha}{M}\subseteq M$. Thus
  $V_\alpha=\bigcup_{\beta<\alpha}V_\beta\in M$. 
  \qedofClaim\qedskip

  $\fnsp{j(\lambda)}{M}\subseteq M$ by $j(\lambda)<j(\lambda^*)$  and \xitemof{x-bedrock-2}.
  Thus $j$ witnesses \xitemof{x-bedrock-0} for $\lambda$. 
  \smallskip

  \assertof{2}: follows from \assertof{1}. \assertof{3}: is trivial.\smallskip

  \assertof{4}: If $\kappa$ satisfies \xitemof{x-bedrock-1} then $\kappa$ is clearly 
  tightly $\calP$-gen.\ hyperhuge.

  Suppose that $\kappa$ is tightly $\calP$-gen.\ hyperhuge. 
  For $\lambda>\kappa$, let 
  $\lambda^*$ be \st\ 
  \begin{xitemize}
  \xitem[x-bedrock-4] $\lambda^*\geq\cardof{V_\lambda}^+$ and
  \xitem[x-bedrock-5] $V_{\lambda^*}\prec_{\Sigma_n}\uniV$ for a sufficiently large natural 
    number $n$. 
  \end{xitemize}
  By \xitemof{x-bedrock-1}, there is $\poQ\in\calP$ \st\ for $(\uniV,\poQ)$-generic $\genH$ 
  there are $j$, $M\subseteq\uniV[\genH]$ \st\ 
  \begin{xitemize}
    \xitem[x-bedrock-5-0] $\Elembed{j}{\uniV}{M}{\kappa}$,
  \xitem[x-bedrock-6] $\lambda^*<j(\kappa)$,
  \xitem[x-bedrock-7] $\cardof{\poQ}\leq j(\kappa)$, and
  \xitem[x-bedrock-8] $j\imageof{j(\lambda^*)}$, $\genH\in M$.
  \end{xitemize}
  By \xitemof{x-bedrock-8}, and Lemma 2.5, \assertof{5} in \cite{sfetal-II}, we have
  ${V_{j(\lambda^*)}}^\uniV\in M$. Hence ${V_{j(\lambda^*)}}^\uniV[\genH]\in M$ by 
  \xitemof{x-bedrock-8}. By \xitemof{x-bedrock-5} and \xitemof{x-bedrock-8},
  ${V_{j(\lambda^*)}}^\uniV[\genH]={V_{j(\lambda^*)}}^{\uniV[\genH]}$, and thus  
  ${V_{j(\lambda)}}^{\uniV[\genH]}\in M$.

  This shows that our $j$ and $M$ witnesses \xitemof{x-bedrock-1} for $\lambda$. \smallskip

  \assertof{5}: Similarly to \assertof{4}.\smallskip

  \assertof{6}: Since $\calP$ is iterable, the \po\ $\poP\ast\utpoQ$ in the definition of
  $\calP$-Laver-gen.\ hyperhugeness is an element of $\calP$. \smallskip

  \assertof{7}: By \assertof{5}. 
\qedofLemma\qedskip\fi}

For a cardinal $\kappa$, a ground $\uniW$ of the universe $\uniV$ is called 
a {\It$\LE\kappa$-ground\/} if there is a \po\ $\poP\in\uniW$ of cardinality $\LE\kappa$ 
(in the sense of $\uniV$) and $(\uniW,\poP)$-generic filter $\genG$ 
\st\ $\uniV=\uniW[\genG]$.  Note that this definition of $\LE\kappa$-ground diverges from the 
convention of ``$\calP$-ground'' in that ``of cardinality $\LE\kappa$'' is meant here not in the 
sense of $\uniW$ but 
rather of $\uniV$. In \Propof{p-bedrock-9} we will show that in our context, this actually 
implies ``$\LE\kappa$-ground'' in line with the definition of $\calP$-ground. 

Let 
\begin{xitemize}
\xitem[x-bedrock-9] $\overline{\uniW}:=\bigcap\setof{\uniW}{\uniW
  \mbox{ is a }\LE\kappa\mbox{-ground}}$. 
\end{xitemize}
Since there are only set many $\LE\kappa$-grounds, $\overline{\uniW}$ contains a ground by 
Theorem 1.3 in \cite{usuba}. 
We shall call $\overline{\uniW}$ the {\It$\LE\kappa$-mantle of\/ 
  $\uniV$}. 

By \Lemmaof{p-bedrock-0}, \assertof{3}, the following theorem generalizes 
Theorem 1.6 in \cite{usuba}. 

\begin{Thm}
  \Label{p-bedrock-1} Suppose that $\calP$ is a class of \pos.
  If $\kappa$ is a tightly $\calP$-gen.\ hyperhuge cardinal, then 
  the $\LE\kappa$-mantle is the smallest ground of\/ $\uniV$ (i.e.\ it is the bedrock of\/
  $\uniV$) and it is also a $\LE\kappa$-ground. 
\end{Thm}
\prf Suppose that $\kappa$ is tightly $\calP$-gen.\ hyperhuge and let $\overline{\uniW}$ 
be the $\LE\kappa$-mantle for this $\kappa$. 

By Theorem 1.3 in \cite{usuba} mentioned above, 
it is enough to show that, for any ground $\uniW\subseteq\overline{\uniW}$ is actually 
a $\LE\kappa$-ground and hence $\uniW=\overline{\uniW}$ holds. 

Let $\uniW\subseteq\overline{\uniW}$ be a ground. Let $\mu$ be the cardinality (in the 
sense of $\uniV$) of a \po\ $\poS\in\uniW$ \st\ there is a $(\uniW,\poS)$-generic $\genF$ 
\st\ $\uniV=\uniW[\genF]$. \Wolog, $\mu\geq\kappa$. 

By \Thmof{p-intro-1}, there is $r\in\uniV$ \st\ $\uniW=\Phi(\cdot,r)^\uniV$. 

Let $\theta\geq\mu$ be \st\ $r\in V_\theta$, and for a sufficiently large natural 
number $n$, we have 
\begin{xitemize}
\xitem[x-bedrock-9-1] ${V_\theta}^\uniV\prec_{\Sigma_n}\uniV$.
\end{xitemize}
By the choice of $\theta$, we have 
\begin{xitemize}
\xitem[x-bedrock-9-2] 
  $\Phi(\cdot,r)^{{V_\theta}^\uniV}=\Phi(\cdot,r)^\uniV\cap {V_\theta}^\uniV
  =\uniW\cap {V_\theta}^\uniV={V_\theta}^\uniW$. 
\end{xitemize}

There is a $\poQ\in\calP$ \st\ for $(\uniV,\poQ)$-generic $\genH$, there are $j$,
$M\subseteq\uniV[\genH]$ with 
\begin{xitemize}
\xitem[x-bedrock-10] 
  $\Elembed{j}{\uniV}{M}{\kappa}$, $\theta<j(\kappa)$,
\xitem[x-bedrock-11] 
  $\cardof{\poQ}\leq j(\kappa)$,
\xitem[x-bedrock-12] 
  ${V_{j(\theta)}}^{\uniV[\genH]}\subseteq M$, and 
\xitem[x-bedrock-13] 
  $\genH$, $j\imageof{j(\theta)}\in M$
\end{xitemize}
by \Lemmaof{p-bedrock-0},\,\assertof{4}. 

Since ${V_\theta}^\uniV\equiv {V_{j(\theta)}}^M$ by elementarity, \xitemof{x-bedrock-9-1} 
implies that ${V_{j(\theta)}}^M$ ($={V_{j(\theta)}}^{\uniV[\genH]}$) satisfies a 
sufficiently large fragment of \ZFC. Thus it follows that 
\begin{xitemize}
\xitem[x-bedrock-16] 
  ${V_{j(\theta)}}^\uniV$ satisfies a sufficiently large fragment of \ZFC\ and also
  ${V_{j(\theta)}}^\uniW$ satisfies still a sufficiently large fragment of \ZFC. 
\end{xitemize}
Thus we have 
\begin{xitemize}
\xitem[x-bedrock-14] 
  ${V_{j(\theta)}}^M\ubecause{=}{}{by \xitemof{x-bedrock-12}}
  {V_{j(\theta)}}^{\uniV[\genH]}\obecause{=}{}{by \Lemmaof{p-Lg-RcA-0-0} and the remark above}
  {V_{j(\theta)}}^\uniV[\genH]
  \ubecause{=}{}{by \Lemmaof{p-Lg-RcA-0-0} and the remark above}
{V_{j(\theta)}}^\uniW[\genF][\genH]$.
\end{xitemize}

Since 
\begin{xitemize}
\xitem[x-bedrock-16-0] 
  ${V_{j(\theta)}}^M\prec_{\Sigma_n}M$
\end{xitemize}
by elementarity and $V$ is a ground of $M$ by 
Grigorieff's Theorem, it follows from \Lemmaof{p-Lg-RcA-4-0},\,\assertof{2} that 
\begin{xitemize}
\xitem[x-bedrock-15] 
  ${V_{j(\theta)}}^\uniV\prec_{\Sigma_{n_0}}\uniV$ for a still large enough $n_0\leq n$. 
\end{xitemize}
\memo{sub-subscript の $n_0$ を調節する．}

\begin{Claim}\Label{cl-bedrock-1}
  \wassertof{1} ${V_{j(\theta)}}^\uniV$ is a generic extension of ${V_{j(\theta)}}^\uniW$ 
  by a \po\ of size $\LE j(\kappa)$.\smallskip

  \wassert{2} $V[\genH]$ is a generic extension of $\uniW$ by a \po\ of size
  $\LE j(\kappa)$.\smallskip

  \wassert{3} ${V_{j(\theta)}}^M$ is a generic extension of ${V_{j(\theta)}}^\uniW$ by a 
  \po\ of size $\LE j(\kappa)$. \smallskip

  \wassert{4} $M$ is a generic extension of $j(W):=\Phi(\cdot, j(r))^M$ by a \po\ of size
  $\LE j(\mu)$. 
\end{Claim}
\noindent
\prfofClaim \assertof{1}: $\uniV$ is a generic extension of $W$ by a \po\ of size
$\mu<\lambda<j(\kappa)$. Thus the assertion follows from \xitemof{x-bedrock-15}.
\smallskip

\assertof{2}: We have $\uniV[\genH]=\uniW[\genF][\genH]$ and
$\cardof{\poS\ast\poQ}\leq j(\kappa)$ by $\cardof{\poS}\leq\mu<\lambda<j(\kappa)$ and by 
\xitemof{x-bedrock-11}. \smallskip

\assertof{3}: By \xitemof{x-bedrock-14} and \assertof{1}.\smallskip

\assertof{4}: By definition of $\mu$ and elementarity of $j$. 
\qedofClaim\qedskip

\begin{Claim}
  \Label{cl-bedrock-1-0}
  ${V_{j(\theta)}}^{j(\uniW)}=j({V_\theta}^\uniW)
  \subseteq j({V_\theta}^{\overline{\uniW}})\subseteq {V_{j(\theta)}}^\uniW
  \subseteq {V_{j(\theta)}}^\uniV$. 
\end{Claim}
\noindent
\prfofClaim
By \xitemof{x-bedrock-12}, we have
${V_{j(\theta)}}^M={V_{j(\theta)}}^{\uniV[\genH]}$. By 
\Claimof{cl-bedrock-1},\,\assertof{3}, it follows that 
\begin{xitemize}
\xitem[x-bedrock-16-1] 
  ${V_{j(\theta)}}^\uniW$ is a
  $\LE j(\kappa)$-ground of ${V_{j(\theta)}}^{M}$. 
\end{xitemize}
By elementarity of $j$ and \xitemof{x-bedrock-16-0}, 
\begin{xitemize}
\xitem[x-bedrock-18]   $j({V_\theta}^{\overline{\uniW}})$ is a
  $\LE j(\kappa)$-mantle of $V_{j(\theta)}^M$. 
\end{xitemize}
Thus we have
\begin{xitemize}
\item[] 
  \smash{${V_{j(\theta)}}^{j(\uniW)}=j({V_\theta}^\uniW)
  \subseteq j({V_\theta}^{\overline{\uniW}})\obecause{\subseteq}{}{\qquad\qquad by
    \xitemof{x-bedrock-16-1}, \xitemof{x-bedrock-18} and by the definition of
    $\LE j(\kappa)$-mantle } {V_{j(\theta)}}^\uniW
  \subseteq {V_{j(\theta)}}^\uniV$}.\qedofClaim\qedskip
\end{xitemize}

The proof of the following claim poses the greatest technical challenge of the entire proof 
of the present theorem. 

The main ingredient of this proof is Solovay's trick with a stationary partition.

\begin{Claim}
  \Label{cl-bedrock-2} $j\imageof\lambda\in j(\uniW)$ for every $\lambda<j(\theta)$. 
\end{Claim}
\noindent
\prfofClaim Since $\theta$ is a limit cardinal, it is enough to show the claim for all 
regular  
$\lambda<j(\theta)$ with $\lambda>j(\mu)$. Let $\lambda$ be one of such cardinals, and let 
$\seqof{S_\alpha}{\alpha<\lambda}$ be a partition of $(E^\lambda_\omega)^\uniW$ into 
stationary sets in $\uniW$. Since $\cardof{\poS}\leq\mu<\lambda$ each $S_\alpha$ is 
stationary subset  
of $\lambda$ in $\uniV$. Since $\uniV[\genH]$ is a generic extension of $\uniV$ by $\poQ$ 
of size $\leq j(\kappa)$, and $\lambda>j(\mu)\geq j(\kappa)$, each $S_\alpha$, 
$\alpha<\lambda$ is a stationary subset of $\lambda$ in $\uniV[\genH]$.
Let
\begin{xitemize}
\xitem[x-bedrock-17]
  $\seqof{S'_\alpha}{\alpha<j(\lambda)}:=j(\seqof{S_\alpha}{\alpha<\lambda})\in j(\uniW)$. 
\end{xitemize}
\begin{Subclaim}\Label{sub-bedrock-a}\mbox{}\\
  \mbox{}\hfill$j\imageof\lambda=\setof{\beta<\sup(j\imageof\lambda)}{
    S'_\beta\cap \sup(j\imageof\lambda)
      \mbox{ is stationary subset of }\sup(j\imageof\lambda)\mbox{ in }M}.
    $
\end{Subclaim}
\prfofClaim  Suppose $\alpha\in\lambda$. To show that
$S'_{j(\alpha)}\cap\sup(j\imageof{\lambda})=j(S_\alpha)\cap\sup(j\imageof{\lambda})$ is 
a stationary subset of $\sup(j\imageof{\lambda})$ in $M$, suppose 
$C\subseteq\sup(j\imageof{\lambda})$ is a club subset of $\sup(j\imageof{\lambda})$ 
in $M$. We have to show that $C$ intersects with $j(S_\alpha)$.

Since $\cardof{\poQ}\leq j(\kappa)\leq j(\mu)<\lambda$ by the choice of $\lambda$, and 
$M\subseteq\uniV[\genH]$, there is an unbounded $D\subseteq\lambda$, $D\in\uniW$ \st\ for 
any $\xi_0$, $\xi_1\in D$ with $\xi_0<\xi_1$, $[j(\xi_0),j(\xi_1)]\cap C\not=\emptyset$. 
Since $S_\alpha$ is stationary, there is $\eta\in S_\alpha\cap\lim(D)$. Since
$\cf(\eta)=\omega$, $j(S_\alpha)\ni j(\eta)=\sup(j\imageof{\eta})$ the right-most side of 
this is an element of $C$ by definition of $D$. This shows that
$j(S_\alpha)\cap C\not=\emptyset$.

Suppose now that $S'_\beta\cap\sup(j\imageof{\lambda})$ is stationary subset of
$\sup(j\imageof{\lambda})$ in $M$ for some $\beta<\sup(j\imageof{\lambda})$. We want to 
show that $\beta=j(\alpha)$ for some $\alpha<\lambda$.

Since
$j\imageof{\lambda}=j\imageof{j(\theta)}\cap\sup(j\imageof{\lambda})\in M$ by 
\xitemof{x-bedrock-13}, and thus, also
$\setof{\pairof{\alpha,j(\alpha)}}{\alpha<\lambda}\in M$, we have
$M\models\cf(j\imageof{\lambda})=\lambda$, and hence also
$\uniV[\genH]\models\cf(j\imageof{\lambda})=\lambda$.

By \xitemof{x-bedrock-12}, $\psof{\lambda}^M=\psof{\lambda}^{\uniV[\genH]}$. Hence 
stationarity of a subset of $\sup(j\imageof{\lambda})$ is absolute between $M$ 
and $V[\genH]$. Thus, our assumption implies that $S'_\beta\cap\sup(j\imageof{\lambda})$ is 
stationary subset of $\sup(j\imageof{\lambda})$ in $\uniV[\genH]$.

By \xitemof{x-bedrock-9-1} and elementarity of $j$, $j({V_\theta}^{\overline{\uniW}})$ is a
$\LE j(\kappa)$-mantle of ${V_{j(\theta)}}^M$ ($={V_{j(\theta)}}^{V[\genH]}$ by 
\xitemof{x-bedrock-12}).

\begin{Subsubclaim}
  \Label{subsub-bedrock-0} For every $\eta<\lambda$, if $\cf^{j(\uniW)}(\eta)=\omega$ then 
  $\cf^\uniV(\eta)=\omega$.
\end{Subsubclaim}
\prfofClaim 
By \Claimof{cl-bedrock-1-0}. 
\qedofSubsubclaim
\qedskip

\begin{Subsubclaim}\Label{subsub-bedrock-1}
  There is an unbounded $E\subseteq\lambda$ \st\ $j\restr E\in j(\uniW)$. 
\end{Subsubclaim}
\prfofClaim By \Claimof{cl-bedrock-1},\,\assertof{4} and $j(\mu)<\lambda$ by choice of
$\lambda$, $M$ is a generic extension of $j(\uniW)$ by the $\lambda$-c.c.\ 
\po\ $j(\poS)$.\footnote{Since $\cardof{j(\poS)}\leq j(\mu)<\lambda$ in $M$, $j(\poS)$ has 
  the $\lambda$-c.c.\ in $M$ and hence it has the $\lambda$-c.c.\ also in
  $j(\uniW)\subseteq M$.} 
Since $j\imageof{\lambda}\in M$, there is a $j(\poS)$-name $\utilde{j}$ 
for $j\restr\lambda$.

Working in $j(\uniW)$, suppose that $\gamma<\lambda$. By the $\lambda$-c.c.\ of $j(\poS)$, 
we can find an  
increasing sequences $\seqof{\alpha_n}{n<\omega}$,
$\seqof{\beta_n}{n\in\omega}\in j(\uniW)$ \st\ $\gamma<\alpha_0$, $\alpha_n<\lambda$ for all
$n\in\omega$ and $\forces{j(\poS)}{\utilde{j}(\alpha_n)<\beta_n<\utilde{j}(\alpha_{n+1})}$ 
for all $n\in\omega$.

Let $\alpha:=\sup_{n\in\omega}\alpha_n$ and
$\beta:=\sup_{n\in\omega}\beta_n$. Then $\cf^{j(\uniW)}(\alpha)=\omega$ and hence
$\cf^\uniV(\alpha)=\omega$ by \Subsubclaimof{subsub-bedrock-0}. 
Thus
$\forces{j(\poS)}{\utilde{j}(\alpha)=\sup_{n\in\omega}\utilde{j}(\alpha_n)=\sup_{n\in\omega}\beta_n=\beta}$.

This shows that 
\begin{xitemize}
\item[] $E=\setof{\alpha<\lambda}{\mbox{there is }\beta\mbox{ \st\ }
  \forces{j(\poS)}{\utilde{j}(\check{\alpha})=\check{\beta}}}\in j(\uniW)$
\end{xitemize}
is a cofinal subset of $\lambda$ and $j\restr E\in j(\uniW)$. 
\qedofSubsubclaim
\qedskip

Now returning to the proof of the second-half of the \Subclaimof{sub-bedrock-a}, 
since $S'_\beta\cap\sup(j\imageof{\lambda})$ is stationary and $\lim(j\imageof{E})$ is club 
in $\sup(j\imageof{\lambda})$, there is $\eta\in S'_\beta\cap \lim(j\imageof{E})$.
Then $\cf^{j(\uniW)}(\eta)=\omega=\cf^\uniV(\eta)$ (the last equality by 
\Subsubclaimof{subsub-bedrock-0}). Let $\zeta<\lambda$ be minimal with 
$\eta\leq j(\zeta)$.

We have $\sup(j\imageof{\zeta})=\sup(j\imageof{(\zeta\cap E)})=\eta\leq j(\zeta)$. Since the 
cofinality of $\eta$ is $\omega$ and $j\restr(\eta\cap E)$  is increasing and 
$\eta\cap E$ is cofinal in $\eta$, we have $j(\zeta)=\sup(j\imageof{\zeta})=\eta$, and 
hence $j(\zeta)\in S'_\beta$. By elementarity, it follows that $\zeta\in S_\alpha$ for some
$\alpha<\lambda$. Then $j(\zeta)\in S'_{j(\alpha)}\cap S'_\beta$. Since $S'_\xi$'s are 
pairwise disjoint, it follows that $j(\alpha)=\beta$. Thus $\beta\in j\imageof{\lambda}$ as 
desired. 
\qedofSubclaim
\qedskip

By \Subclaimof{sub-bedrock-a} we have 
\begin{xitemize}
\item[] 
  $j\imageof\lambda=\setof{\beta<\sup(j\imageof\lambda)}{{}\begin{array}[t]{@{}l}

  S'_\beta\cap \sup(j\imageof\lambda)
  \mbox{ is stationary subset of }\sup(j\imageof\lambda)\\\mbox{ in }M}.
  \end{array}
  $
\end{xitemize}
By \Claimof{cl-bedrock-1},\,\assertof{4}, and since $j(\mu)<\lambda$, it follows that 
\begin{xitemize}
\item[] 
  $j\imageof\lambda=\setof{\beta<\sup(j\imageof\lambda)}{{}\begin{array}[t]{@{}l}

  S'_\beta\cap \sup(j\imageof\lambda)
  \mbox{ is stationary subset of }\sup(j\imageof\lambda)\\\mbox{ in }j(\uniW)}.
  \end{array}
  $ 
\end{xitemize}
{\ifprivate\darkred\fi It is apparent that} the right side of the equality is a definition 
of an element of $j(\uniW)$.
\qedofClaim
\qedskip

\begin{Claim}
  \Label{cl-bedrock-3} ${V_{j(\theta)}}^{\uniW}={V_{j(\theta)}}^{j(\uniW)}$. 
\end{Claim}
\noindent
\prfofClaim ``$\supseteq$'': By \Claimof{cl-bedrock-1-0}.

``$\subseteq$'': 
For $\lambda<j(\theta)$, if $X\in\psof{\lambda}\cap\uniW$, then
$j(X)\in j(\uniW)$,   
and $j\imageof{\lambda}\in j(\uniW)$ by \Claimof{cl-bedrock-2}. Thus
$j\imageof{X}=j(X)\cap j\imageof{\lambda}\in j(\uniW)$. Since we also have
$j\restr\lambda\in j(\uniW)$, it follows that $X\in j(\uniW)$.

Thus, by induction on $\alpha\leq j(\theta)$, we can prove that
${V_\alpha}^\uniW\subseteq j(\uniW)$. 
\qedofClaim\qedskip

Now, we have
\begin{xitemize}
\item[] 
 $j({V_\theta}^\uniW)\ubecause{\subseteq}{}{by
  ${V_\theta}^\uniW\subseteq {V_\theta}^{\overline{\uniW}}$ and by elementarity of $j$} 
  j({V_\theta}^{\overline{\uniW}})
  \obecause{\subseteq}{}{by \xitemof{x-bedrock-18} and \Claimof{cl-bedrock-1},\,\assertof{3}}
  {V_{j(\theta)}}^\uniW
  \ubecause{=}{}{\qquad\qquad by \Claimof{cl-bedrock-3}}{V_{j(\theta)}}^{j(\uniW)}
  =j({V_\theta}^\uniW)$.
\end{xitemize}

It follows $j({V_\theta}^\uniW)={V_{j(\theta)}}^\uniW$.
Since $j({V_\theta}^\uniV)={V_{j(\theta)}}^M$ is a generic extension of
${V_{j(\theta)}}=j({V_\theta}^\uniW)$ by a \po\ of size $\LE j(\kappa)$ 
(by \Claimof{cl-bedrock-1},\,\assertof{3}). Thus, by elementarity of $j$, 
${V_\theta}^\uniV$ is a generic extension of ${V_\theta}^\uniW$ by a \po\  of size
$\LE\kappa$. 

By definition of $\overline{\uniW}$ and by \xitemof{x-bedrock-9-1}, it follows that 
${V_\theta}^{\overline{\uniW}}\subseteq {V_\theta}^\uniW$. Since $\theta$ was arbitrary, 
this implies  
that $\overline{\uniW}\subseteq\uniW$ and hence $\overline{\uniW}=\uniW$ as 
desired. 
\memox{Scan\_2023-06-16--18.56-resurrection - annotated.pdf, pp.70--80\\
recurrence 5?.pdf, p.7-}
\memo{INSERT A REMARK saying that it is still open if the assumption of generic hyperhuge 
  is weakened referring [usuba2].}
\qedofThm\qedskip

As already noticed \Thmabove\ is a generalization of Theorem 1.6 in \cite{usuba} which 
states that there is the bedrock (of the universe $\uniV$) if there is a hyperhuge cardinal. 
In Theorem 1.3 in \cite{usuba2} the assumption of hyperhuge cardinal in this theorem is 
weakened to the existence of an extendible cardinal.

Though our proof of  \Thmof{p-bedrock-1} has a global structure quite similar to that of 
the proof of Theorem 1.3 in \cite{usuba2}, it seems to be necessary here that $\ol{\uniW}$ 
in our proof is the 
$\LE\kappa$-mantle while, in \cite{usuba2}, a similar proof could be done with $\ol{\uniW}$ which 
is the $\LT\kappa$-mantle, and that difference seems to force us to assume the stronger 
consistency strength of a hyperhuge cardinal.


\begin{Thm}
  \Label{p-bedrock-2} Suppose that $\calP$ is a class of \pos.
  If $\kappa$ is a tightly $\calP$-gen.\ hyperhuge cardinal, then $\kappa$ is a hyperhuge 
  cardinal in the bedrock $\overline{\uniW}$ of\/ $\uniV$. 
\end{Thm}
\prf In $\uniV$, let $\lambda>\kappa$ be arbitrary. Let $\theta>\lambda$ be sufficiently 
large \st\ it satisfies ${V_\theta}^\uniV\prec_{\Sigma_n}\uniV$ 
for sufficiently large natural number $n$ (this is just the condition 
\xitemof{x-bedrock-9-1}
in the proof of \Thmof{p-bedrock-1}). Let $\poQ\in\calP$ be \st\ for a
$(\uniV,\poQ)$-generic $\genH$, there are $j$, $M\subseteq\uniV[\genH]$ \st\
$\Elembed{j}{\uniV}{M}{\kappa}$, $\theta<j(\kappa)$, $j\imageof{j(\theta)}\in M$, and
${V_{j(\theta)}}^{\uniV[\genH]}={V_{j(\theta)}}^M$. Note that $\uniW$ in the proof of 
\Thmof{p-bedrock-1} coincides with the bedrock $\overline{\uniW}$ of $\uniV$. 

Thus we have\quad 
\ixitem[x-bedrock-19] $j\imageof{j(\lambda^*)}\in j(\overline{\uniW})$\quad for all
  $\lambda^*<j(\theta)$\quad by \Claimof{cl-bedrock-2}, and \ixitem[x-bedrock-20]
  $j({V_\theta}^{\overline{\uniW}})={V_{j(\theta)}}^{\overline{\uniW}}$\quad by 
\Claimof{cl-bedrock-3}.  

Let $\seqof{X_\alpha}{\alpha<\gamma}$ be an enumeration 
of $\psof{\psof{j(\lambda)}}^{\overline{\uniW}}$ in ${\overline{\uniW}}$. Note that $\gamma<j(\theta)$ by the choice 
of $\theta$.

By \xitemof{x-bedrock-20}, we have $\seqof{X_\alpha}{\alpha<\gamma}\in j({\overline{\uniW}})$ 
and $j\imageof{j(\lambda)}\in j({\overline{\uniW}})$ by \xitemof{x-bedrock-19}. 

Also for any $\lambda^*<j(\theta)$ and $A\in\calH((\lambda^*)^+)$ we have
$j\restr A\in j(\overline{\uniW})$ by \xitemof{x-bedrock-19} and Lemma 2.5 in 
\cite{sfetal-II}. It follows that
\begin{xitemize}
\xitem[x-bedrock-21] 
  $I:=\setof{\alpha<\lambda}{j\imageof{j(\lambda)}\in j(X_\alpha)}\in j(\overline{\uniW})$. 
\end{xitemize}

Let $U=\setof{X_\alpha}{\alpha\in I}$. Then $U\in j({\overline{\uniW}})$. 
Furthermore we have
\begin{xitemize}
\xitem[] $U\in {V_{j(\theta)}}^{j({\overline{\uniW}})}=j({V_\theta}^{\overline{\uniW}})
  \obecause{=}{}{by \xitemof{x-bedrock-20}}{V_{j(\theta)}}^{\overline{\uniW}}\subseteq {\overline{\uniW}}$.
\end{xitemize}

It is a routine to check that $U$ is a $\kappa$-complete normal ultrafilter 
over $\psof{j(\lambda)}^{\overline{\uniW}}$ and  
$\setof{x\in\psof{j(\lambda)}}{x\cap\kappa\in\kappa, \otp(x\cap j(\kappa))=\kappa,\,
  \otp(x)=\lambda}\in U$ holds. 

Since $\lambda$ was arbitrary, this implies that $\kappa$ is hyperhuge in ${\overline{\uniW}}$ by 
\Lemmaof{p-Lg-RcA-1-0}.\qedofThm

\begin{Cor}
  \Label{p-bedrock-3} Suppose that $\calP$ is an arbitrary class of \pos\ and 
  $\kappa$ is a tightly $\calP$-gen.\ hyperhuge cardinal.  Then \wassertof{1} there are 
  cofinally many huge cardinals. \smallskip

  \wassert{2} \SCH\ holds above $\kappa$. 
\end{Cor}
\prf Suppose that $\kappa$ is a tightly $\calP$-gen.\ hyperhuge cardinal. 
By \Thmof{p-bedrock-1} there is the bedrock $\overline{\uniW}$ and $\kappa$ is hyperhuge 
cardinal in $\overline{\uniW}$. \smallskip

\assertof{1}: 
Since the existence of a hyperhuge cardinal implies the 
existence of cofinally many huge cardinals (it is easy to show that the target $j(\kappa)$ 
of hyperhuge embedding for a sufficiently large inaccessible $\lambda$ is a huge cardinal), 
there are cofinally many huge cardinals in 
$\overline{\uniW}$. Since $\uniV$ is attained by a set forcing starting 
from $\overline{\uniW}$, a final segment of these huge cardinals survive in $\uniV$. 
\smallskip

\assertof{2}: By Theorem 20.8 in \cite{millennium-book}, \SCH\ holds above $\kappa$ in
$\overline{\uniW}$. Since $\uniV$ is a set generic extension of $\overline{\uniW}$ by a 
forcing of size $\LE\kappa$. \SCH\ 
should hold above $\kappa$. 
\qedofCor\qedskip

Compare \Corabove,\,\assertof{1} with \Corof{p-Lg-RcA-1-a-0} and \Corof{p-bedrock-17}.

For iterable stationary preserving $\calP$ containing all proper \pos, 
\Thmof{p-bedrock-3},\,\assertof{2} 
holds already under the $\calP$-Laver-gen.\ supercompactness of $\kappa$. The reason is 
that in such case \PFA\ holds (see Theorem 5.7 in \cite{sfetal-II}), and by Viale 
\cite{viale}, \SCH\ follows from it.  

In the following Corollary, we adopt the notation of Ikegami-Trang in \cite{ikegami-trang} on 
their version of Maximality Principle. 

\begin{Cor}
  \Label{p-bedrock-4} Suppose that $\calP$ is the class of semi-proper 
  \pos. 
  If $\kappa$ is a tightly $\calP$-Laver gen.\ hyperhuge cardinal, then 
  $\MP_{\Pi_2}(\omega_1,\mbox{all \pos})$ holds. 
\end{Cor}
\prf Theorem 1.8, (1) in \cite{ikegami-trang} states that the conclusion of the present 
corollary holds under $\MM^{++}$ and proper class many Woodin cardinals. 

If $\kappa$ is a $\calP$-Laver gen.\ supercompact cardinal, then $\kappa=\aleph_2$ and
$\MM^{++}$ holds (see \cite{sfetal-II}). By \Corof{p-bedrock-3}, tightly $\calP$-Laver 
gen.\ hyperhugeness implies that there are proper class many Woodin cardinals. 
\qedofCor

\begin{Cor}
  \Label{p-bedrock-5} Suppose that $\calP$ is the class of all \pos. 
  Then the following theories are 
  equiconsistent:\smallskip
  
  \wassert{a} \ZFC\ $+$ ``there is a hyperhuge cardinal''.

  \wassert{b} \ZFC\ $+$ ``there is a tightly $\calP$-Laver gen.\ hyperhuge cardinal''.

  \wassert{c} \ZFC\ $+$ ``there is a tightly $\calP$-gen.\ hyperhuge cardinal''.

  \wassert{d} \ZFC\ $+$ ``bedrock\/ $\overline{\uniW}$ exists and $\omega_1$ is a hyperhuge 
  cardinal in $\overline{\uniW}$''.\qed 
\end{Cor}

\begin{Cor}
  \Label{p-bedrock-6} Suppose that $\calP$ is one of the following classes of \pos: 
  all semi-proper \pos; all proper \pos; all ccc \pos; all $\sigma$-closed \pos. 
  Then the following theories are 
  equiconsistent:\smallskip
  
  \wassert{a} \ZFC\ $+$ ``there is a hyperhuge cardinal''.

  \wassert{b} \ZFC\ $+$ ``there is a tightly $\calP$-Laver gen.\ hyperhuge cardinal''.

  \wassert{c} \ZFC\ $+$ ``there is a tightly $\calP$-gen.\ hyperhuge cardinal''.

  \wassert{d} \ZFC\ $+$ ``bedrock\/ $\overline{\uniW}$ exists and $\kappa_\refl$ is a hyperhuge 
  cardinal in $\overline{\uniW}$''. \qed
\end{Cor}

The proof of \Thmof{p-bedrock-2} can be modified to obtain the following:

\begin{Thm}
  \Label{p-bedrock-2-0} Suppose that $\calP$ is a class of \pos.
  If a definable cardinal $\kappa$ is a tightly super-$C^{(\infty)}$-$\calP$-gen.\ 
  hyperhuge cardinal, then 
  $\kappa$ is super-$C^{(\infty)}$-hyperhuge 
  in the bedrock $\overline{\uniW}$ of $\uniV$. \qed
\end{Thm}

The definability of the cardinal $\kappa$ (e.g.\ as $\omega_1$, $\continuum$ etc.) in 
\Thmabove\ actually is needed so that the  
conclusion of the theorem is formalizable (in infinitely many formulas). 


\begin{Cor}
  \Label{p-bedrock-8} Suppose that $\calP$ is one of the following classes of \pos: 
  all semi-proper \pos; all proper \pos; all ccc \pos; all $\sigma$-closed \pos. 
  Then the following theories are 
  equiconsistent:\smallskip
  
  \wassert{a} \ZFC\ $+$ ``$c$ is a super-$C^{(\infty)}$ hyperhuge cardinal'' where $c$ is a 
  new constant symbol but ``... is super-$C^{(\infty)}$ hyperhuge ...'' is formulated as an 
  infinite collection of formulas in $\Lin$. \smallskip

  \wassert{b} \ZFC\ $+$ ``there is a tightly super-$C^{(\infty)}$-$\calP$-Laver gen.\ 
  hyperhuge cardinal''.\smallskip

  \wassert{c} \ZFC\ $+$ ``bedrock\/ $\overline{\uniW}$ exists and ${\kappa_\refl}^\uniV$ is a super
  $C^{(\infty)}$-hyperhuge cardinal in $\overline{\uniW}$''. \qed
\end{Cor}

\begin{Thm}
  \Label{p-bedrock-8-0} Suppose that $\kappa$ is tightly 
  super-$C^{(\infty)}$-$\calP$-Laver gen.\ hyperhuge for an iterable $\calP$. Then, for 
  each $n\in\natnums$, there are stationarily many super-$C^{(n)}$-hyperhuge cardinals. 
\end{Thm}
\prf Let $n\in\natnums$. 
By a modification of \Thmof{p-bedrock-2}, we can show that $\kappa$ is super-
$C^{(\infty)}$-hyperhuge in the bedrock $\overline{\uniW}$. 
Thus there are stationarily many super-$C^{(n)}$-hyperhuge cardinals in $\overline{\uniW}$ 
by \Corof{p-Lg-RcA-1-a-0}. Since $\uniV$ is a set generic extension of $\overline{\uniW}$, 
the stationarity of these classes is preserved by the generic extension.
\qedofThm

\begin{Cor}
  \Label{p-bedrock-8-1} The consistency strength of the existence of a tightly 
  super-$C^{(\infty)}$-$\calP$-Laver gen.\ hyperhuge cardinal for one of the 
  iterable $\calP$'s in \Thmof{p-Lg-RcA-4} is strictly between that of the existence of 
  a super $C^{(n)}$-hyperhuge cardinal, and that of the existence of a 2-huge cardinal.
\end{Cor}
\prf Suppose that $\kappa$ is tightly 
  super-$C^{(\infty)}$-$\calP$-Laver gen.\ hyperhuge for an iterable $\calP$. Then by 
  \Thmof{p-bedrock-8-0}, there is a super-$C^{(n)}$-hyperhuge $\lambda>\kappa$ with
  $V_\lambda\prec_{\Sigma_{n'}}$ for sufficiently large $n'>n$. any 
  super-$C^{(n)}$-hyperhuge $\lambda_0<\lambda$ is super-$C^{(n)}$-hyperhuge 
  in $V_\lambda$ by elementarity. Thus $V_\lambda$ is a model of \ZFC\ $+$ ``there is a super-
  $C^{(n)}$-hyperhuge cardinal''.

  If $\kappa$ is 2-huge then by \Lemmaof{p-Lg-RcA-2} and \Thmof{p-Lg-RcA-4}, there is a set 
  model with a tightly 
  super-$C^{(\infty)}$-$\calP$-Laver gen.\ hyperhuge cardinal. \qedofCor

\section{Bedrock and Laver genericity} 
\Label{bedrock-Lg}
\begin{Prop}\Labelx{p-bedrock-9}{\qquad\qquad\qquad} Suppose that $\calP$ is a class of 
  \pos\ and 
  $\kappa$ is  
    tightly $\calP$-gen.\ ultrahuge cardinal. By \Thmof{p-bedrock-1} there is the bedrock
    $\overline{\uniW}$ of\/ $\uniV$.

    \quad\ 
    We have $(\kappa^+)^{\overline{\uniW}}=(\kappa^+)^\uniV$ and 
    $\uniV$ is a set generic extension $\overline{\uniW}$ by some  
    \po\ $\poP\in\overline{\uniW}$ \st\ $\overline{\uniW}\models\cardof{\poP}\leq\kappa$. 
\end{Prop}
\memo{Scan\_2023-06-16--18.56-resurrection - annotated.pdf, pp.68--69\\
  recurrence 5?.pdf, p.4-\\
  Scan\_2023-06-16--18.56-resurrection - annotated.pdf, pp.127--131}
\prf Let $\poP_0\in{\overline{\uniW}}$ be \st\ there is 
a $({\overline{\uniW}},\poP_0)$-generic $\genG_0$ \st\ 
$\uniV=\overline{\uniW}[\genG_0]$. By the proof of \Thmof{p-bedrock-1}, $\poP_0$ can be 
chosen \st\ \ixitemx[x-bedrock-21-0]{\quad\quad} $\uniV\models\cardof{\poP_0}\leq\kappa$.

Let \ixitem[x-bedrock-22] $\theta>\kappa+\cardof{\poP_0}$ 
(in $\uniV$) be large enough, and \st\ it satisfies \xitemof{x-bedrock-9-1}
in the proof of \Thmof{p-bedrock-1} for sufficiently large natural number $n$. 
Let $\poQ\in\calP$ be \st\ for $(\uniV,\poQ)$-generic $\genH$ there are $j$,
$N\subseteq\uniV[\genH]$ \st\ $\Elembed{j}{\uniV}{M}{\kappa}$, $j(\kappa)>\theta$,
\ixitem[x-bedrock-22-0] $\cardof{\poQ}\leq j(\kappa)$,\quad $\genH$,
$j\imageof{j(\theta)}\in M$, and 
${V_{j(\theta)}}^{\uniV[\genH]}\subseteq M$. 

Since $\cardof{\poP_0}<\theta<j(\kappa)$ we have
\begin{xitemize}
  \xitem[x-bedrock-22-1] 
  $\overline{\uniW}\models\cardof{\poP_0}<j(\kappa)$.
\end{xitemize}

Thus\vspace{-2ex}
\begin{xitemize}
  \xitem[x-bedrock-23] 
  $(j(\kappa)^+)^{j(\overline{\uniW})}
    =(j(\kappa)^+)^{{\uniV_{j(\theta)}}^{j(\overline{\uniW})}}\!\!
  \ubecause{=}{}{by \Claimof{cl-bedrock-3}}(j(\kappa)^+)^{{\uniV_{j(\theta)}}^{\overline{\uniW}}}
  \obecause{=}{}{by \xitemof{x-bedrock-22-1}}
  (j(\kappa)^+)^{{\uniV_{j(\theta)}}^{\uniV}}$\\[-2ex]
  $\obecause{=}{}{by \xitemof{x-bedrock-22-0}}
  (j(\kappa)^+)^{{\uniV_{j(\theta)}}^{\uniV[\genH]}}. $
\end{xitemize}
By elementarity it follows that $(\kappa^+)^\uniW=(\kappa^+)^\uniV$. By 
\xitemof{x-bedrock-21-0} this implies that\\
$\overline{\uniW}\models\cardof{\poP_0}\leq\kappa$. 
\qedofProp

\memo{Scan\_2023-06-16--18.56-resurrection - annotated.pdf, pp.129〜}
\begin{Prop}
  \Label{p-bedrock-10} Suppose that $\kappa$ is tightly $\calP$-gen.\ hyperhuge for a 
  class $\calP$ of \pos. By \Propof{p-bedrock-9}, there is a \po\ $\poP$ in the 
  bedrock $\overline{\uniW}$ with $\overline{\uniW}\models\cardof{\poP}\leq\kappa$ \st\
  $\uniV=\overline{\uniW}[\genG]$ for a $(\overline{\uniW},\poP)$-generic $\genG$.

  \quad\ 
  For any bounded $b\subseteq\kappa$ in $\uniV$, there is $\poP_b\in\overline{\uniW}$ with
  $\overline{\uniW}\models\cardof{\poP_b}<\kappa$, $\poP_b\subseteq\poP$, \st\
  $\genG\cap\poP_b$ is $(\overline{\uniW},\poP_b)$-generic and 
  $\overline{\uniW}[\genG\cap\poP_b]\ni b$. 
\end{Prop}
\prf Since $\overline{\uniW}\modelof{\cardof{\poP}\leq\kappa}$ we may assume, \wolog, that 
the underlying set of $\poP$ is $\kappa$. Let $\theta>\cardof{\poP}$ be large enough and 
\st\ $V_\theta\prec_{\Sigma_n}\uniV$ for a large enough $n$ in accordance with 
\Lemmaof{p-Lg-RcA-0-0}.

Let $\poQ\in\calP$ be \st\ for $(\uniV,\poQ)$-generic $\genH$, there are $j$,
$M\subseteq\uniV[\genH]$ \st\ $\Elembed{j}{\uniV}{M}{\kappa}$, $\theta< j(\kappa)$,
$\cardof{\poQ}\leq j(\kappa)$, ${V_{j(\theta)}}^{\uniV[\genH]}\subseteq M$ and
$j\imageof{j(\theta)}\in M$ (see \Lemmaof{p-bedrock-0},\,\assertof{5}).

By \Claimof{cl-bedrock-3}, we have \ixitem[x-bedrock-24] $j({V_\theta}^{\overline{\uniW}})$
($= {V_{j(\theta)}}^{j(\overline{\uniW})}$) $={V_{j(\theta)}}^{\overline{\uniW}}$. By the 
choice  
of $\theta$, we have \ixitem[x-bedrock-25]
${V_{j(\theta)}}^\uniV={V_{j(\theta)}}^{\overline{\uniW}}[\genG]$.  

Since the underlying set of $\poP$ is $\kappa$ we have $j(\poP)\cap\kappa=\poP$ and
$j(\genG)\cap\kappa=\genG$. By $b=j(b)\in {V_{j(\theta)}}^\uniV
\obecause{=}{}{by \xitemof{x-bedrock-24} and \xitemof{x-bedrock-25}}
j({V_\theta}^{\overline{\uniW}})[\genG]=j({V_\theta}^{\overline{\uniW}})[j(\genG\cap\kappa)]$,
and since $\cardof{\poP}<\theta<j(\kappa)$, we have
\begin{xitemize}
\item[] $M\modelof{\,
  \begin{array}[t]{@{}l}
    \mbox{there is }P\in j(\overline{\uniW})
    \mbox{ with }P\subseteq j(\poP)\mbox{ and }
    j(\overline{\uniW})\modelof{\cardof{P}<j(\kappa)}\\
    \mbox{\st\ }j(\overline{\uniW})[j(\genG)\cap P]\ni j(b)}.
  \end{array}$
\end{xitemize}
By elementarity, it follows that
\begin{xitemize}
\item[] $\uniV\modelof{\,
  \begin{array}[t]{@{}l}
    \mbox{there is }P\in\overline{\uniW}\mbox{ and }G\subseteq P
    \mbox{ with }P\subseteq \poP\mbox{ and }
    \overline{\uniW}\modelof{\cardof{P}<\kappa}\\
    \mbox{\st\ }G\mbox{ is }(\overline{\uniW},P)\mbox{-generic filter and }
    \overline{\uniW}[\genG\cap P]\ni j(b)}
  \end{array}$
\end{xitemize}
as desired. \qedofProp
\qedskip

\begin{Cor}
  \Label{p-bedrock-11} Suppose $\kappa$ is a 
  tightly $\calP$-gen.\ hyperhuge cardinal for a class $\calP$ of \pos. Then 
  we have $2^{\LT\kappa}=\kappa$.\ifextended\else\qed\fi
\end{Cor}
{\ifextended\extendedcolor
\prf $\kappa$ is hyperhuge in the bedrock $\overline{\uniW}$ of $\uniV$ by 
\Thmof{p-bedrock-2}.  
Thus, $\overline{\uniW}\models2^{\LT\kappa}=\kappa$. It follows that essentially there are 
at most $\kappa$ many (nice) $\poP$-names of elements of $[\kappa]^{\LT\kappa}$ for \pos\
$\poP$  
of size $\LT\kappa$. Since all elements of $([\kappa]^{\LT\kappa})^\uniV$ are realizations 
of names of this kind by \Propof{p-bedrock-10}, it follows that $\uniV\models2^{\LT\kappa}=\kappa$. 
\qedofCor\qedskip\fi}

Hamkins \cite{hamkins} proved that if $\calP$ contains $\Col(\omega,\lambda)$ for every
$\lambda$,  and $(\calP,\calH(\aleph_1))$-\RcA\ holds, then we have $L_\kappa\prec L$ for
$\kappa=\omega_1$. Practically the same proof concludes that the mantle $\uniW$ (the 
intersection of all grounds which is shown to be a model of \ZFC\ in \cite{usuba}) also 
satisfies ${V_{\kappa}}^\uniW\prec\uniW$ for $\kappa={\omega_1}^\uniV$.
Note that, since grounds are downward directed, the mantle of the universe $\uniV$ is also 
the mantle of any ground $\uniW$ of $\uniV$.

Actually, we can say a little bit more.
\begin{Lemma}\Label{p-bedrock-12} Let\/ $\uniW^*$ be the mantle of\/ $\uniV$. 
  \wassertof{1} Suppose that 
  $(\calP,\calH(\continuum)^\uniV)$-\RcA\  
  holds for  
  a class $\calP$ of \pos\ \st\ either $\calP$ 
  contains \pos\ collapsing $\lambda$ to be countable for cofinally many 
  cardinals $\lambda$, or it contains \pos\ forcing $\continuum$ arbitrary large adding 
  reals without collapsing cardinals below the number of the reals added. 

  \quad\ Then ${V_\alpha}^{\uniW^*}$ for all $\alpha<\continuum$ is of cardinality
  $\LT\continuum$ in $\uniV$.\smallskip

  \wassert{2} Suppose that
  $(\calP,\calH(\omega_1)^\uniV)$-\RcA\  
  holds for  
  a class $\calP$ of \pos\ \st\ $\calP$ 
  contains \pos\ collapsing $\lambda$ to be countable for cofinally many 
  cardinals $\lambda$.  

  \quad\ Then ${V_\alpha}^{\uniW^*}$ for all $\alpha<(\omega_1)^\uniV$ is countable 
  in $\uniV$.\smallskip 
  
\end{Lemma}
\prf \assertof{1}: For $\alpha<\continuum$, 
``\,${V_\alpha}^{W^*}$ is of cardinality $\LT\continuum$\,'' can be formulated as 
an $\Lin$-formula with the the parameter $\alpha\in\calH(\continuum)$, and it is forcable 
by a \po\ in 
$\calP$. By $(\calP,\calH(\continuum)^\uniV)$-\RcA, there is a ground $\uniW$ of $\uniV$
($\uniW^*\subseteq\uniW\subseteq\uniV$)  
\st\ the statement above holds in $\uniW$. Then
$\uniW\modelof{{V_\alpha}^{W^*}\mbox{ is of cardinality }\LT\continuum}$, and hence 
$\uniV\modelof{{V_\alpha}^{W^*}\mbox{ is of cardinality }\LT\continuum}$. \smallskip

\assertof{2}: A proof similar to that of \assertof{1} will do.

For $\alpha<\omega_1$, ``${V_\alpha}^{W^*}$ is of cardinality $\LT\aleph_1$'' is a 
statement represented as 
an $\Lin$-formula with the parameter $\alpha\in\calH(\aleph_1)$, and it is forcable by a 
\po\ in 
$\calP$. By $(\calP,\calH(\omega_1)^\uniV)$-\RcA, there is a ground $\uniW$ of $\uniV$
($\uniW^*\subseteq\uniW\subseteq\uniV$)  
\st\ the statement above holds in $\uniW$. Then 
$\uniW\modelof{{V_\alpha}^{W^*}\mbox{ is countable}}$. 
and hence 
$\uniV\modelof{{V_\alpha}^{W^*}\mbox{ is countable}}$. 
\qedofLemma\qedskip

Compare the following proposition with \Lemmaof{p-Lg-RcA-1-a}: 
\begin{Prop}
  \Label{p-bedrock-13} Let $\uniW^*$ be the mantle of\/ $\uniV$. 
  \wassert{1} Suppose that $(\calP,\calH(\continuum)^\uniV)$-\RcA\  
  holds for $\calP$ as in \Lemmaof{p-bedrock-12},\assertof{1}. That is, 
  either $\calP$ 
  contains \pos\ collapsing $\lambda$ to be countable for cofinally many 
  cardinals $\lambda$, or it contains \pos\ forcing $\continuum$ arbitrary large adding 
  reals without collapsing small cardinals. 

  \quad\ 
  Then we have ${V_{\kappa}}^{\uniW^*}\prec\uniW^*$ for $\kappa=(\continuum)^\uniV$.
  \smallskip

  \wassert{2} Suppose that
  $(\calP,\calH((\aleph_1)^\uniV))$-\RcA\  
  holds for  
  a class $\calP$ of \pos\ \st\ $\calP$ 
  contains \pos\ collapsing $\lambda$ to make it countable for cofinally many 
  cardinals $\lambda$ in $\On$.

  \quad\
  Then we have ${V_{\kappa}}^{\uniW^*}\prec\uniW^*$ for $\kappa={\omega_1}^\uniV$.
  
\end{Prop}
\prf \assertof{1}: We show that $(V_{(\continuum)^\uniV})^{\uniW^*}$ in $\uniW^*$ passes the 
Vaught's test. Suppose 
$\uniW^*\models\varphi(\overline{a},b)$  
where $\overline{a}\in(V_{(\continuum)^\uniV})^{\uniW^*}$ and $b\in\uniW^*$. We have
$\overline{a}\in\calH(\continuum)^\uniV$ by \Lemmaof{p-bedrock-12},\,\assertof{1}. The  
statement 
\begin{xitemize}
\item[] 
 $\psi:=\exists y(y\in\uniW^*\land\cardof{\trcl(y)}<\continuum
\land\uniW^*\models\varphi(\overline{a},y))$ 
\end{xitemize}
is an $\Lin$-formula the parameters $\overline{a}\in\calH(\continuum)^\uniV$ and 
forcable by a \po\ in $\calP$ (just by collapsing a large enough cardinal to countable). By
$(\calP,\calH(\continuum))$-\RcA, it follows that there  
is a ground $\uniW$ of $\uniV$ (so $\uniW^*\subseteq\uniW\subseteq\uniV$) \st\ 
\ixitem[x-bedrock-26] $\uniW\models\psi$.

Let $b'\in\uniW$ be a witness for \xitemof{x-bedrock-26}. Then we have
\ixitem[x-bedrock-27] $\uniW\models\cardof{\trcl(b')}<\continuum$, and
$\uniW^*\models\varphi(\overline{a},b')$.

Since $(\continuum)^\uniW\leq(\continuum)^\uniV$, 
it follows that
$\uniW^*\models\cardof{\trcl(b')}<(\continuum)^\uniV$ 
by \xitemof{x-bedrock-27}, and thus
$b'\in (V_{(\continuum)^\uniV})^{\uniW^*}$. \smallskip

\assertof{2}: can be proved similarly to \assertof{1}. {\ifextended\extendedcolor
\par
We show that $(V_{{\omega_1}^\uniV})^{\uniW^*}$ in $\uniW^*$ passes the 
Vaught's test. Suppose 
$\uniW^*\models\varphi(\overline{a},b)$  
where $\overline{a}\in(V_{{\omega_1}^\uniV})^{\uniW^*}$ and $b\in\uniW^*$. We have
$\overline{a}\in\calH(\aleph_1)^\uniV$ by \Lemmaof{p-bedrock-12},\,\assertof{2}. The  
statement 
\begin{xitemize}
\item[] 
 $\psi:=\exists y(y\in\uniW^*\land\cardof{\trcl(y)}<\aleph_1
\land\uniW^*\models\varphi(\overline{a},y))$ 
\end{xitemize}
is $\Lin$-formula with the parameters $\overline{a}\in\calH(\aleph_1)^\uniV$ and 
forcable by a \po\ in $\calP$. By $(\calP,\calH(\aleph_1)^\uniV)$-\RcA, it follows that 
there  
is a ground $\uniW$ of $\uniV$ (so $\uniW^*\subseteq\uniW\subseteq\uniV$) \st\ 
\ixitema[x-bedrock-26-0] $\uniW\models\psi$.

Let $b'\in\uniW$ be a witness for \xitemof{x-bedrock-26-0}. Then we have
\ixitema[x-bedrock-27-0] $\uniV\models\cardof{\trcl(b')}<\aleph_1$, and
$\uniW^*\models\varphi(\overline{a},b')$.

Since ${\omega_1}^\uniW\leq{\omega_1}^\uniV$, it follows that
$\uniW^*\models\cardof{\trcl(b')}<(\aleph_1)^\uniV$ by \xitemof{x-bedrock-27-0}, and thus
$b'\in (V_{{\omega_1}^\uniV})^{\uniW^*}$. 
\fi}\qedofProp\qedskip

\begin{Thm}\Label{p-bedrock-16} Suppose that $\calP$ is a class  of 
  \pos\ and $\kappa:=(\omega_1)^\uniV$ is tightly $\calP$-gen.\ hyperhuge. 
  Then 
  \tfae:\smallskip 

  \wassert{a} $(\mbox{all \pos},\calH(\kappa))$-\RcA\ holds.\smallskip

  \wassert{b} $(\mbox{all \pos},\calH(\kappa))$-\RcAp\ holds.\smallskip

  \wassert{c} ${V_\kappa}^{\overline{\uniW}}\prec\overline{\uniW}$ where
  $\overline{\uniW}$  
  is the bedrock of\/ $\uniV$.
\end{Thm}
\prf \assertof{a} $\Leftrightarrow$ \assertof{b}: is trivial.
\assertof{a} $\Rightarrow$ \assertof{c}: By \Propof{p-bedrock-13},\,\assertof{2}.\smallskip

\assertof{c} $\Rightarrow$ \assertof{b}: Assume
\ixitem[x-bedrock-29] ${V_\kappa}^{\overline{\uniW}}\prec\ol{\uniW}$, and
$\forces{\poP}{\varphi(\ol{a}\checked)}$, for a \po\ $\poP$, an $\Lin$-formula
$\varphi=\varphi(\ol{x})$, and 
$\ol{a}\in\calH(\kappa)^{\uniV}$.

By \Propof{p-bedrock-10}, 
there are \po\
$\poQ\in \ol{\uniW}$ with $\ol{\uniW}\models\cardof{\poQ}<\kappa$, and
$(\ol{\uniW}, \poQ)$-generic 
$\genG\in\uniV$ \st\ $\ol{a}\in\ol{W}[\genG]$. Let $\utilde{\ol{a}}$ be a $\poQ$-name of
$\ol{a}$.

By the choice of $\varphi$ and $\ol{a}$, we have
\begin{xitemize}
\item[] 
  $\ol{\uniW}\models\forces{\poQ}{\mbox{there is a \po\ }P
  \mbox{ which forces }\varphi(\utilde{\ol{a}}\checked)}$. 
\end{xitemize}

By the elementarity \xitemof{x-bedrock-29}, it follows that
${V_\kappa}^{\ol{\uniW}}\models\forces{\poQ}{\xmbox{there is a \po\ }P 
  \xmbox{ which forces }\varphi(\utilde{\ol{a}}\checked)}$. 
Thus, there is a $\poQ$-name
$\utpoR\in{V_\kappa}^{\ol{\uniW}}$ \st\
$\forces{\poQ\ast\utpoR}{\varphi((\utilde{\ol{a}}\checked)}$.

We have 
$\poQ\ast\utpoR, \psof{\poQ\ast\utpoR}^{\ol{\uniW}}
\in {V_\kappa}^{\ol{\uniW}}
=\calH(\kappa)^{\ol{\uniW}}\subseteq\calH(\kappa)^\uniV$, and $\kappa$ is $\omega_1$ in
$\uniV$. Hence,  we can construct a
$(\ol{\uniW}[\genG],\utpoR[\genG])$-generic $\genG'\in\uniV$ in $\omega$ steps in $\uniV$. 
$\ol{\uniW}[\genG][\genG']$ is a ground in $\uniV$ and
$\ol{\uniW}[\genG][\genG']\models\varphi(\ol{a})$. \qedofThm\qedskip  


The following corollary is proved similarly to \Corof{p-Lg-RcA-1-a-0}. 

\begin{Cor}\Label{p-bedrock-17} Suppose that $\calP$ is a class  of 
  \pos\ and $\omega_1$ is tightly $\calP$-gen.\ hyperhuge. If
  $(\mbox{all \pos}, \calH(\aleph_1))$-\RcA\ holds then there are stationarily many hyperhuge 
  cardinals. More precisely, under this condition, for any club subclass $\calC$ of $\On$ 
  defined with a parameter, there is a hyperhuge cardinal in $\calC$. In particular, there 
  are class may hyperhuge cardinals.\ifextended\else\qed\fi
\end{Cor}
{\ifextended\extendedcolor
\prf Since $\uniV$ is a set generic extension of the bedrock $\ol{\uniW}$, it is enough to 
show  
that there are stationarily many hyperhuge cardinals in $\ol{\uniW}$ (in the same sense as in 
the statement of the corollary).

Suppose this is not the case. Then there is an $\Lin$-formula $\Phi=\Phi(x,y)$ \st\ 
\begin{xitemize}
\item[] $\ol{\uniW}\modelof{\exists y\,(\,{}
  \begin{array}[t]{@{}l}
    \Phi(\cdot,y)\mbox{ is a club in }\On\\
    \mbox{but }\Phi(\cdot,y)\mbox{ does not contain any hyperhuge cardinal})}.
  \end{array}$
\end{xitemize}

Since we have ${V_{(\omega_1)^\uniV}}^{\overline{\uniW}}\prec\overline{\uniW}$ by 
\Thmof{p-bedrock-16},  
it follows that 

\begin{xitemize}
\xitemA[x-bedrock-30] ${V_{(\omega_1)^\uniV}}^{\overline{\uniW}}\modelof{\exists y\,(\,{}
  \begin{array}[t]{@{}l}
    \Phi(\cdot,y)\mbox{ is a club in }\On\\
    \mbox{but }\Phi(\cdot,y)\mbox{ does not contain any hyperhuge cardinal})}.
  \end{array}$
\end{xitemize}
Let $b\in{V_{(\omega_1)^\uniV}}^{\overline{\uniW}}$ be a witness of \xitemAof{x-bedrock-30}. 
Then we have $(\omega_1)^\uniV\in\Phi(\cdot,b)$ by the closedness of $\Phi(\cdot,b)$. But 
this is a contradiction since $(\omega_1)^\uniV$ is a hyperhuge 
cardinal by \Thmof{x-bedrock-2}. \qedofCor\qedskip\fi}

The following proposition is a variation of Theorem 5.7 in \cite{sfetal-II}. 
\begin{Prop}\Label{p-bedrock-14} Suppose that $\kappa$ is tightly $\calP$-Laver gen.\  
  supercompact for an iterable class $\calP$ of \pos. Then we have $\MA(\calP,\LT\kappa)$. 
\end{Prop}
\prf Suppose that $\poP\in\calP$ and $\calD$ is a family of dense subsets of $\poP$ with
$\cardof{\calD}<\kappa$. Let $\utpoQ$ be a
$\poP$-name of a \po\ \st, for 
$(\uniV,\poP\ast\utpoQ)$-generic $\genH$, there are $j$, $M\subseteq\uniV[\genH]$ \st\
$\Elembed{j}{\uniV}{M}{\kappa}$, $\calP$, $\genH\in M$,
$j\imageof{\lambda}\in M$, and $\cardof{\poP\ast\poQ}\leq j(\kappa)$. 

Note that $j(\calD)=\setof{j(D)}{D\in\calD}$. Let $\genG$ be the $\poP$ part of $\genH$. We 
have $\genG\in M$. Thus
\begin{xitemize}
\item[] $\genG^*=\setof{\condp\in j(\poP)}{j(\condq)\leq_{j(\poP)}\condp
  \mbox{ for some }\condq\in\genG}$
\end{xitemize}
is an element in $M$. Since $\genG^*$ is $j(\calD)$-generic filter over $j(\poP)$, 
\begin{xitemize}
\item[] $M\models\exists G\,(G\mbox{ is a }j(\calD)\mbox{-generic filter over }j(\poP))$. 
\end{xitemize}
By elementarity, it follows that 
\begin{xitemize}
\item[] $\uniV\models\exists G\,(G\mbox{ is a }\calD\mbox{-generic filter over }\poP)$. 
\qedofProp
\end{xitemize}
\qedskip

\begin{Lemma}\Label{p-bedrock-15}
  Suppose that $\kappa$ is a tightly $\calP$-gen.\ hyperhuge cardinal for a 
  class $\calP$ of \pos, and $a\in\calH(\kappa)$ is \st\
  $\uniV\models\psi(a)$ for some $\Lin$-formula $\psi=\psi(x)$. 
  Let\/ $\overline{\uniW}$ be  
  the bedrock of\/ $\uniV$.

  \quad\ Then there is
  $\poP^*\in {V_\kappa}^{\overline{\uniW}}$  
  with $\overline{\uniW}\models\cardof{\poP^*}<\kappa$,  
  and $(\overline{\uniW},\poP^*)$-generic $\genG^*\in\uniV$ \st\
  $a\in\overline{\uniW}[\genG^*]$,  
  $\overline{\uniW}[\genG^*]\modelof{
    \psi(a)}$, and\/
  $\overline{\uniW}[\genG^*]$ is a $\calP$-ground of\/ $\uniV$. 
\end{Lemma}
\prf Assume that $\uniV=\overline{\uniW}[\genG]$ where $\genG$ is 
a $(\overline{\uniW},\poP)$-generic filter over a $\poP\in\overline{\uniW}$ with
$\overline{\uniW}\modelof{\cardof{\poP}\leq\kappa}$ (\Propof{p-bedrock-9}). 
\Wolog, we shall assume that the underlying set of $\poP$ is $\kappa$.

Let $\lambda>\kappa$ 
be \st\ \ixitem[x-bedrock-27-1] $V_\lambda\prec_{\Sigma_n}\uniV$ for a sufficiently 
large $n$.\footnote{Here ``sufficiently large $n$'' refers, among other things, largeness 
  in terms of \Lemmaof{p-Lg-RcA-0-0}, and the absoluteness of ``$\poP\in\calP$''. } 

Let $\poQ\in\calP$ be \st, for $(\uniV,\poQ)$-generic $\genH$, there are $j$,
$M\subseteq\uniV[\genH]$ \st\ $\Elembed{j}{\uniV}{M}{\kappa}$, $j(\kappa)>\lambda$,\quad
\ixitem[x-bedrock-28] ${V_{j(\lambda)}}^\uniV[\genH]\in M$, and $\cardof{\poQ}=j(\kappa)$ 
(see  
\Lemmaof{p-bedrock-0},\,\assertof{5}).

Note that we have ${V_{j(\lambda)}}^M\prec_{\Sigma_n}M$ by elementarity. By 
\Lemmaof{p-Lg-RcA-0-0}, it follows that
$M\models{V_{j(\lambda)}={V_{j(\lambda)}}^\uniV[\genH]}$.
Also, by \Lemmaof{p-Lg-RcA-0-0} and \xitemof{x-bedrock-28}, we have
${V_{j(\lambda)}}^{\uniV[\genH]}={V_{j(\lambda)}}^\uniV[\genH]={V_{j(\lambda)}}^M$.

Thus, for $\delta<\kappa$, noting $j(\delta)=\delta$, we have 
\begin{xitemize}
\item[] 
  $M\modelof{{}
  \begin{array}[t]{@{}l}
    V_{j(\lambda)}\mbox{ is a }\calP\mbox{-generic extension of a ground }W\mbox{ which is 
      a model of }\\
    \psi(j(a)),
    \mbox{and }W\mbox{ is a }P\mbox{-generic extension of }
    {V_{j(\lambda)}}^{\overline{\uniW}}\mbox{ for some \po\ }\\
    P\mbox{ of size }<j(\kappa).}
  \end{array}$
\end{xitemize}

By elementarity, it follows that 
\begin{xitemize}
\item[] 
  $\uniV\modelof{{}
  \begin{array}[t]{@{}l}
    V_{\lambda}\mbox{ is a }\calP\mbox{-generic extension of a ground }W\mbox{ which is 
      a model of }\\
    \psi(a),
    \mbox{and }W\mbox{ is a }P\mbox{-generic extension of }
    {V_{\lambda}}^{\overline{\uniW}}\mbox{ for some \po\ }\\
    P\mbox{ of size }<\kappa.}
  \end{array}$
\end{xitemize}

Now by the choice \xitemof{x-bedrock-27-1} of $\lambda$, it follows that
\begin{xitemize}
\item[] 
  $\uniV$ is a $\calP$-generic extension of a ground $\uniW$ which is 
  a model of 
  $\psi(a)$, and $\uniW$ is 
  a $\poP^*$-generic extension of 
  $\overline{\uniW}$ for  
  a \po\ $\poP^*$ of size $<\kappa$. \qedofLemma
\end{xitemize}

{\ifextended\extendedcolor
A \po\ $\poP$ has pre-caliber $\kappa$ if for any $A\in[\poP]^{\GE\kappa}$ there is
$B\in[A]^{\GE\kappa}$ \st\ $B$ is centered (i.e.\ $b$ has a lower bound in $\poP$ for and
$b\in[B]^{\LT\aleph_0}$). 
\begin{LemmaA}{\rm (Lemma III, 3.35, in \cite{kunen-2011})}\Label{p-bedrock-18}
  $\MA(\aleph_1)$ implies that all ccc \po\ $\poP$ has pre-caliber $\aleph_1$. 
\end{LemmaA}
\prf
\memo{Scan\_2023-06-16--18.56-resurrection.pdf p.137}
Suppose $\condp_\alpha\in\poP$ for $\alpha<\omega_1$. It is enough to find a filter 
containing uncountably many of $\condp_\alpha$'s.

For each $\alpha<\omega_1$, let
$D_\alpha=\setof{\condq\in \poP}{\condq\leq_\poP\condp_\beta
\mbox{ for some }\beta\geq\alpha}$. 

If there is no $\conds\in\poP$ \st\ uncountably many $D_\alpha$'s are dense below $\pos$, 
we can construct an strictly increasing sequence $\seqof{\beta_\alpha}{\alpha<\omega_1}$ in
$\omega_1$ \st\ $\condp_{\beta_\alpha}$ is incompatible with all $\condp_{\beta_\gamma}$,
$\gamma<\alpha$. But then t $\setof{\condp_{\beta_\alpha}}{\alpha<\omega_1}$ is a pairwise 
incompatible uncountable set which is a contradiction to the ccc of $\poP$.

Thus, there is $\conds\in\poP$ \st\ there are uncountably many $\alpha$'s \st\ 
$\calD_\alpha$ is dense below $\conds$. Then since $D_\alpha$, $\alpha<\omega_1$ build a 
decreasing sequence,  $D_\alpha$ is dense below $\conds$ for all $\alpha<\omega_1$.
Let $\calD:=\setof{D_\alpha}{\alpha<\omega_1}$ and let $\genG$ be a $\calD$-generic filter 
over $\poP$. Then $\condp_\alpha\in\genG$ for all $\alpha\in I$ as desired. 
\qedofLemmaA

\begin{CorA}\Label{p-bedrock-19}
  $\MA(\aleph_1)$ implies that, for all ccc \pos\ $\poP$, $\poQ$, the product
  $\poP\times\poQ$ has the ccc. \qed 
\end{CorA}
\prf It is enough to show that $\poP\times\poQ$ has pre-caliber $\aleph_1$ by 
\LemmaAof{p-bedrock-18}. Suppose that 
$\pairof{\condp_\alpha,\condq_\alpha}\in\poP\times\poQ$ for $\alpha<\omega_1$. Since $\poP$ 
has pre-caliber $\aleph_1$, there is $I_0\in[\omega_1]^{\aleph_1}$ \st\
$\setof{\condp_\alpha}{\alpha\in I_0}$ is centered in $\poP$. Since $\poQ$ also has 
pre-caliber $\aleph_1$ by \LemmaAof{p-bedrock-18}, there is $I_1\subseteq[I_0]^{\aleph_1}$ 
\st\ 
$\setof{\condq_\alpha}{\alpha\in I_1}$ is centered in $\poQ$. 
$\setof{\pairof{\condp_\alpha,\condq_\alpha}}{\alpha\in I_1}$ is centered inc
$\poP\times\poQ$. \qedofCorA\qedskip\fi}

The following Lemma is classical:

\begin{Lemma}\Label{p-bedrock-20} {\ifextended\extendedcolor\wassertof{1}
    Suppose that $\uniV\modelof{\poP\times\poQ\mbox{ is ccc}}$ and 
    $\genG$ is 
    a $(\uniV,\poP)$-generic filter, then $\uniV[\genG]\modelof{\poQ\mbox{ is ccc}}$. 
    \smallskip

    \wassert{2}
    \fi}
  Suppose that $\uniV\models\MA(\aleph_1)$. If
  $\uniV\modelof{\poP\mbox{ and }\poQ\mbox{ are ccc}}$ and $\genG$ is 
  a $(\uniV,\poP)$-generic filter, then $\uniV[\genG]\modelof{\poQ\mbox{ is ccc}}$. 
  \ifextended\else\qed\fi
\end{Lemma}
{\ifextended\extendedcolor 
\prf \assertof{1}: Suppose otherwise and there is a $f\in\uniV[\genG]$ \st\
$V[\genG]\modelof{\mapping{f}{\omega_1}{\poQ}
  \xmbox{ is \st, }f(\alpha), f(\beta)\xmbox{ for distinct }\alpha,\beta<\omega_1
  \xmbox{ are pairwise incompatible}}$.

Let $\utilde{f}$ be a $\poP$-name of $f$ and let $\condp_\alpha\in\poP$,
$\condq_\alpha\in\poQ$  
for $\alpha<\omega_1$ are 
\st\ $\condp_\alpha\forces{\poP}{\utilde{f}(\alpha)=\check{\condq}_0}$. 
Then $\setof{\pairof{\condp_\alpha,\condq_\alpha}}{\alpha<\omega_1}$ is pairwise 
incompatible in $\poP\times\poQ$. This is a contradiction to the assumption. \smallskip

\assertof{2}: By \assertof{1} and \CorAof{p-bedrock-19}.
\qedofLemma
\fi}

\begin{Thm}\Label{p-bedrock-21}
  Suppose that $\kappa$ is tightly ccc-Laver-gen.\ hyperhuge. Then \tfae:\smallskip

  \wassert{a} $(ccc, \calH(\kappa))$-\RcA\ holds.\smallskip

  \wassert{b} $(ccc, \calH(\kappa))$-\RcAp\ holds.\smallskip

  \wassert{c} ${V_\kappa}^{\overline{\uniW}}\prec\overline{\uniW}$ where
  $\overline{\uniW}$  
  is the bedrock of\/ $\uniV$.
\end{Thm}
\prf Note that $\kappa=\continuum$ by Theorem 5.8 in \cite{sfetal-II}. Thus, by 
\Propof{p-bedrock-14}, $\MA$ holds.

\assertof{b} $\Rightarrow$ \assertof{a}: is trivial. \assertof{a} $\Rightarrow$ 
\assertof{c}: By \Propof{p-bedrock-13},\,\assertof{1}.\smallskip

\assertof{c} $\Rightarrow$ \assertof{b}: Suppose (in $\uniV$) that \ixitem[x-bedrock-31]
$\forces{\poP}{\varphi(\ol{a})}$ for a ccc \po\ $\poP$, $\Lin$-formula
$\varphi=\varphi(\ol{x})$, and $\ol{a}\in\calH(\kappa)$.

By \Lemmaof{p-bedrock-15} there is $\poP^*\in\ol{\uniW}$ with
$\ol{\uniW}\models\cardof{\poP^*}<\kappa$ 
and $(\ol{\uniW},\poP^*)$-generic $\genG^*\in\uniV$ \st\ \ixitem[x-bedrock-32]
$\ol{a}\in\ol{\uniW}[\genG^*]$, \ixitemx[x-bedrock-33]{\qquad\qquad} $\ol{\uniW}[\genG^*]\models\MA$, and 
\ixitemx[x-bedrock-34]{\mbox{}\\} $\ol{\uniW}$ is a ccc-ground of $\uniV$. 

By \xitemof{x-bedrock-31} and \xitemof{x-bedrock-34}, we have
\begin{xitemize}
  \item[] 
  $\ol{\uniW}[\genG^*]\modelof{\mbox{there is a ccc \po\ }P
    \mbox{ \st\ }\forces{P}{\varphi(\ol{a})}}$. 
\end{xitemize}

By the assumption \assertof{c}, e have \ixitem[x-bedrock-35]
${V_\kappa}^{\ol{\uniW}[\genG^*]}\prec\ol{\uniW}[\genG^*]$. Hence
\begin{xitemize}
  \xitem[x-bedrock-36] 
  ${V_\kappa}^{\ol{\uniW}[\genG^*]}\modelof{\mbox{there is a ccc \po\ }P
    \mbox{ \st\ }\forces{P}{\varphi(\ol{a})}}$. 
\end{xitemize}
Let $\poQ^*$ be a witness of \xitemof{x-bedrock-36}. Then
$\ol{\uniW}[\genG^*]\models\cardof{\poQ^*}<\kappa$ and
$\ol{\uniW}[\genG^*]\modelof{\poQ^*\mbox{ is ccc}}$ by \xitemof{x-bedrock-35}. By 
\xitemof{x-bedrock-34} and \Lemmaof{p-bedrock-20}, we have
$\uniV\modelof{\poQ^*\mbox{ is ccc}}$.
Since $\psof{\poQ^*}^{\ol{\uniW}[\genG^*]}\in {V_\kappa}^{\ol{\uniW}[\genG^*]}$ and
$\kappa=\continuum$ in $\uniV$, \MA\ (in $\uniV$) implies that there is a
$(\ol{\uniW},\poQ^*)$-generic $\genH^*\in\uniV$.

$\ol{\uniW}[\genG^*][\genH^*]\models\varphi(\ol{a})$ since $\poQ^*$ is a witness of 
\xitemof{x-bedrock-36}, and $\ol{\uniW}[\genG^*][\genH^*]$ is a ccc-ground of $\uniV$ by 
\Lemmaof{p-bedrock-20}.
\qedofThm

\section{The Laver-Generic Maximum}
\Label{LGM} 
\memox{Janos paper の Toward the Laver-Generic Maximum を rewrite したものをここに加える．
  では ... という Axioms の組合せがほどんどすべての知られている Axioms や Principles を
  統合するものになっている，ということを \sectionof{bedrock-Lg} で指摘する．}
Suppose that $\kappa_0<\kappa_1<\kappa^*$ are regular cardinals \st\
$V_{\kappa^*}\models\ZFC$ and $\kappa_0$, $\kappa_1$ are super-$C^{(\infty)}$-hyperhuge 
cardinals in $V_{\kappa^*}$ . Note that existence of a 2-huge cardinal implies the 
consistency of this constellation (see \Lemmaof{p-Lg-RcA-2}).

Consider the following construction: 
\begin{xitemize}
  \item[] 
  first make $\kappa_0$ a tightly super-$C^{(\infty)}$-all 
  \pos-Laver gen.\ hyperhuge by a \po\ of size $\kappa_0$ as in 
  \Thmof{p-Lg-RcA-4},\,\assertof{4}. Then we force $\kappa_1$ to be a tightly super-
  $C^{(\infty)}$-semi-proper-Laver gen.\ hyperhuge by a \po\ of size $\kappa_1$ as in 
  \Thmof{p-Lg-RcA-4},\,\assertof{2'}.
\end{xitemize}
The resulting model satisfies:
\begin{xitemize}
  \xitem[x-bedrock-37] \ZFC\\
  $+$ ``\parbox[t]{0.83\textwidth}{\,$\omega_2=\continuum$ is the tightly 
    super-$C^{(\infty)}$-semi-proper-Laver gen.\  
  hyperhuge cardinal''}\bigskip\\
  $+$ ``\parbox[t]{0.83\textwidth}{\,There  is a semi-proper ground $\uniW$ of the universe
    $\uniV$ \st\ $(\continuum)^\uniW={\omega_1}={\omega_1}^\uniV$ is the tightly 
  super-$C^{(\infty)}$-all-\pos-Laver gen.\ hyperhuge cardinal in $\uniW$''.}
\end{xitemize}

We want to call \xitemof{x-bedrock-37} Laver Generic Maximum (\LGM), of cause not because of 
the maximality of possible  
consistency strength among similar assertions (this is not true since we can also switch in 
some other  
notion of large cardinal stronger than the hyperhugeness) but rather because this 
combination of the properties implies that practically all set-theoretic assertions known 
to be consistent with \ZFC\ are realized either as consequences of \xitemof{x-bedrock-37} 
or as theorems in (many of) grounds of $\uniV$ or some other inner models of $\uniV$.

So we have under the \LGM\ \xitemof{x-bedrock-37} that
\begin{xitemize}
\item[]\hspace{-2em}{\bf --}\ \ The bedrock $\ol{\uniW}$ exists\hfill (\Thmof{p-bedrock-1}).
\item[]\hspace{-2em}{\bf --}\ \ ${\omega_1}^\uniV$ and ${\omega_2}^\uniV$ are 
  super-$C^{(\infty)}$-hyperhuge cardinals in $\ol{\uniW}$
  \hfill (\Thmof{p-bedrock-2-0}). 
\item[]\hspace{-2em}{\bf --}\ \ ${V_{(\omega_1)^\uniV}}^{\ol{\uniW}}\prec\ol{\uniW}$. 
  \hfill (\Thmof{p-Lg-RcA-1-a}). 
\item[]\hspace{-2em}{\bf--}\ \ ${V_{(\omega_2)^\uniV}}^{\ol{\uniW}}\prec\ol{\uniW}$. 
  \hfill (\Thmof{p-Lg-RcA-1-a}). 
\item[]\hspace{-2em}{\bf--}\ \ $(\mbox{semi-proper}, \calH(\continuum))$-\RcAp\ holds 
  \hfill (\Thmof{p-Lg-RcA-5}). 
\item[]\hspace{-2em}{\bf--}\ \ $(\mbox{all \pos}, \calH(\aleph_1)^{\ol{\uniW}})$-\RcAp\ 
  holds \qquad (\Thmof{p-Lg-RcA-5},  
  and since 
  \RcAp\ for all\\
  \hfill\pos\ with the same parameters is preserved by generic extensions).
\item[]\hspace{-2em}{\bf--}\ \ For each natural number $n$ there are stationarily many 
  super-$C^{(n)}$-hyperhuge 
  cardinals\hfill (this holds in $\ol{\uniW}$ by \Corof{p-Lg-RcA-1-a-0}, and since 
  \\
  \hfill the statement is preserved by set forcing, it also holds in $\uniV$).
\item[]\hspace{-2em}{\bf--}\ \ For any natural number $n$ and any $a\in\calH(\continuum)$, 
  there is a semi-proper ground $\uniW$ with $a\in\uniW$ \st\ 
  $\uniW\modelof{\continuum\mbox{ is the tightly super-}C^{(n)}\xmbox{-ccc-Laver gen.\ hyperhuge
      cardinal}}$. \hfill(there is a super- $C^{(n)}$-hyperhuge cardinal\\
  \hfill\parbox[t]{0.8\textwidth}{by the previous 
  item and hence we can force the statement by a ccc forcing (by a variation of 
  \Thmof{p-Lg-RcA-4},\,\assertof{3}). 
  $(\mbox{semi-proper},\calH(\continuum))$-\RcAp\ now implies the existence of a 
  semi-proper ground 
  $\uniW$ with $a\in\uniW$ satisfying the statement). }
  \\
\item[]\hspace{-2em}{\bf--}\ \ Unbounded Resurrection Axiom of Tsaprounis in 
  \cite{tsaprounis1}  
  for semi-proper\\
  \hfill \mbox{(see \cite{future}).}
\item[]\hspace{-2em}{\bf--}\ \ $\MM^{++}$\hfill (by Theorem 5.7 in \cite{sfetal-II}). 
\item[]\hspace{-2em}{\bf--}\ \ $\continuum=\omega_2$\hfill (either by $\MM^{++}$ or by 
  \Thmof{p-Lg-RcA-1},\,\assertof{4}). 
\item[]\hspace{-2em}{\bf--}\ \ For any known instance {\sf CM} of Cichoń's Maximum (even one of 
  those in  which some mild large cardinals are involved) and any
  $a\in\calH(\aleph_2)$ 
  there is a semi-proper ground $\uniW$ with $a\in\uniW$ \st\ $\uniW\models{\sf CM}$\\
  \hfill(by
  $(\mbox{semi-proper}, \calH(\continuum))$-\RcAp).  
\item[]\hspace{-2em}{\bf--}\ \ $\cdots$
\end{xitemize}

If we can accept the tightly super-$C^{(\infty)}$-semi-proper-Laver gen.\ hyperhuge 
continuum as a natural and/or  
even desirable set-theoretic assumption, \xitemof{x-bedrock-37} may be considered as a sort 
of the final solution to the continuum problem (and actually much more) in terms of the 
properties listed above --- for discussions about arguments supporting naturalness of 
the tightly super-$C^{(\infty)}$-semi-proper-Laver gen.\ hyperhuge continuum, see also 
section 2 of \cite{janos}.  

On the other hand, it also seems that \Thmof{p-bedrock-16} and \Thmof{p-bedrock-21},  
together with the intuition that the universe of set theory should accommodate as many 
prominent grounds as possible, suggest that 
each of the combinations below are reasonable one: 

\begin{xitemize}
  \xitem[x-bedrock-37-0]  \ZFC\\
  $+$ ``\parbox[t]{0.83\textwidth}{%
    $\omega_1=\continuum$ is the tightly super-$C^{(\infty)}$-all \pos-Laver gen.\ 
    hyperhuge cardinal''.
  }
\end{xitemize}

\begin{xitemize}
  \xitem[x-bedrock-38]  \ZFC\\
  $+$ ``\parbox[t]{0.83\textwidth}{\,$\continuum$ is the tightly super-$C^{(\infty)}$-ccc 
    Laver gen.\  
    hyperhuge cardinal''}\bigskip\\
  $+$ ``\parbox[t]{0.83\textwidth}{\,There is a ccc-ground $\uniW$ of the universe $\uniV$ 
    \st\ $(\continuum)^\uniW={\omega_2}^\uniW$ is the tightly semi-proper-Laver gen.\ 
    hyperhuge cardinal in $\uniW$''}\bigskip\\
  $+$ ``\parbox[t]{0.83\textwidth}{\,There  is a semi-proper-ground $\uniW'$ of the universe $\uniV$ 
    \st\ $(\continuum)^{\uniW'}={\omega_1}={\omega_1}^\uniV$ is the tightly 
  super-$C^{(\infty)}$-all-\pos-Laver gen.\ hyperhuge cardinal in $\uniW'$\,''.}
\end{xitemize}
The combinations of axioms \xitemof{x-bedrock-37-0} and \xitemof{x-bedrock-38} will be also examined 
further in the subsequent papers.  
Note that \xitemof{x-bedrock-38} implies that the continuum is extremely large (weakly 
Mahlo and much more,  
see \cite{sfetal-II}, \cite{fuchino-sakai-2}), and that the Fodor-type Reflection Principle 
(\FRP) holds (\FRP\ holds in the ccc-ground $\uniW$ with the tightly semi-proper-Laver gen.\ 
hyperhuge continuum, since this implies $\MM^{++}$ in $\uniW$ (by Theorem 5.7 in 
\cite{sfetal-II}). Since \FRP\ is preserved by ccc generic extension (Theorem 3.4 in 
\cite{fjetal}), 
$\uniV$ also satisfies it). 
\bigskip

\noindent{\bf Added on June 20, 2025:}  After this paper has be submitted the first author 
noticed that Laver-generic versions of extendibility and super-$C^{(\infty)}$-extendibility 
can be considered and ``ultrahuge'' in Theorem 3.1 and Theorem 4.10 are replaced 
by ``extendible'' with almost the same proofs. This is an enormous improvement since the 
consistency of the existence of a tightly $\calP$-Laver-gen.\ extendible cardinal for
$\Sigma_2$-definable class $\calP$ of \pos\  with a reasonable notion of transfinite 
iteration follows from an extendible cardinal, while the consistency of a tightly super-
$C^{(\infty)}$- $\calP$-Laver-gen.\ extendible cardinal for class $\calP$ of \pos\ with a 
reasonable notion of transfinite iteration follows from an almost huge cardinal \cite{future2}.

\end{document}